\documentclass[thmsa,a4paper,oneside,final,notitlepage,onecolumn,german]{article}
\usepackage{amssymb}
\usepackage{epsfig}

\usepackage{sw20jart}

\evensidemargin\oddsidemargin
\input{tcilatex}

\begin{document} 
% THE ZIB PAGE
%\ZTPAuthor{Dimitrios I. Dais, Utz-Uwe Haus and Martin Henk}
%\ZTPTitle{On crepant resolutions of 2-parameter series of Gorenstein cyclic quotient singularities}
%\ZTPPreprint
%\ZTPNumber{98-12}
%\ZTPMonth{March}
%\ZTPYear{1998}
%---------------------------------------------------------------------------
\author{\textbf{Dimitrios I. Dais, Utz-Uwe Haus and Martin Henk}}
\title{\textbf{On crepant resolutions of 2-parameter series }\\
\textbf{of Gorenstein cyclic quotient singularities }}
\date{}
%%%%
%\zibtitlepage
%\newpage\hphantom{space}\vphantom{space}\newpage \setcounter{page}{1}
%%%

\maketitle

\begin{abstract}
\noindent An immediate generalization of the classical McKay correspondence
for Gorenstein quotient spaces $\Bbb{C}^{r}/G$ in dimensions $r\geq 4$ would
primarily demand the existence of projective, crepant, full
desingularizations. Since this is not always possible, it is natural to ask
about special classes of such quotient spaces which would satisfy the above
property. In this paper we give explicit necessary and sufficient conditions
under which 2-parameter series of Gorenstein cyclic quotient singularities
have torus-equivariant resolutions of this specific sort in all dimensions.
\end{abstract}

\section{Introduction\label{INTRO}}

\noindent Let $Y$ be a Calabi-Yau threefold and $G$ a finite group of
analytic automorphisms of $Y$, such that for all isotropy groups $G_{y}$, $%
y\in Y$, we have $G_{y}\subset $ SL$\left( \mathsf{T}_{y}^{\text{h}}\left(
Y\right) \right) $, where $\mathsf{T}_{y}^{\text{h}}\left( Y\right) $
denotes the holomorphic tangent space of $Y$ at $y$. In the framework of the
study of the ``index'' of the physical orbifold theory \cite{DHVW} Dixon,
Harvey, Vafa and Witten introduced for the orbit-space $Y/G$ a ``stringy''
analogue $\chi _{\text{str}}\left( Y,G\right) $ of the Euler-Poincar\'{e}
characteristic. Working out some concrete examples, they verified that $\chi
_{\text{str}}\left( Y,G\right) $ is equal to the usual Euler-Poincar\'{e}
characteristic $\chi \left( W\right) $ of a new Calabi-Yau threefold $W$,
which is nothing but the overlying space of a projective desingularization $%
f:W\rightarrow Y/G$ of $Y/G$, as long as $f$ does not affect the triviality
of the dualizing sheaf (in other words, $f$ preserves the ``holomorphic
volume form'' or is ``crepant''). Historically, this was the starting-point
for several mathematical investigations about crepant resolutions of
Gorenstein quotient singularities, because, in its local version, the above
property can be seen to be implied by a generalization of the so-called
``McKay-correspondence'' in dimension $r=3$ (cf. Hirzebruch-H\"{o}fer \cite
{HH}). Meanwhile there are lots of proposed very promising approaches to the
problem of generalizing McKay's $2$-dimensional bijections for arbitrary $%
r\geq 3$ (cf. Reid \cite{Reid4}), though the list of closely related open
questions seems to remain rather long.\medskip 

\noindent $\bullet $ In dimension $2$ the classification of quotient
singularities $\left( \Bbb{C}^{2}/G,\left[ \mathbf{0}\right] \right) $ (for $%
G$ a finite subgroup of SL$\left( 2,\Bbb{C}\right) $), as well as the
minimal resolution of $\left[ \mathbf{0}\right] $ by a tree of $\Bbb{P}_{%
\Bbb{C}}^{1}$'s is well-known (Klein \cite{Klein1}, Du Val \cite{DuVal}); in
fact, these rational, Gorenstein singularities can be alternatively
characterized as the $A$-$D$-$E$ hypersurface singularities being embedded
in $\Bbb{C}^{3}$ (see e.g. Lamotke \cite{Lamotke}). Roughly formulated, the 
\textit{classical} $2$-dimensional McKay-correspondence for the quotient
space $X=\Bbb{C}^{2}/G$ provides a bijection between the set of irreducible
representations of $G$ and a basis of $H^{\ast }\left( \widehat{X},\Bbb{Z}%
\right) $ (or dually, between the set of conjugacy classes of $G$ and a
basis of $H_{\ast }\left( \widehat{X},\Bbb{Z}\right) $), where $f:$ $%
\widehat{X}\rightarrow X$ $\ $is the minimal ($=$ crepant)\
desingularization of $X$. (See McKay \cite{McK}, Gonzalez-Sprinberg, Verdier 
\cite{GSV}, Kn\"{o}rrer \cite{Kn} and Ito-Nakamura \cite{Ito-Nak}).\medskip

\noindent$\bullet $ The most important aspects of the \textit{%
three-dimensional} generalization of McKay's bijections for $\Bbb{C}^{3}/G$%
's, $G\subset $ SL$\left( 3,\Bbb{C}\right) $, were only recently clarified
by the paper \cite{Ito-Reid} of Ito and Reid; considering a canonical
grading on the Tate-twist of the acting $G$'s by the so-called ``ages'',
they proved that for any projective crepant resolution $f:$ $\widehat{X}%
\rightarrow X=\Bbb{C}^{3}/G$, there are one-to-one correspondences between
the elements of $G$ of age $1$ and the exceptional prime divisors of $f$,
and between them and the members of a basis of $H^{2}\left( \widehat{X},\Bbb{%
Q}\right) $, respectively. On the other hand, the \textit{existence }of
crepant $f$'s was proved by Markushevich, Ito and Roan:

\begin{theorem}[Existence-Theorem in Dimension 3]
\label{MARIR}The underlying spaces of all $3$-dimensional Gorenstein
quotient singularities possess crepant resolutions.
\end{theorem}

\noindent(The $f$'s are unique only up to ``isomorphism in codimension $1$%
'', and to win projectivity, one has to make particular choices). From the
point of view of birational geometry, Ito and Reid \cite{Ito-Reid} proved,
in addition, the following theorem working in all dimensions.

\begin{theorem}[Ito-Reid Correspondence]
Let $G$ be a finite subgroup of \emph{SL}$\left( r,\Bbb{C}\right) $ acting
linearly on $\Bbb{C}^{r},\,r\geq 2$, and $X=\Bbb{C}^{r}/G$. Then there is a
canonical one-to-one correspondence between the junior conjugacy classes in $%
G$ and the crepant discrete valuations of $X$.
\end{theorem}

\noindent In dimensions $r\geq 4$, however, there are already from the very
beginning certain qualitative obstructions to generalize thm. \ref{MARIR},
and Reid's question (\cite{Reid3}, \cite{Ito-Reid}, \S\ 4.5, \cite{Reid4},
5.4) still remains unanswered:\medskip

\noindent \underline{$\bullet $ \textbf{Main question}}\textbf{\ : Under
which conditions on the acting groups} $G\subset \mathbf{SL}\left( r,\Bbb{C}%
\right) ,\ r\geq 4$\textbf{, do the quotient spaces }$\Bbb{C}^{r}/G$ \textbf{%
have projective crepant desingularizations?\medskip }

\noindent$\bullet $ An immediate generalization of the classical McKay
correspondence for $r\geq 4$ cannot avoid the use of a satisfactory answer
to the above question! The reason is simple. Using, for instance, only 
\textit{partial} crepant desingularizations, we obtain overlying spaces of
uncontrollable cohomology dimensions even in trivial examples. (In contrast
to this, Batyrev \cite{Bat} proved recently the
invariance of \textit{all} cohomology dimensions (over $\Bbb{Q}$) of the
overlying spaces of all ``\textit{full}'' crepant desingularizations for
arbitrary $r\geq 2$. For the case of abelian acting groups, see \cite{BD},
5.4)\medskip .

\noindent$\bullet $ It is worth mentioning that the existence of \textit{%
terminal} Gorenstein singularities implies automatically that \textit{not all%
} Gorenstein quotient spaces $\Bbb{C}^{r}/G$, $r\geq 4$, can have such
desingularizations (see Morrison-Stevens \cite{Mo-Ste}).\medskip

\noindent $\bullet $ On the other hand, as it was proved in \cite{DHZ1} by
making use of Watanabe's classification of all abelian quotient
singularities $\left( \Bbb{C}^{r}/G,\left[ \mathbf{0}\right] \right) $, $%
G\subset $ SL$\left( r,\Bbb{C}\right) $, (up to analytic isomorphism) whose
underlying spaces are embeddable as complete intersections (``c.i.'s'') of
hypersurfaces into an affine complex space, and methods of toric and
discrete geometry,

\begin{theorem}
\label{WATI}The underlying spaces of all abelian quotient c.i.-singularities
admit torus-equivariant projective, crepant resolutions \emph{(}and
therefore \emph{smooth minimal models) }in all dimensions.
\end{theorem}

\noindent Hence, the expected answer(s) to the main question will be surely
of special nature, depending crucially on the generators of the acting
groups or at least on the properties of the ring $\Bbb{C}\left[ x_{1},\ldots
,x_{r}\right] ^{G}$ of invariants.\medskip\ 

\noindent$\bullet $ In view of theorem \ref{WATI}, it is natural to ask
which would be the behaviour of abelian Gorenstein \textit{non-c.i.}
quotient singularities with respect to existence-problem of these specific
resolutions. In \cite{DH} we studied the Gorenstein $1$-parameter cyclic
quotient singularities, i.e., those of type 
\[
\frac{1}{l}\,\left( \stackunder{\left( r-1\right) \text{-times}}{\underbrace{%
1,\ldots ,1},}l-\left( r-1\right) \right) ,\ \ \ \ \ \ \ l=\left| G\right|
\geq r\geq 4\ , 
\]
generalizing the classical example of the affine cone which lies over the $r$%
-tuple Veronese embedding of $\Bbb{P}_{\Bbb{C}}^{r-1}$. Using toric
geometry, this means that the (pure) junior elements lie on a \textit{%
straight line}. In the present paper we shall extend our results in the case
in which also $2$ ``free'' parameters are allowed. Exploiting the \textit{%
coplanarity} of the corresponding junior lattice points, it is possible to
give again a definitive answer to the above main question by an explicit
arithmetical criterion involving only the weights of the defining types. As
expected from the algorithmic point of view, this criterion can be regarded
as a somewhat more complicated ``Hirzebruch-Jung-procedure'' working in all
dimensions. It should be finally pointed out, that this kind of criteria
seem to be again interesting for mathematical physics, this time in the
framework of the theory of ``D-branes'' (see Mohri \cite{Mohri} and example 
\ref{Mohris} below).\bigskip\ \newline
$\bullet $ The paper is organized as follows: In \S \ref{Prel} we recall
some basic concepts from toric geometry and fix our notation. A detailed
study of $2$-dimensional rational s.c.p.cones and a method for the
determination of the vertices of the corresponding support polygons by
Kleinian approximations are presented in \S \ref{BRUCH}. In section \ref{USP}
we explain how the underlying spaces of abelian quotient singularities are
to be regarded as affine toric varieties and recall the geometric condition
under which the quotient spaces $\Bbb{C}^{r}/G$ are Gorenstein. The first
part of section \ref{LATR} is devoted to the hierarchy of lattice
triangulations and to the simple combinatorial mechanism leading to the
hierarchy of partial (resp. full) crepant desingularizations of $\Bbb{C}%
^{r}/G$. Our main theorems on the existence of projective, crepant, full
resolutions are formulated and discussed in the second part of \S \ref{LATR}
and proved in \S \ref{MAIN}. Finally, in section \ref{COHOM} we give
concrete formulae for the computation of the dimensions of the cohomology
groups of the spaces resolving fully the $2$-parameter Gorenstein cyclic
quotient singularities by any crepant birational morphism.\bigskip

\noindent\textit{General terminology.} We always work with normal complex
varieties, i.e., with normal, integral, separated schemes over $\Bbb{C}$.
Sing$\left( X\right) $ denotes the \textit{singular locus} of such a variety 
$X$. (By the word \textit{singularity} we intimate either a singular point
or the germ of a singular point, but the meaning will be in each case clear
from the context). As in \cite{BD}, \cite{DH}, \cite{DHZ1}, by a \textit{%
desingularization} (or \textit{resolution of singularities}) $f:\widehat{X}%
\rightarrow X$ of a non-smooth $X$, we mean a ``full'' or ``overall''
desingularization (if not mentioned), i.e., Sing$\left( \widehat{X}\right)
=\varnothing $. When we refer to \textit{partial} desingularizations, we
mention it explicitly. A partial desingularization $f:X^{\prime }\rightarrow
X$ of a $\Bbb{Q}$-Gorenstein complex variety $X$ (with global index $j$) is
called\ \textit{non-discrepant} or simply \textit{crepant} if $\omega _{X}^{%
\left[ j\right] }\cong f_{\ast }\left( \omega _{X^{\prime }}^{\otimes
j}\right) $, or, in other words, if the (up to rational equivalence uniquely
determined) difference $jK_{X^{\prime }}-f^{\ast }\left( jK_{X}\right) $
contains exceptional prime divisors which have vanishing multiplicities. ($%
\omega _{X},K_{X}$ and $\omega _{X^{\prime }},K_{X^{\prime }}$ denote here
the dualizing sheaves and the canonical divisors of $X$ and $Y$
respectively). Furthermore, $f:X^{\prime }\rightarrow X$ is \textit{%
projective} if $X^{\prime }$ admits an $f$-ample Cartier divisor. The terms 
\textit{canonical} and \textit{terminal singularity} are to be understood in
the usual sense (see Reid \cite{Reid1}, \cite{Reid2}).

\section{Preliminaries from the theory of toric varieties\label{Prel}}

\noindent At first we recall some basic facts from the theory of toric
varieties. We shall mainly use the same notation as in \cite{DH}, \cite{DHZ1}%
, \cite{DHZ2}. Our standard references on toric geometry are the textbooks
of Oda \cite{Oda}, Fulton \cite{Fulton} and Ewald \cite{Ewald}, and the
lecture notes \cite{KKMS}. \medskip

\noindent \textsf{(a)} The \textit{linear hull, }the\textit{\ affine hull},
the \textit{positive hull} and \textit{the convex hull} of a set $B\subset $ 
$\Bbb{R}^{r}$, $r\geq 1,$ will be denoted by lin$\left( B\right) $, aff$%
\left( B\right) $, pos$\left( B\right) $ (or $\Bbb{R}_{\geq 0}\,B$) and conv$%
\left( B\right) $ respectively. The \textit{dimension} dim$\left( B\right) $
of a $B\subset \Bbb{R}^{r}$ is defined to be the dimension of its affine
hull. \newline
\newline
\textsf{(b) }Let $N\cong \Bbb{Z}^{r}$ be a free $\Bbb{Z}$-module of rank $%
r\geq 1$. $N$ can be regarded as a \textit{lattice }in $N_{\Bbb{R}%
}:=N\otimes _{\Bbb{Z}}\Bbb{R}\cong \Bbb{R}^{r}$. (For fixed identification,
we shall represent the elements of $N_{\Bbb{R}}$ by column-vectors in $\Bbb{R%
}^{r}$). If $\left\{ n_{1},\ldots ,n_{r}\right\} $ is a $\Bbb{Z}$-basis of $%
N $, then 
\[
\text{det}\left( N\right) :=\left| \text{det}\left( n_{1},\ldots
,n_{r}\right) \right| 
\]
is the \textit{lattice determinant}. An $n\in N$ is called \textit{primitive}
if conv$\left( \left\{ \mathbf{0},n\right\} \right) \cap N$ contains no
other points except $\mathbf{0}$ and $n$.\smallskip

Let $N\cong \Bbb{Z}^{r}$ be as above, $M:=$ Hom$_{\Bbb{Z}}\left( N,\Bbb{Z}%
\right) $ its dual lattice, $N_{\Bbb{R}},M_{\Bbb{R}}$ their real scalar
extensions, and $\left\langle .,.\right\rangle :N_{\Bbb{R}}\times M_{\Bbb{R}%
}\rightarrow \Bbb{R}$ the natural $\Bbb{R}$-bilinear pairing. A subset $%
\sigma $ of $N_{\Bbb{R}}$ is called \textit{strongly convex polyhedral cone }%
(\textit{s.c.p.cone}, for short), if there exist $n_{1},..,n_{k}\in N_{\Bbb{R%
}}$, such that $\sigma =$ pos$\left( \left\{ n_{1},\ldots ,n_{k}\right\}
\right) $ and $\sigma \cap \left( -\sigma \right) =\left\{ \mathbf{0}%
\right\} $. Its \textit{relative interior }int$\left( \sigma \right) $
(resp. its \textit{relative boundary }$\partial \sigma $) is the usual
topological interior (resp. the usual topological boundary) of it,
considered as subset of lin$\left( \sigma \right) $. The \textit{dual cone}
of $\sigma $ is defined by 
\[
\sigma ^{\vee }:=\left\{ \mathbf{x}\in M_{\Bbb{R}}\ \left| \ \left\langle 
\mathbf{x},\mathbf{y}\right\rangle \geq 0,\ \forall \mathbf{y},\ \mathbf{y}%
\in \sigma \right. \right\} 
\]
and satisfies: $\sigma ^{\vee }+\left( -\sigma ^{\vee }\right) =M_{\Bbb{R}}$
and dim$\left( \sigma ^{\vee }\right) =r$. \ A subset $\tau $ of a s.c.p.
cone $\sigma $ is called a \textit{face} of $\sigma $ (notation: $\tau \prec
\sigma $), if $\tau =\left\{ \mathbf{y}\in \sigma \ \left| \ \left\langle
m_{0},\mathbf{y}\right\rangle =0\right. \right\} $, for some $m_{0}\in
\sigma ^{\vee }$. A s.c.p. $\sigma =$ pos$\left( \left\{ n_{1},\ldots
,n_{k}\right\} \right) $ is called \textit{simplicial} (resp. \textit{%
rational}) if $n_{1},\ldots ,n_{k}$ are $\Bbb{R}$-linearly independent
(resp. if $n_{1},\ldots ,n_{k}\in N_{\Bbb{Q}}$, where $N_{\Bbb{Q}}:=N\otimes
_{\Bbb{Z}}\Bbb{Q}$).\newline
\newline
\textsf{(c) }If $\sigma \subset N_{\Bbb{R}}$ is a rational s.c.p. cone, then 
$\sigma $ has $\mathbf{0}$ as its apex and the subsemigroup $\sigma \cap N$
of $N$ is a monoid. The following two propositions describe the fundamental
properties of this monoid $\sigma \cap N$ and their proofs go essentially
back to Gordan \cite{Gordan}, Hilbert \cite{Hilbert} and van der Corput \cite
{vanderCorput1}, \cite{vanderCorput2}.

\begin{proposition}[Gordan's lemma]
\label{Gorlem}$\sigma \cap N$ is finitely generated as additive semigroup,
i.e., there exist 
\[
\ n_{1},n_{2},\ldots ,n_{\nu }\in \sigma \cap N\ ,\text{ \ such that\ \ \ }%
\sigma \cap N=\Bbb{Z}_{\geq 0}\ n_{1}+\Bbb{Z}_{\geq 0}\ n_{2}+\cdots +\Bbb{Z}%
_{\geq 0}\ n_{\nu }\ .
\]
\end{proposition}

\begin{proposition}[Minimal generating system]
\label{MINGS}Among all the systems of generators of $\sigma \cap N$, there
is a system $\mathbf{Hlb}_{N}\left( \sigma \right) $ of \emph{minimal
cardinality}, which is uniquely determined \emph{(}up to the ordering of its
elements\emph{)} by the following characterization \emph{:\smallskip } 
\begin{equation}
\mathbf{Hlb}_{N}\left( \sigma \right) =\left\{ n\in \sigma \cap \left(
N\smallsetminus \left\{ \mathbf{0}\right\} \right) \ \left| \ 
\begin{array}{c}
n\ \text{\emph{cannot be expressed as the sum }} \\ 
\text{\emph{of two other vectors belonging }} \\ 
\text{\emph{\ to\ } }\sigma \cap \left( N\smallsetminus \left\{ \mathbf{0}%
\right\} \right) 
\end{array}
\right. \right\}   \label{Hilbbasis}
\end{equation}
In particular, if $\sigma =$ \emph{pos}$\left( \left\{ n_{1},\ldots
,n_{k}\right\} \right) $, then
\begin{equation}
\mathbf{Hlb}_{N}\left( \sigma \right) \subseteq \left\{ n\in \sigma \cap
\left( N\smallsetminus \left\{ \mathbf{0}\right\} \right) \ \left| \
n=\tsum\nolimits_{i=1}^{k}\delta _{i}n_{i}\right. \text{\emph{, with\ } }%
0\leq \delta _{i}<1,\ \forall i,\ 1\leq i\leq k\right\} .  \label{HilbP}
\end{equation}
\end{proposition}

\noindent $\mathbf{Hlb}_{N}\left( \sigma \right) $ is called\emph{\ }\textit{%
the Hilbert basis of }$\sigma $ w.r.t. $N.$ About \textit{algorithms} for
the determination of Hibert bases of pointed rational cones, we refer to
Henk-Weismantel \cite{HEW}, and to the other references therein.\smallskip
\medskip 

\noindent \textsf{(d)} For a lattice $N\cong \Bbb{Z}^{r}$ having $M$ as its
dual, we define an $r$-dimensional \textit{algebraic torus }$T_{N}\cong
\left( \Bbb{C}^{\ast }\right) ^{r}$ by : 
\[
T_{N}:=\text{Hom}_{\Bbb{Z}}\left( M,\Bbb{C}^{\ast }\right) =N\otimes _{\Bbb{Z%
}}\Bbb{C}^{\ast }. 
\]
Every $m\in M$ assigns a character $\mathbf{e}\left( m\right)
:T_{N}\rightarrow \Bbb{C}^{\ast }$. Moreover, each $n\in N$ determines an $1$%
-parameter subgroup 
\[
\gamma _{n}:\Bbb{C}^{\ast }\rightarrow T_{N}\ \ \ \text{with\ \ \ }\gamma
_{n}\left( \lambda \right) \left( m\right) :=\lambda ^{\left\langle
m,n\right\rangle }\text{, \ \ for\ \ \ }\lambda \in \Bbb{C}^{\ast },\ m\in
M\ .\ \ 
\]
\ We can therefore identify $M$ with the character group of $T_{N}$ and $N$
with the group of $1$-parameter subgroups of $T_{N}$. On the other hand, for
a rational s.c.p.c. $\sigma $ with 
\[
M\cap \sigma ^{\vee }=\Bbb{Z}_{\geq 0}\ m_{1}+\Bbb{Z}_{\geq 0}\ m_{2}+\cdots
+\Bbb{Z}_{\geq 0}\ m_{k}, 
\]
we associate to the finitely generated, normal, monoidal $\Bbb{C}$%
-subalgebra $\Bbb{C}\left[ M\cap \sigma ^{\vee }\right] $ of $\Bbb{C}\left[ M%
\right] $ an affine complex variety 
\[
U_{\sigma }:=\text{Max-Spec}\left( \Bbb{C}\left[ M\cap \sigma ^{\vee }\right]
\right) , 
\]
which can be identified with the set of semigroup homomorphisms : 
\[
U_{\sigma }=\left\{ u:M\cap \sigma ^{\vee }\ \rightarrow \Bbb{C\ }\left| 
\begin{array}{c}
\ u\left( \mathbf{0}\right) =1,\ u\left( m+m^{\prime }\right) =u\left(
m\right) \cdot u\left( m^{\prime }\right) ,\smallskip \  \\ 
\text{for all \ \ }m,m^{\prime }\in M\cap \sigma ^{\vee }
\end{array}
\right. \right\} \ , 
\]
where $\mathbf{e}\left( m\right) \left( u\right) :=u\left( m\right) ,\
\forall m,\ m\in M\cap \sigma ^{\vee }\ $ and\ $\forall u,\ u\in U_{\sigma }$%
. In the analytic category, $U_{\sigma }$, identified with its image under
the injective map $\left( \mathbf{e}\left( m_{1}\right) ,\ldots ,\mathbf{e}%
\left( m_{k}\right) \right) :U_{\sigma }\hookrightarrow \Bbb{C}^{k}$, can be
regarded as an analytic set determined by a system of equations of the form:
(monomial) = (monomial). This analytic structure induced on $U_{\sigma }$ is
independent of the semigroup generators $\left\{ m_{1},\ldots ,m_{k}\right\} 
$ and each map $\mathbf{e}\left( m\right) $ on $U_{\sigma }$ is holomorphic
w.r.t. it. In particular, for $\tau \prec \sigma $, $U_{\tau }$ is an open
subset of $U_{\sigma }$. Moreover, if we define $d:=\#\left( \mathbf{Hlb}%
_{M}\left( \sigma ^{\vee }\right) \right) \ \left( \leq k\right) $, then $d$
is nothing but the \textit{embedding dimension} of $U_{\sigma }$, i.e. the 
\textit{minimal} number of generators of the maximal ideal of the local $%
\Bbb{C}$-algebra $\mathcal{O}_{U_{\sigma },\ \left( \mathbf{0}\in \Bbb{C}%
^{d}\right) }$ (cf. \cite{Oda}, 1.2 \& 1.3).\bigskip

\noindent \textsf{(e)} A \textit{fan }w.r.t.\textit{\ }$N\cong \Bbb{Z}^{r}$
is a finite collection $\Delta $ of rational s.c.p. cones in $N_{\Bbb{R}}$,
such that :\smallskip \newline
(i) any face $\tau $ of $\sigma \in \Delta $ belongs to $\Delta $,
and\smallskip \newline
(ii) for $\sigma _{1},\sigma _{2}\in \Delta $, the intersection $\sigma
_{1}\cap \sigma _{2}$ is a face of both $\sigma _{1}$ and $\sigma
_{2}.\smallskip $\newline
The union $\left| \Delta \right| :=\cup \left\{ \sigma \ \left| \ \sigma \in
\Delta \right. \right\} $ is called the \textit{support }of $\Delta $.
Furthermore, we define 
\[
\Delta \left( i\right) :=\left\{ \sigma \in \Delta \ \left| \ \text{dim}%
\left( \sigma \right) =i\right. \right\} \ ,\ \text{for \ \ }0\leq i\leq r\
. 
\]
If $\varrho \in \Delta \left( 1\right) $, then there exists a unique
primitive vector $n\left( \varrho \right) \in N\cap \varrho $ with $\varrho =%
\Bbb{R}_{\geq 0}\ n\left( \varrho \right) $ and each cone $\sigma \in \Delta 
$ can be therefore written as 
\[
\sigma =\sum\Sb \varrho \in \Delta \left( 1\right)  \\ \varrho \prec \sigma 
\endSb \ \Bbb{R}_{\geq 0}\ n\left( \varrho \right) \ \ . 
\]
The set Gen$\left( \sigma \right) :=\left\{ n\left( \varrho \right) \ \left|
\ \varrho \in \Delta \left( 1\right) ,\varrho \prec \sigma \right. \right\} $
is called the\textit{\ set of minimal generators }(within the pure first
skeleton) of $\sigma $. For $\Delta $ itself one defines analogously 
\[
\text{Gen}\left( \Delta \right) :=\bigcup_{\sigma \in \Delta }\text{ Gen}%
\left( \sigma \right) \ .\medskip 
\]
\textsf{(f) }The \textit{toric variety X}$\left( N,\Delta \right) $
associated to a fan\textit{\ }$\Delta $ w.r.t. the lattice\textit{\ }$N$ is
by definition the identification space 
\[
X\left( N,\Delta \right) :=\left( \left( \bigcup_{\sigma \in \Delta }\
U_{\sigma }\right) \ /\ \sim \right) 
\]
with $U_{\sigma _{1}}\ni u_{1}\sim u_{2}\in U_{\sigma _{2}}$ if and only if
there is a $\tau \in \Delta ,$ such that $\tau \prec \sigma _{1}\cap \sigma
_{2}$ and $u_{1}=u_{2}$ within $U_{\tau }$. As complex variety, $X\left(
N,\Delta \right) $ turns out to be irreducible, normal, Cohen-Macaulay and
to have at most rational singularities (cf. \cite{Fulton}, p. 76, and \cite
{Oda}, thm. 1.4, p. 7, and cor. 3.9, p. 125). $X\left( N,\Delta \right) $ is
called \textit{simplicial }if all cones of $\Delta $ are simplicial.\medskip 
\newline
\textsf{(g)} $X\left( N,\Delta \right) $ admits a canonical $T_{N}$-action
which extends the group multiplication of $T_{N}=U_{\left\{ \mathbf{0}%
\right\} }$ : 
\begin{equation}
T_{N}\times X\left( N,\Delta \right) \ni \left( t,u\right) \longmapsto
t\cdot u\in X\left( N,\Delta \right)  \label{torus action}
\end{equation}
where, for $u\in U_{\sigma }$, $\left( t\cdot u\right) \left( m\right)
:=t\left( m\right) \cdot u\left( m\right) ,\ \forall m,\ m\in M\cap \sigma
^{\vee }$ . The orbits w.r.t. the action (\ref{torus action}) are
parametrized by the set of all the cones belonging to $\Delta $. For a $\tau
\in \Delta $, we denote by orb$\left( \tau \right) $ (resp. by $V\left( \tau
\right) $) the orbit (resp. the closure of the orbit) which is associated to 
$\tau $. If $\tau \in \Delta $, then $V\left( \tau \right) =X\left( N\left(
\tau \right) ,\text{ Star}\left( \tau ;\Delta \right) \right) $ is itself a
toric variety w.r.t. 
\[
N\left( \tau \right) :=N/N_{\tau }\ ,\ \ \ \ \ \text{Star}\left( \tau
;\Delta \right) :=\left\{ \overline{\sigma }\ \left| \ \sigma \in \Delta ,\
\tau \prec \sigma \right. \right\} \ , 
\]
where $N_{\tau }$ is the sublattice of $N$ generated (as subgroup) by $N\cap 
$lin$\left( \Bbb{\tau }\right) $ and $\overline{\sigma }\ =\left( \sigma
+\left( N_{\tau }\right) _{\Bbb{R}}\right) /\left( N_{\tau }\right) _{\Bbb{R}%
}$ denotes the image of $\sigma $ in $N\left( \tau \right) _{\Bbb{R}}=N_{%
\Bbb{R}}/\left( N_{\tau }\right) _{\Bbb{R}}$.\medskip

\noindent \textsf{(h)}\textit{\ }A \textit{map of fans\ }$\varpi :\left(
N^{\prime },\Delta ^{\prime }\right) \rightarrow \left( N,\Delta \right) $
is a $\Bbb{Z}$-linear homomorphism $\varpi :N^{\prime }\rightarrow N$ whose
scalar extension $\varpi =\varpi _{\Bbb{R}}:N_{\Bbb{R}}^{\prime }\rightarrow
N_{\Bbb{R}}$ satisfies the property: 
\[
\forall \sigma ^{\prime },\ \sigma ^{\prime }\in \Delta ^{\prime }\ \ \text{ 
}\exists \ \sigma ,\ \sigma \in \Delta \ \ \text{ with\ \ }\varpi \left(
\sigma ^{\prime }\right) \subset \sigma \,. 
\]
$\varpi \otimes _{\Bbb{Z}}$id$_{\Bbb{C}^{\ast }}:T_{N^{\prime }}=N^{\prime
}\otimes _{\Bbb{Z}}\Bbb{C}^{\ast }\rightarrow T_{N}=N\otimes _{\Bbb{Z}}\Bbb{C%
}^{\ast }$ is a homomorphism from $T_{N^{\prime }}$ to $T_{N}$ and the
scalar extension $\varpi ^{\vee }:M_{\Bbb{R}}\rightarrow M_{\Bbb{R}}^{\prime
}$ of the dual $\Bbb{Z}$-linear map $\varpi ^{\vee }:M\rightarrow M^{\prime
} $ induces canonically an \textit{equivariant holomorphic map }$\varpi
_{\ast }:X\left( N^{\prime },\Delta ^{\prime }\right) \rightarrow X\left(
N,\Delta \right) $. This map is\textit{\ proper} if and only if $\varpi
^{-1}\left( \left| \Delta \right| \right) =\left| \Delta ^{\prime }\right| .$
In particular, if $N=N^{\prime }$ and $\Delta ^{\prime }$ is a refinement of 
$\Delta $, then id$_{\ast }:X\left( N,\Delta ^{\prime }\right) \rightarrow
X\left( N,\Delta \right) $ is \textit{proper}\emph{\ }and \textit{birational 
}(cf. \cite{Oda}, thm. 1.15 and cor. 1.18). \medskip

\noindent \textsf{(i)}\textit{\ }Let $N\cong \Bbb{Z}^{r}$ be a lattice of
rank $r$ and $\sigma \subset N_{\Bbb{R}}$ a simplicial, rational s.c.p.c. of
dimension $k\leq r$. $\sigma $ can be obviously written as $\sigma =\varrho
_{1}+\cdots +\varrho _{k}$, for distinct $1$-dimensional cones $\varrho
_{1},\ldots ,\varrho _{k}$. We denote by 
\[
\mathbf{Par}\left( \sigma \right) :=\left\{ \mathbf{y}\in \left( N_{\sigma
}\right) _{\Bbb{R}}\ \left| \ \mathbf{y}=\sum_{j=1}^{k}\ \varepsilon _{j}\
n\left( \varrho _{j}\right) ,\ \text{with\ \ }0\leq \varepsilon _{j}<1,\
\forall j,\ 1\leq j\leq k\right. \right\} 
\]
the \textit{fundamental }(\textit{half-open})\textit{\ parallelotope }which
is associated to\textit{\ }$\sigma $. The \textit{multiplicity} mult$\left(
\sigma ;N\right) $ of $\sigma $ with respect to $N$ is defined as 
\[
\text{mult}\left( \sigma ;N\right) :=\#\left( \mathbf{Par}\left( \sigma
\right) \cap N_{\sigma }\right) =\text{Vol}\left( \mathbf{Par}\left( \sigma
\right) ;N_{\sigma }\right) \ , 
\]
where Vol$\left( \mathbf{Par}\left( \sigma \right) \right) $ denotes the
usual volume (Lebesgue measure) of $\mathbf{Par}\left( \sigma \right) $ and 
\[
\text{Vol}\left( \mathbf{Par}\left( \sigma \right) ;N_{\sigma }\right) :=%
\frac{\text{Vol}\left( \mathbf{Par}\left( \sigma \right) \right) }{\text{det}%
\left( N_{\sigma }\right) }=\frac{\text{det}\left( \Bbb{Z}\,n\left( \varrho
_{1}\right) \oplus \cdots \oplus \Bbb{Z}\,n\left( \varrho _{k}\right)
\right) }{\text{det}\left( N_{\sigma }\right) } 
\]
its the relative volume w.r.t. $N_{\sigma }$. If mult$\left( \sigma
;N\right) =1$, then $\sigma $ is called a \textit{basic cone} \textit{w.r.t.}
$N$.

\begin{proposition}[Smoothness Criterion]
\label{SMCR}The affine toric variety $U_{\sigma }$ is smooth iff $\sigma $
is basic \textit{w.r.t.} $N$. \emph{(}Correspondingly, an arbitrary toric
variety $X\left( N,\Delta \right) $ is smooth if and only if it is
simplicial and each s.c.p. cone $\sigma \in \Delta $ is basic.\emph{)}
\end{proposition}

\noindent \textit{Proof. }It follows from \cite{Oda}, thm. 1.10, p. 15. $%
_{\Box }\medskip $

\noindent By Carath\'{e}odory's theorem concerning convex polyhedral cones
(cf. Ewald \cite{Ewald}, III 2.6, p. 75, \& V 4.2, p. 158) one can choose a
refinement $\Delta ^{\prime }$ of any given fan $\Delta $, so that $\Delta
^{\prime }$ becomes simplicial. Since further subdivisions of $\Delta
^{\prime }$ reduce the multiplicities of its cones, we may arrive (after
finitely many subdivisions) at a fan $\widetilde{\Delta }$ having only basic
cones.

\begin{theorem}[Existence of Desingularizations]
For every toric variety $X\left( N,\Delta \right) $ there exists a
refinement $\widetilde{\Delta }$ of $\Delta $ consisting of exclusively
basic cones w.r.t. $N$, i.e., such that 
\[
f=\emph{id}_{\ast }:X\left( N,\widetilde{\Delta }\right) \longrightarrow
X\left( N,\Delta \right) =U_{\sigma }
\]
is a \emph{(}full\emph{) }desingularization.
\end{theorem}

\noindent Though this theorem can be also treated in terms of blow-ups of
(not necessarily reduced) subschemes of $X\left( N,\Delta \right) $
supporting Sing$\left( X\left( N,\Delta \right) \right) $ (cf. \cite{KKMS},
\S I.2), it is a theorem of purely ``existential nature'', because it does
not provide any intrinsic characterization of a ``canonical'' way for the
construction of the required subdivisions. On the other hand, there is
another theorem due to Seb\"{o} (\cite{Sebo}, thm. 2.2) which informs us
that such an intrinsically geometric choice of subdivisions is indeed
possible in \textbf{low} dimensions by considering the elements of the
Hilbert basis of a cone as the set of minimal generators of the
desingularizing fan.

\begin{theorem}[Seb\"{o}'s Theorem]
\label{Seb}Let $N$ be a lattice of rank $r$ and $\sigma \subset N_{\Bbb{R}}$
a rational, simplicial s.c.p. cone of dimension $r$. Moreover, let $\Delta $
denote the fan consisting of $\sigma $ together with all its faces. For $%
r\leq 3$ there exists a \emph{(}full\emph{) }desingularization \ \ $f=$ 
\emph{id}$_{\ast }:X\left( N,\widetilde{\Delta }\right) \longrightarrow
X\left( N,\Delta \right) $, such that 
\[
\emph{Gen}\left( \widetilde{\Delta }\right) =\mathbf{Hlb}_{N}\left( \sigma
\right) \,\ .
\]
\end{theorem}

\noindent This theorem was reproved independently by Aguzzoli \& Mundizi ( 
\cite{A-M}, prop. 2.4) and by Bouvier \& Gonzalez-Sprinberg (\cite{B-GS1},
th. 1, \cite{B-GS2}, th. 2.9) by using slightly different methods. Moreover,
the latter authors constructed two concrete counterexamples showing that the
statement in \ref{Seb} is not true (in general) for dimensions $r\geq 4$.
(See also rem. \ref{FIRLA} below).

\section{Finite continued fractions and two-dimensional rational cones\label%
{BRUCH}}

\noindent As we have already pointed out in the introduction, for the study
of the existence of torus-equivariant, crepant, full resolutions of $2$%
-parameter-series of Gorenstein cyclic quotient singularities we shall
exploit the particular property of the corresponding lattice points of the
junior simplex to be \textit{coplanar} in an essential way. Consequently we
shall reduce the whole problem to a $2$-dimensional one (after appropriate
lattice transformations), and we shall therefore make use of the precise
structure of $2$-dimensional rational s.c.p.cones. That's why we take a
closer look, in this section, at the interplay between the
``lattice-geometry'' of two-dimensional rational s.c.p. cones and the
continued fraction expansions of rational numbers realized by their
``parametrizing integers'' $p$ and $q$.\medskip

\noindent $\bullet $ It was Hirzebruch \cite{Hirz1} in the early fifties who
first described the minimal resolution of \textit{any} $2$-dimensional
cyclic quotient singularity by means of the ``modified'' euclidean division
algorithm and $\left( -\right) $-sign continued fractions, by blowing up
points and by generalizing some previous partial results of Jung \cite{Jung}%
. Later Cohn \cite{Cohn} proposed the use of support polygons with lattice
points as vertices (instead of complex coordinate charts), in order to
simplify considerably the resolution-procedure. Around the same time,
Saint-Donat (cf. \cite{KKMS}, pp. 16-19 \& 35-38) gave the first proofs
exclusively in terms of toric geometry. A much more detailed explanation of
the role of negative-regular continued fractions in the study of $2$%
-dimensional rational cones (via polar polyhedra) is contained in Oda's book 
\cite{Oda}, \S 1.6.\medskip\ \newline
$\bullet $ Though we shall mostly adopt here Oda's notation and basic
concepts (so that the proofs remain valid for arbitrary lattices of rank $2$%
), our purpose is not to say anything about the resolution itself, which is
well-known and not crucial for later use, but to give an explicit method for
the determination of the \textit{vertices }of the support polygons. We are
mainly motivated by the original works of Klein \cite{Klein2}, \cite{Klein3}%
, \cite{Klein4}, written over hundred years ago, in which his ``Umri\ss %
polygone'' (also known as \textit{Kleinian polygons}, cf. Finkel'shtein \cite
{Fink}) were used to approximate real (not necessarily rational) numbers by
only \textit{regular }continued fractions. From the algorithmic point of
view, Kleinian approximations are more ``economic'' (see remarks \ref{RELE}, 
\ref{length1} and \ref{length2} below). Technically, it seems to be more
convenient to work simultaneously with both primary and dual cones $\sigma $
and $\sigma ^{\vee }$ and to interchange their combinatorial data by
specific transfer-rules.\bigskip\ 

\noindent \textit{Notation. }We shall henceforth use the following extra
notation. For $\nu \in \Bbb{N}$, $\mu \in \Bbb{Z}$, we denote by $\left[
\,\mu \,\right] _{\nu }$ the (uniquely determined) integer for which 
\[
0\leq \left[ \,\mu \,\right] _{\nu }<\nu ,\ \ \ \,\mu \,\equiv \left[ \,\mu
\,\right] _{\nu }\left( \text{mod }\nu \right) . 
\]
``gcd'' will be abbreviation for greatest common divisor. If $q\in \Bbb{Q}$,
we define $\left\lfloor q\right\rfloor $ (resp. $\left\lceil q\right\rceil
\, $) to be the greatest integer number $\leq q$ (resp. the smallest integer 
$\geq q$).\medskip

\noindent \textsf{(a) }Let $\kappa $ and $\lambda $ be two given relatively
prime positive integers. Suppose that $\frac{\kappa }{\lambda }$ can be written
as 
\begin{equation}
\frac{\kappa }{\lambda }=a_{1}+\frac{\varepsilon _{1}}{a_{2}+\dfrac{%
\varepsilon _{2}}{a_{3}+\dfrac{\varepsilon _{3}}{%
\begin{array}{ll}
\ddots &  \\ 
& a_{\nu -1}+\dfrac{\varepsilon _{\nu -1}}{a_{\nu }}
\end{array}
}}}  \label{conf}
\end{equation}
The right-hand side of (\ref{conf}) is called \textit{semi-regular continued
fraction }for $\frac{\kappa }{\lambda }$ (and $\nu $ its \textit{length}) if
it has the following properties:\medskip \newline
(i) $a_{j}$ is an integer for all $j$, $1\leq j\leq \nu $,\smallskip \newline
(ii) $\varepsilon _{j}\in \left\{ -1,1\right\} $ for all $j$, $1\leq j\leq
\nu -1$,\smallskip \newline
(iii) $a_{j}\geq 1$ for $j\geq 2$ and $a_{\nu }\geq 2\,$(!), and\smallskip 
\newline
(iv) if $a_{j}=1$ for some $j$, $1<j<\nu $, then $\varepsilon
_{j}=1.\bigskip $\newline
In particular, if $\varepsilon _{j}=1$ (resp. $\varepsilon _{j}=-1$) for all 
$j$, $1\leq j<\nu $, \smallskip then we write $\frac{\kappa }{\lambda }=%
\left[ a_{1},a_{2},\ldots ,a_{\nu }\right] $ (resp. $\frac{\kappa }{\lambda }%
=\left[ \!\left[ a_{1},a_{2},\ldots ,a_{\nu }\right] \!\right] $). This is
the\textit{\ regular }(resp. the \textit{negative-regular}$)$\textit{\
continued fraction expansion} of $\frac{\kappa }{\lambda }$. These two
expansions are unique (in this form) and can be obtained by the usual and
the modified euclidean algorithm, respectively, depending on the choice of
the kind of the associated remainders. The next two lemmas outline their
main properties.

\begin{lemma}
\label{POS}Let $\lambda ,\kappa $ be two integers with $0<\lambda <\kappa ,$ 
\emph{gcd}$\left( \lambda ,\kappa \right) =1$. Then\emph{:\smallskip }%
\newline
\emph{(i)} There exists always a uniquely determined regular continued
fraction expansion 
\begin{equation}
\frac{\kappa }{\lambda }=\left[ a_{1},a_{2},\ldots ,a_{\nu }\right] 
\label{posconf}
\end{equation}
of $\frac{\kappa }{\lambda }$ \emph{(}with $a_{j}\geq 1$ for all $j\geq 1$
and $a_{\nu }\geq 2$\emph{).\smallskip } \newline
\emph{(ii)} Defining the finite sequences $\left( \mathsf{P}_{i}\right)
_{-1\leq i\leq \nu }$ and $\left( \mathsf{Q}_{i}\right) _{-1\leq i\leq \nu }$
by 
\begin{equation}
\left\{ 
\begin{array}{lll}
\mathsf{P}_{-1}=0, & \mathsf{P}_{0}=1, & \mathsf{P}_{i}=a_{i}\ \mathsf{P}%
_{i-1}+\mathsf{P}_{i-2},\ \forall i,\ \ 1\leq i\leq \nu , \\ 
& \  &  \\ 
\mathsf{Q}_{-1}=1, & \mathsf{Q}_{0}=0, & \mathsf{Q}_{i}=a_{i}\ \mathsf{Q}%
_{i-1}+\mathsf{Q}_{i-2},\ \forall i,\ \ 1\leq i\leq \nu ,
\end{array}
\right.   \label{PQ}
\end{equation}
we obtain for all $i$, $1\leq i\leq \nu $, 
\[
\frac{\mathsf{P}_{i}}{\mathsf{Q}_{i}}=\left\{ 
\begin{array}{ll}
\left[ a_{1},a_{2},\ldots ,a_{i}\right] , & \text{\emph{if \ }}a_{i}\geq 2%
\text{\emph{\ }} \\ 
\left[ a_{1},a_{2},\ldots ,a_{i-2},a_{i-1}+1\right] , & \text{\emph{if \ }}%
a_{i}=1
\end{array}
\right. 
\]
\emph{(This is the so-called }$i$\textit{-th convergent}\emph{\ of the given
continued fraction (\ref{posconf})).\smallskip }\newline
\emph{(iii)} The above sequences satisfy the following conditions\emph{%
:\smallskip }\newline
\emph{(a) }For all $i$, $1\leq i\leq \nu $ \emph{(resp. }$2\leq i\leq \nu $%
\emph{)}, 
\begin{equation}
\frac{\mathsf{P}_{i}}{\mathsf{P}_{i-1}}=\left\{ 
\begin{array}{ll}
\left[ a_{i},a_{i-1},\ldots ,a_{2},a_{1}\right] , & \text{\emph{if \ }}%
a_{1}\geq 2\text{\emph{\ }} \\ 
\left[ a_{i},a_{i-1},\ldots ,a_{3},a_{2}+1\right] , & \text{\emph{if \ }}%
a_{1}=1
\end{array}
\right. \   \label{fquot}
\end{equation}
\emph{resp.} 
\begin{equation}
\ \frac{\mathsf{Q}_{i}}{\mathsf{Q}_{i-1}}=\left\{ 
\begin{array}{ll}
\left[ a_{i},a_{i-1},\ldots ,a_{3},a_{2}\right] , & \text{\emph{if \ }}%
a_{2}\geq 2\text{\emph{\ }} \\ 
\left[ a_{i},a_{i-1},\ldots ,a_{4},a_{3}+1\right] , & \text{\emph{if \ }}%
a_{2}=1
\end{array}
\right.   \label{fquot2}
\end{equation}
\emph{(b)} For all $i$\emph{, }$0\leq i\leq \nu $\emph{,} 
\begin{equation}
\mathsf{P}_{i}\ \mathsf{Q}_{i-1}-\mathsf{P}_{i-1}\ \mathsf{Q}_{i}=\left(
-1\right) ^{i}\smallskip   \label{ano-kato}
\end{equation}
In particular, for $i$, $1\leq i\leq \nu $, we have \emph{gcd}$\left( 
\mathsf{P}_{i}\ ,\mathsf{Q}_{i}\right) =1$, and\smallskip\ \smallskip\ 
\begin{equation}
\frac{\mathsf{P}_{1}}{\mathsf{Q}_{1}}<\frac{\mathsf{P}_{3}}{\mathsf{Q}_{3}}%
<\cdots \leq \frac{\mathsf{P}_{\nu }}{\mathsf{Q}_{\nu }}=\frac{\kappa }{%
\lambda }\leq \cdots <\frac{\mathsf{P}_{4}}{\mathsf{Q}_{4}}<\frac{\mathsf{P}%
_{2}}{\mathsf{Q}_{2}}  \label{de-in}
\end{equation}
$\emph{(iv)}$ For $a\in \Bbb{Z}$, setting 
\[
\mathcal{M}^{+}\left( a\right) :=\left( 
\begin{array}{ll}
a & 1 \\ 
1 & 0
\end{array}
\right) \in \text{ }\emph{SL}\left( 2,\Bbb{Z}\right) 
\]
we get 
\begin{equation}
\left( 
\begin{array}{ll}
\mathsf{P}_{i} & \mathsf{P}_{i-1} \\ 
\mathsf{Q}_{i} & \mathsf{Q}_{i-1}
\end{array}
\right) =\mathcal{M}^{+}\left( a_{1}\right) \cdot \mathcal{M}^{+}\left(
a_{2}\right) \cdot \mathcal{M}^{+}\left( a_{3}\right) \ \cdots \ \mathcal{M}%
^{+}\left( a_{i}\right)   \label{psuc}
\end{equation}
\end{lemma}

\noindent \textit{Proof. }For (i) one uses the standard euclidean division
algorithm by defining successively 
\[
y_{1}:=\frac{\kappa }{\lambda },\ \ \ y_{i+1}:=\left( y_{i}-\left\lfloor
y_{i}\right\rfloor \right) ^{-1},\ \ a_{i}:=\left\lfloor y_{i}\right\rfloor
,\ \ i=1,2,\ldots 
\]
(ii) is obvious by the definition of the recurrence relation (\ref{PQ}%
).\smallskip\ \newline
(iii) For the proof of (\ref{fquot}) write this quotient as 
\[
\frac{\mathsf{P}_{i}}{\mathsf{P}_{i-1}}=a_{i}+\frac{\mathsf{P}_{i-2}}{%
\mathsf{P}_{i-1}}=a_{i}+\ 1\ /\ \left( \frac{\mathsf{P}_{i-1}}{\mathsf{P}%
_{i-2}}\right) \ . 
\]
(\ref{fquot2}) is similar. (\ref{ano-kato}) can be shown by induction on $i$
and implies gcd$\left( \mathsf{P}_{i}\ ,\mathsf{Q}_{i}\right) =1$ and (\ref
{de-in}).\smallskip \newline
(iv) Write 
\[
\left( 
\begin{array}{ll}
\mathsf{P}_{i} & \mathsf{P}_{i-1} \\ 
\mathsf{Q}_{i} & \mathsf{Q}_{i-1}
\end{array}
\right) =\left( 
\begin{array}{ll}
\mathsf{P}_{i-1} & \mathsf{P}_{i-2} \\ 
\mathsf{Q}_{i-1} & \mathsf{Q}_{i-2}
\end{array}
\right) \cdot \mathcal{M}^{+}\left( a_{i}\right) 
\]
and use induction. $_{\Box }$

\begin{remark}
\label{dioph}\emph{For }$\lambda ,\kappa $\emph{\ as in \ref{POS}, all
integer solutions of the linear diophantine equation} 
\[
\kappa \cdot \frak{x}+\lambda \cdot \frak{x}^{\prime }=1
\]
\emph{can be read off directly from the regular continued fraction expansion
(\ref{posconf}) of }$\frac{\kappa }{\lambda }$\emph{. They are of the form } 
\[
\frak{x}=\frak{x}_{0}+\xi \cdot \lambda ,\ \ \ \frak{x}^{\prime }=\frak{x}%
_{0}^{\prime }-\xi \cdot \kappa ,\ \ \ \xi \in \Bbb{Z},
\]
\emph{where} 
\[
\frak{x}_{0}=\varepsilon \ \mathsf{Q}_{\nu -1}=\left\{ 
\begin{array}{ll}
\frac{\varepsilon \cdot \lambda }{\left[ a_{\nu },a_{\nu -1},\ldots
,a_{3},a_{2}\right] } & \text{\emph{if} \ }a_{2}\geq 2 \\ 
\  &  \\ 
\frac{\varepsilon \cdot \lambda }{\left[ a_{\nu },a_{\nu -1},\ldots
,a_{4},a_{3}+1\right] } & \text{\emph{if} \ }a_{2}=1
\end{array}
\right. 
\]
\emph{and} 
\[
\frak{x}_{0}^{\prime }=-\varepsilon \ \mathsf{P}_{\nu -1}=\left\{ 
\begin{array}{ll}
\frac{-\varepsilon \cdot \kappa }{\left[ a_{\nu },a_{\nu -1},\ldots
,a_{2},a_{1}\right] } & \text{\emph{if} \ }a_{1}\geq 2 \\ 
\  &  \\ 
\frac{-\varepsilon \cdot \kappa }{\left[ a_{\nu },a_{\nu -1},\ldots
,a_{3},a_{2}+1\right] } & \text{\emph{if} \ }a_{1}=1
\end{array}
\right. 
\]
\emph{with} 
\[
\varepsilon :=\left\{ 
\begin{array}{rr}
1, & \text{\emph{if} \ }\nu \ \text{\emph{is even}} \\ 
-1, & \text{\emph{if \ }\ }\nu \ \text{\emph{is odd}}
\end{array}
\right. 
\]
\emph{(This is an immediate consequence of the equalities (\ref{fquot}), (%
\ref{fquot2}) and (\ref{ano-kato}) for }$i=\nu $\emph{).}
\end{remark}

\begin{remark}
\label{RELE}\emph{\ As it follows from a theorem of \ Lam\'{e} \cite{Lame},
the length }$\nu $\emph{\ of (\ref{posconf}) is smaller than }$5$ \emph{%
multiplied by the number of digits in the decimal expansion of }$\lambda $%
\emph{. Fixing }$\kappa $\emph{, we may more precisely say that, since the
smallest positive integer }$\lambda $\emph{\ for which the regular continued
fraction expansion (\ref{posconf}) of }$\frac{\kappa }{\lambda }$\emph{\
takes a given value }$\nu =j$\emph{, is the }$\left( j+1\right) $\emph{%
-Fibonacci number }$\mathsf{Fib}\left( j+1\right) $\emph{,} 
\[
\mathsf{Fib}\left( j\right) :=\frac{\frak{k}^{\,j}-\left( -\frak{k}\right)
^{-j}}{\frak{k}+\frak{k}^{\,-1}},\,\,\,\,\,\frak{k}:=\frac{1+\sqrt{5}}{2}\,,
\]
\emph{the length }$\nu $\emph{\ can be bounded from above by } 
\begin{equation}
\nu \leq \frac{1+\text{\emph{log}}_{10}\left( \lambda \right) }{\text{\emph{%
log}}_{10}\left( \frak{k}\right) }=\left( 2.0780...\right) \left( 1+\text{%
\emph{log}}_{10}\left( \lambda \right) \right) <5\cdot \text{\emph{log}}%
_{10}\left( \lambda \right)   \label{estimation}
\end{equation}
\emph{(For the early history on the estimations of }$\nu $\emph{\ the reader
is referred to Shallit \cite{Shallit}. For better approximations of }$\nu $ 
\emph{in certain number-regions, see Dixon \cite{Dixon} and Kilian \cite
{Kilian}.) }
\end{remark}

\noindent The analogue of lemma \ref{POS} for the \textit{negative}-regular
continued fractions is formulated as follows.

\begin{lemma}
Let $\lambda ,\kappa $ be two integers with $0<\lambda <\kappa ,$ \emph{gcd}$%
\left( \lambda ,\kappa \right) =1$. Then\emph{:\smallskip }\newline
\emph{(i)} There exists always a uniquely determined negative-regular
continued fraction expansion 
\begin{equation}
\text{ }\frac{\kappa }{\lambda }=\left[ \!\left[ c_{1},c_{2},\ldots ,c_{\rho
}\right] \!\right]   \label{mconf}
\end{equation}
of $\frac{\kappa }{\lambda }$ \emph{(}with $c_{j}\geq 2$\emph{,}$\ \forall
j,\ \ 1\leq j\leq \rho $\emph{).}\smallskip \newline
\emph{(ii)} Defining the finite sequences $\left( \mathsf{R}_{i}\right)
_{-1\leq i\leq \rho }$ and $\left( \mathsf{S}_{i}\right) _{-1\leq i\leq \rho
}$ by 
\begin{equation}
\left\{ 
\begin{array}{llll}
\mathsf{R}_{-1}=0, & \mathsf{R}_{0}=1, & \mathsf{R}_{1}=c_{1}, & \mathsf{R}%
_{i}=c_{i}\ \mathsf{R}_{i-1}-\mathsf{R}_{i-2},\ \forall i,\ \ 2\leq i\leq
\rho , \\ 
&  & \  &  \\ 
\mathsf{S}_{-1}=-1, & \mathsf{S}_{0}=0, & \mathsf{S}_{1}=1, & \mathsf{S}%
_{i}=c_{i}\ \mathsf{S}_{i-1}-\mathsf{S}_{i-2},\ \forall i,\ \ 2\leq i\leq
\rho ,
\end{array}
\right.   \label{RS}
\end{equation}
we obtain 
\[
\frac{\mathsf{R}_{i}}{\mathsf{S}_{i}}=\left[ \!\left[ c_{1},c_{2},\ldots
,c_{i}\right] \!\right] ,\ \forall i,\ \ 1\leq i\leq \rho .
\]
\emph{(This is the corresponding }$i$\textit{-th convergent}\emph{\ of the
continued fraction (\ref{mconf})).\smallskip }\newline
\emph{(iii)} The above sequences satisfy the following conditions\emph{%
:\smallskip }\newline
\emph{(a) }For all $i$, $1\leq i\leq \rho $ \emph{(resp. }$2\leq i\leq \rho $%
\emph{)}, 
\begin{equation}
\frac{\mathsf{R}_{i}}{\mathsf{R}_{i-1}}=\left[ \!\left[ c_{i},c_{i-1},\ldots
,c_{2},c_{1}\right] \!\right]   \label{fquot3}
\end{equation}
\emph{resp.} 
\begin{equation}
\frac{\mathsf{S}_{i}}{\mathsf{S}_{i-1}}=\left[ \!\left[ c_{i},c_{i-1},\ldots
,c_{3},c_{2}\right] \!\right]   \label{fquot4}
\end{equation}
\emph{(b)} For all $i$, $0\leq i\leq \rho $, 
\begin{equation}
\mathsf{R}_{i-1}\ \mathsf{S}_{i}-\mathsf{R}_{i}\ \mathsf{S}%
_{i-1}=1.\smallskip   \label{anw}
\end{equation}
In particular, for $i$, $1\leq i\leq \rho $, we have \emph{gcd}$\left( 
\mathsf{R}_{i}\ ,\mathsf{S}_{i}\right) =1$, and\smallskip \smallskip\
\smallskip\ 
\begin{equation}
\frac{\kappa }{\lambda }=\frac{\mathsf{R}_{\rho }}{\mathsf{S}_{\rho }}<\frac{%
\mathsf{R}_{\rho -1}}{\mathsf{S}_{\rho -1}}<\cdots <\frac{\mathsf{R}_{3}}{%
\mathsf{S}_{3}}<\frac{\mathsf{R}_{2}}{\mathsf{S}_{2}}<\frac{\mathsf{R}_{1}}{%
\mathsf{S}_{1}}  \label{sdec}
\end{equation}
$\emph{(iv)}$ For $c\in \Bbb{Z}$, setting 
\[
\mathcal{M}^{-}\left( c\right) :=\left( 
\begin{array}{rr}
c & -1 \\ 
1 & 0
\end{array}
\right) \in \text{ }\emph{SL}\left( 2,\Bbb{Z}\right) 
\]
we get for all $i$, $1\leq i\leq \rho $, 
\begin{equation}
\left( 
\begin{array}{ll}
\mathsf{R}_{i} & -\mathsf{R}_{i-1} \\ 
\mathsf{S}_{i} & -\mathsf{S}_{i-1}
\end{array}
\right) =\mathcal{M}^{-}\left( c_{1}\right) \cdot \mathcal{M}^{-}\left(
c_{2}\right) \cdot \mathcal{M}^{-}\left( c_{3}\right) \ \cdots \ \mathcal{M}%
^{-}\left( c_{i}\right)   \label{nsuc}
\end{equation}
\end{lemma}

\noindent \textit{Proof. }For (i) use the modified euclidean division
algorithm and define successively 
\begin{equation}
y_{1}:=\frac{\kappa }{\lambda },\ \ \ y_{i+1}:=\left( \left\lceil
y_{i}\right\rceil -y_{i}\right) ^{-1},\ \ c_{i}:=\left\lceil
y_{i}\right\rceil ,\ \ i=1,2,\ldots  \label{MEA}
\end{equation}
The proofs of the remaining assertions (ii)-(iv) are similar to those of \ref
{POS}. $_{\Box }$

\begin{remark}
\label{length1}\emph{It is easy to show that the length }$\rho $\emph{\ of (%
\ref{mconf}) equals }$\kappa -1$\emph{\ if and only if }$c_{i}=2$\emph{, for}
all $i$, $1\leq i\leq \rho $\emph{, and }$\rho \leq \QOVERD\lfloor \rfloor
{k-1}{2}$ \emph{otherwise. (For }$\kappa \gg 10\cdot $\emph{log}$_{10}\left(
\lambda \right) $\emph{, (\ref{estimation}) implies that }$\rho $\emph{\
might take values }$\gg \nu $\emph{; that this happens in the most
``generic'' case will become clear in rem. \ref{length2}).\medskip }
\end{remark}

\noindent \textsf{(b) }Already Lagrange \cite{Lagrange} and M\"{o}bius \cite
{Mob} knew for a wide palette of examples of ``non-standard'' continued
fraction expansions how one should modify them to get regular, i.e.,
``usual'' continued fractions by special substitutions. However, a
systematic ``mechanical method'' for expressing a given semi-regular
continued fraction via a regular one was first developed by Minkowski (cf. 
\cite{Minkowski}, pp. 116-118) and was further elucidated in Perron's
classical book \cite{Perron}, Kap. V, \S 40. In fact, in order to write down
(\ref{conf}) as a regular continued fraction, one has to define $\varepsilon
_{0}:=1,\varepsilon _{\nu }:=1$, to add a $\frac{1}{1+\cdots }$ in front of
every negative $\varepsilon _{i}$, to replace all minus signs by plus signs,
and finally to substitute $a_{i}-\frac{1}{2}\left( \left( 1-\varepsilon
_{i-1}\right) +\left( 1-\varepsilon _{i}\right) \right) $ for each $a_{i}$,
for all $i$, $1\leq i\leq \nu $. Since this method demands some extra
restrictive operations, like the elimination of the zeros and the
short-cutting of the probably occuring $1$'s at the last step, it seems to
be rather laborious. Moreover, since it is important for our later arguments
to provide a detailed relationship between the regular and the
negative-regular continued fraction expansion (by assuming the first one as
given), we shall prefer to use another, slightly different approach.
Proposition \ref{PL-M} describes the transfer process explicitly. This is
actually a ``folklore-type-result'' and can be found (without proof) in
Myerson's paper \cite{Myerson}, p. 424, who examines the average length of
negative-regular expansions, as well as in works of Hirzebruch and Zagier
(cf. \cite{Hirz2}, p. 215; \cite{Zagier}, p. 131) in a version concerning
infinite continued fractions and serves as auxiliary tool for studying
positive definite binary forms. Our proof is relatively short (compared with
the above mentioned repetitive substitutions) as it makes use only of formal
multiplications of certain $2\times 2$ unimodular matrices.

\begin{proposition}
\label{PL-M}If $\lambda ,\kappa \in \Bbb{Z}$ with $0<\lambda <\kappa ,$ 
\emph{gcd}$\left( \lambda ,\kappa \right) =1$, and 
\[
\frac{\kappa }{\lambda }=\left[ a_{1},a_{2},\ldots ,a_{\nu }\right] =\left[
\!\left[ c_{1},c_{2},\ldots ,c_{\rho }\right] \!\right] 
\]
are the regular and negative-regular continued fraction expansions of $\frac{%
\kappa }{\lambda }$, respectively, then $\left( c_{1},c_{2},\ldots ,c_{\rho
}\right) $, as ordered $\rho $-tuple, equals 
\[
\left( c_{1},c_{2},\ldots ,c_{\rho }\right) =
\]
\[
=\left\{ 
\begin{array}{ll}
\left( a_{1}+1,\stackunder{\left( a_{2}-1\right) \text{\emph{-times}}}{%
\underbrace{2,...,2}},a_{3}+2,\stackunder{\left( a_{4}-1\right) \text{\emph{%
-times}}}{\underbrace{2,...,2}},a_{5}+2,\ldots ,a_{\nu -1}+2,\stackunder{%
\left( a_{\nu }-1\right) \text{\emph{-times}}}{\underbrace{2,...,2}}\right) ,
& \text{\emph{if \ }}\nu \ \text{ }\emph{even} \\ 
\  &  \\ 
\left( a_{1}+1,\stackunder{\left( a_{2}-1\right) \text{\emph{-times}}}{%
\underbrace{2,...,2}},a_{3}+2,\stackunder{\left( a_{4}-1\right) \text{\emph{%
-times}}}{\underbrace{2,...,2}},a_{5}+2,\ldots ,\stackunder{\left( a_{\nu
-1}-1\right) \text{\emph{-times}}}{\underbrace{2,...,2}}a_{\nu }+1\right) ,
& \text{\emph{if\ \ }}\nu \ \text{ }\emph{odd.}
\end{array}
\right. 
\]
\end{proposition}

\noindent \textit{Proof. }At first define the matrices 
\[
\mathcal{U}:=\left( 
\begin{array}{ll}
1 & 1 \\ 
0 & 1
\end{array}
\right) \text{, \ \ \ }\mathcal{V}:=\left( 
\begin{array}{ll}
1 & 0 \\ 
1 & 1
\end{array}
\right) \text{, \ \ \ }\mathcal{W}:=\left( 
\begin{array}{rr}
1 & 0 \\ 
1 & -1
\end{array}
\right) 
\]
from SL$\left( 2,\Bbb{Z}\right) $. Obviously, $\mathcal{W}^{2}=$ Id, and for
all $a$, $a\in \Bbb{N}$, we have 
\begin{equation}
\mathcal{M}^{+}\left( a\right) =\mathcal{W}\cdot \mathcal{U}\cdot \mathcal{V}%
^{a-1}\ \   \label{prwti}
\end{equation}
and 
\begin{equation}
\left( \mathcal{M}^{-}\left( 2\right) \right) ^{a}=\mathcal{U}\cdot \mathcal{%
V}^{a}\cdot \mathcal{U}^{-1}  \label{devteri}
\end{equation}
Using (\ref{prwti}) and (\ref{devteri}) we obtain, in addition, 
\begin{equation}
\mathcal{M}^{+}\left( a\right) =\mathcal{M}^{-}\left( a+1\right) \cdot 
\mathcal{W}=\mathcal{W}\cdot \left( \mathcal{M}^{-}\left( 2\right) \right)
^{a-1}\cdot \mathcal{U}  \label{triti}
\end{equation}
and 
\begin{equation}
\mathcal{U}\cdot \mathcal{M}^{-}\left( a\right) =\mathcal{M}^{-}\left(
a+1\right)  \label{tetarti}
\end{equation}
for all $a$, $a\in \Bbb{N}$. If $\nu $ is even, then (\ref{triti}) implies 
\[
\mathcal{M}^{+}\left( a_{1}\right) \cdot \mathcal{M}^{+}\left( a_{2}\right)
\ \cdots \ \mathcal{M}^{+}\left( a_{\nu -1}\right) \cdot \mathcal{M}%
^{+}\left( a_{\nu }\right) =\prod_{j=1}^{\frac{\nu }{2}}\ \left( \mathcal{M}%
^{+}\left( a_{2j-1}\right) \cdot \mathcal{M}^{+}\left( a_{2j}\right) \right)
=\smallskip 
\]
\begin{eqnarray*}
&=&\prod_{j=1}^{\frac{\nu }{2}}\ \left( \left( \mathcal{M}^{-}\left(
a_{2j-1}+1\right) \cdot \mathcal{W}\right) \cdot \left( \mathcal{W}\cdot
\left( \mathcal{M}^{-}\left( 2\right) \right) ^{a_{2j}-1}\cdot \mathcal{U}%
\right) \right) =\smallskip \\
&=&\prod_{j=1}^{\frac{\nu }{2}}\ \left( \mathcal{M}^{-}\left(
a_{2j-1}+1\right) \cdot \left( \left( \mathcal{M}^{-}\left( 2\right) \right)
^{a_{2j}-1}\cdot \mathcal{U}\right) \right) \ \stackrel{\text{by (\ref
{tetarti})}}{=}\smallskip \  \\
&=&\left( \left( \mathcal{M}^{-}\left( a_{1}+1\right) \cdot \left( \mathcal{M%
}^{-}\left( 2\right) \right) ^{a_{2}-1}\right) \cdot \prod_{j=2}^{\frac{\nu 
}{2}}\ \left( \mathcal{M}^{-}\left( a_{2j-1}+2\right) \cdot \left( \left( 
\mathcal{M}^{-}\left( 2\right) \right) ^{a_{2j}-1}\right) \right) \right)
\cdot \mathcal{U}
\end{eqnarray*}
and the desired relation follows directly from the equalities (\ref{psuc})
and (\ref{nsuc}) of the previous lemmas because $\frac{\kappa }{\lambda }=%
\frac{\mathsf{R}_{\rho }}{\mathsf{S}_{\rho }}=\frac{\mathsf{P}_{\nu }}{%
\mathsf{Q}_{\nu }}$. If $\nu $ is odd, then using the above argument for the
first $\nu -1$ factors of the product of the corresponding $\mathcal{M}^{+}$%
-matrices, we get 
\[
\mathcal{M}^{+}\left( a_{1}\right) \cdot \mathcal{M}^{+}\left( a_{2}\right)
\ \cdots \ \mathcal{M}^{+}\left( a_{\nu -1}\right) \cdot \mathcal{M}%
^{+}\left( a_{\nu }\right) =\smallskip 
\]
\begin{eqnarray*}
&=&\left( \left( \mathcal{M}^{-}\left( a_{1}+1\right) \cdot \left( \mathcal{M%
}^{-}\left( 2\right) \right) ^{a_{2}-1}\right) \cdot \prod_{j=2}^{\frac{\nu
-1}{2}}\ \left( \mathcal{M}^{-}\left( a_{2j-1}+2\right) \cdot \left( \left( 
\mathcal{M}^{-}\left( 2\right) \right) ^{a_{2j}-1}\right) \right) \right)
\cdot \mathcal{U}\cdot \mathcal{M}^{+}\left( a_{\nu }\right) =\smallskip \\
&=&\mathcal{M}^{-}\left( a_{1}+1\right) \cdot \left( \mathcal{M}^{-}\left(
2\right) \right) ^{a_{2}-1}\cdot \left( \prod_{j=2}^{\frac{\nu -1}{2}}\ 
\mathcal{M}^{-}\left( a_{2j-1}+2\right) \cdot \left( \mathcal{M}^{-}\left(
2\right) \right) ^{a_{2j}-1}\right) \cdot \left( \mathcal{U}\cdot \mathcal{M}%
^{-}\left( a_{\nu }+1\right) \cdot \mathcal{W}\right) .
\end{eqnarray*}
Since\ $\mathcal{U}\cdot \mathcal{M}^{-}\left( a_{\nu }+1\right) \cdot 
\mathcal{W}=\mathcal{M}^{+}\left( a_{\nu }+1\right) $ has the same first
column-vector as $\mathcal{M}^{-}\left( a_{\nu }+1\right) $ the conclusion
follows again from (\ref{psuc}) and (\ref{nsuc}). $_{\Box }$

\begin{remark}
\label{matr-transf}\emph{From the proof of \ref{PL-M} one easily verifies
that for all indices }$i$\emph{,} \emph{\ }$1\leq i\leq \QTOVERD\lfloor
\rfloor {\nu }{2}$\emph{, } \emph{the convergents of these two continued
fraction expansions of }$\frac{\kappa }{\lambda }$ \emph{are connected by
the following ``matrix-transfer rule'':} 
\[
\left( 
\begin{array}{ll}
\mathsf{P}_{2i} & \mathsf{P}_{2i-1} \\ 
\mathsf{Q}_{2i} & \mathsf{Q}_{2i-1}
\end{array}
\right) =\left( 
\begin{array}{ll}
\mathsf{R}_{a_{2}+a_{4}+\cdots +a_{2i}} & -\mathsf{R}_{a_{2}+a_{4}+\cdots
+a_{2i}-1} \\ 
\mathsf{S}_{a_{2}+a_{4}+\cdots +a_{2i}} & -\mathsf{S}_{a_{2}+a_{4}+\cdots
+a_{2i}-1}
\end{array}
\right) \cdot \mathcal{U}
\]
\end{remark}

\begin{remark}
\label{length2}\emph{By prop. \ref{PL-M} the length }$\rho $\emph{\ of the
negative-regular continued fraction expansion (\ref{mconf}) of }$\frac{%
\kappa }{\lambda }$\emph{\ equals} 
\[
\rho =\left\{ 
\begin{array}{ll}
\sum_{i=1}^{\frac{\nu }{2}}\ a_{2i},\smallskip  & \text{\emph{if}\ }\nu 
\text{ \ \emph{even} } \\ 
\left( \sum_{i=1}^{\frac{\nu -1}{2}}\ a_{2i}\right) +1, & \text{\emph{if\ } }%
\nu \text{ \ \emph{odd}}
\end{array}
\right. 
\]
\emph{If for at least one }$i,$\emph{\ }$1\leq i\leq \rho $\emph{, we have }$%
c_{i}\geq 3$\emph{\ (cf. rem. \ref{length1}), we get\smallskip } 
\[
1\leq \QOVERD\lfloor \rfloor {\nu +1}{2}\leq \rho \leq \QOVERD\lfloor
\rfloor {\kappa -1}{2}\ .
\]
\emph{Obviously, if} \emph{all (or allmost all) }$a_{2i}$\emph{'s are} $\gg 2
$ \emph{(assumption which expresses the ``generic'' case), then }$\rho \gg
\nu $\emph{, and the number }$\rho -\nu $ \emph{of the extra modified
euclidean division algorithms (\ref{MEA}) (coming from the extra ``twos'')
which one needs in order to determine }\textit{directly} \emph{the }$\left(
-\right) $\emph{-sign continued fraction expansion of }$\frac{\kappa }{%
\lambda }$\emph{\ (i.e., without using prop. \ref{PL-M}) may become
tremendously large!\medskip }
\end{remark}

\noindent

\noindent\textsf{(c) }The
``lattice-geometry'' of two-dimensional rational s.c.p. cones is completely
describable by means of just two (relatively prime) integers
(``parameters'').

\begin{lemma}
\label{HNF}Let $N$ be a lattice of rank $2$ and $\sigma \subset N_{\Bbb{R}}$
a two-dimensional rational s.c.polyhedral cone with \emph{Gen}$\left( \sigma
\right) =\left\{ n_{1},n_{2}\right\} $. Then there exist a $\Bbb{Z}$-basis $%
\left\{ \frak{y}_{1},\frak{y}_{2}\right\} $ of $N$\emph{\ }and two integers $%
p=p_{\sigma },\ q=q_{\sigma }\in \Bbb{Z}_{\geq 0}$ with $0\leq p<q$\emph{,
gcd}$\left( p,q\right) =1$, such that 
\[
n_{1}=\frak{y}_{1},\ \ \ n_{2}=p\,\frak{y}_{1}+q\,\frak{y}_{2},\ \ \ q=\text{%
\emph{mult}}\left( \sigma ;N\right) =\frac{\text{\emph{det}}\left( \Bbb{Z\ }%
n_{1}\oplus \Bbb{Z\ }n_{2}\right) }{\text{\emph{det}}\left( N\right) }\ .
\]
Moreover, if $\Phi $ is a $\Bbb{Z}$-module isomorphism $\Phi :N\rightarrow N$%
, then the above property is preserved by the same numbers $p=p_{\Phi _{\Bbb{%
R}}\left( \sigma \right) },\ q=q_{\Phi _{\Bbb{R}}\left( \sigma \right) }$
for the cone $\Phi _{\Bbb{R}}\left( \sigma \right) \subset N_{\Bbb{R}}$ with
respect to the $\Bbb{Z}$-basis $\left\{ \Phi \left( \frak{y}_{1}\right)
,\Phi \left( \frak{y}_{2}\right) \right\} $ of the lattice $N$.
\end{lemma}

\noindent \textit{Proof.} Choose an arbitrary $\Bbb{Z}$-basis $\left\{ \frak{%
n}_{1},\frak{n}_{2}\right\} $ of $N$ with $\frak{n}_{1}=n_{1}$. Since $%
\sigma $ is $2$-dimensional, it is also simplicial; this means that $n_{2}$
may be expressed as a linear combination of the members of this $\Bbb{Z}$%
-basis having the form: 
\[
n_{2}=\lambda _{1}\ \frak{n}_{1}+\lambda _{2}\ \frak{n}_{2},\ \ \ \text{
with \ \ }\lambda _{1}\ \in \Bbb{Z},\ \ \lambda _{2}\in \left( \Bbb{Zr}%
\left\{ 0\right\} \right) \ . 
\]
Now define $q:=\left| \lambda _{2}\right| $ and $p:=\left[ \lambda _{1}%
\right] _{q}$. Since $0\leq p<q$ and 
\[
n_{2}=p\ \frak{n}_{1}+q\ \left( \text{sgn}\left( \lambda _{2}\right) \ \frak{%
n}_{2}+\frac{\lambda _{1}-p}{q}\ \frak{n}_{1}\right) \ , 
\]
it suffices to consider $\frak{y}_{1}:=\frak{n}_{1}$ and $\frak{y}_{2}:=$ sgn%
$\left( \lambda _{2}\right) \frak{n}_{2}+\frac{\lambda _{1}-p}{q}\frak{n}%
_{1} $. Furthermore, gcd$\left( p,q\right) =1$ and $q=$ mult$\left( \sigma
;N\right) $, because $n_{1},n_{2}$ are primitive. The last assertion is
obvious. $_{\Box }$

\begin{definition}
\label{PQC}\emph{If }$N$\emph{\ is a lattice of rank }$2$\emph{\ and }$%
\sigma \subset N_{\Bbb{R}}$\emph{\ a two-dimensional rational s.c.polyhedral
cone with Gen}$\left( \sigma \right) =\left\{ n_{1},n_{2}\right\} $\emph{,
then we call }$\sigma $\emph{\ a }$\left( p,q\right) $\emph{-}\textit{cone
w.r.t. the basis\emph{\ }}$\left\{ \frak{y}_{1},\frak{y}_{2}\right\} $\emph{%
, if } $p=p_{\sigma },\ q=q_{\sigma }$ \emph{as in lemma \ref{HNF}. (To
avoid confusion, we should stress at this point that saying ``w.r.t. the
basis }$\left\{ \frak{y}_{1},\frak{y}_{2}\right\} $\emph{'' we just indicate
the choice of one suitable }$\Bbb{Z}$\emph{-basis of }$N$ \emph{among} 
\textit{all} \emph{its }$\Bbb{Z}$\emph{-bases in order to apply lemma \ref
{HNF} for }$\sigma $\emph{; but, of course, if }$\left\{ \frak{y}_{1},%
\widehat{\frak{y}}_{2}\right\} $ \emph{were a} $\Bbb{Z}$\emph{-basis of }$N$ 
\emph{having} \emph{the} \textit{same }\emph{property, i.e.,} $n_{2}=%
\widehat{\,p}\,\frak{y}_{1}+\widehat{q}\,\widehat{\frak{y}}_{2}$\emph{,} $%
0\leq \widehat{\,p}<\widehat{q}$\emph{, gcd}$\left( \widehat{\,p},\widehat{q}%
\right) =1$\emph{, then obviously }$\,p=\widehat{\,p}$\emph{\ and }$q=%
\widehat{q}$\emph{, i.e.,} $\widehat{\frak{y}}_{2}=\frak{y}_{2}$\emph{!)}
\end{definition}

\begin{remark}
\emph{For }$p,q$\emph{\ as in lemma \ref{HNF}, there is a uniquely
determined integer }$p^{\prime }=p_{\sigma }^{\prime }$\emph{, }$0\leq
p^{\prime }<q$\emph{, such that} 
\[
pp^{\prime }\equiv 1\text{\emph{(mod} }q\text{\emph{), \ \ \ (i.e., \ }}%
\left[ \text{\emph{\ }}pp^{\prime }\right] _{q}=1\text{\emph{)\ .}}
\]
$p^{\prime }$\emph{\ is often called} \textit{the socius} \emph{of }$p$\emph{%
. If }$p\neq 0$ \emph{(which means that }$q\neq 1$\emph{), then using a
formula due to Voronoi (see \cite{U-H}, p. 183), }$p^{\prime }$\emph{\ can
be written as} 
\[
p^{\prime }=\left[ 3-2p+6\ \left( \sum_{j=1}^{p-1}\ \left( \QDOVERD\lfloor
\rfloor {j\ q}{p}\right) ^{2}\right) \right] _{q}\ \ .
\]
\end{remark}

\begin{proposition}
\label{ISO2}Let $N$ be a lattice of rank $2$ and $\sigma ,\tau \subset N_{%
\Bbb{R}}$ two $2$-dimensional rational s.c.p. cones with \emph{Gen}$\left(
\sigma \right) =\left\{ n_{1},n_{2}\right\} $\emph{,} \emph{Gen}$\left( \tau
\right) =\left\{ u_{1},u_{2}\right\} $. Then the following conditions are
equivalent\emph{:\smallskip }\newline
\emph{(i) \ }There exists an isomorphism of germs\emph{: }$\left( U_{\sigma }%
\text{\emph{, orb}}\left( \sigma \right) \right) \cong \left( U_{\tau }\text{%
\emph{, orb}}\left( \tau \right) \right) $.\smallskip \newline
\emph{(ii)} There exists a $\Bbb{Z}$-module isomorphism $\varpi
:N\rightarrow N$, whose scalar extension $\varpi =\varpi _{\Bbb{R}}:N_{\Bbb{R%
}}\rightarrow N_{\Bbb{R}}$ \ has the property\emph{: }$\varpi \left( \sigma
\right) =\tau .\smallskip $\newline
\emph{(iii)} For the numbers $p_{\sigma },$ $p_{\tau },$ $q_{\sigma },$ $%
q_{\tau }$ associated to $\sigma ,\tau $ w.r.t. a $\Bbb{Z}$-basis $\left\{ 
\frak{y}_{1},\frak{y}_{2}\right\} $ of $N$ \emph{(}as in \emph{\ref{HNF})}
we have 
\begin{equation}
q_{\tau }=q_{\sigma }\text{ \ \ \ \emph{and} \ \ \ }\left\{ 
\begin{array}{ll}
\text{\emph{either}} & p_{\tau }=p_{\sigma } \\ 
\text{\emph{or}} & p_{\tau }=p_{\sigma }^{\prime }
\end{array}
\right.   \label{isorel}
\end{equation}
\end{proposition}

\noindent \textit{Proof. }For the equivalence (i)$\Leftrightarrow $(ii) see
Ewald \cite{Ewald}, Ch. VI, thm. 2.11, pp. 222-223.\smallskip \newline
(ii)$\Rightarrow $(iii): Since $\left\{ \varpi \left( \frak{y}_{1}\right)
,\varpi \left( \frak{y}_{2}\right) \right\} $ is a $\Bbb{Z}$-basis of $N$
too, there exists a matrix $\mathcal{A}\in $ SL$\left( 2,\Bbb{Z}\right) $
such that 
\begin{equation}
\left( \varpi \left( \frak{y}_{1}\right) ,\varpi \left( \frak{y}_{2}\right)
\right) =\left( \frak{y}_{1},\frak{y}_{2}\right) \cdot \mathcal{A}
\label{BAS-CH}
\end{equation}
Now $\varpi \left( \sigma \right) =\tau $ implies $\varpi \left( \text{Gen}%
\left( \sigma \right) \right) =$ Gen$\left( \tau \right) $, i.e., 
\[
\text{either \ \ \ }\left( \varpi \left( n_{1}\right) =u_{1}\text{\ \ \& \ }%
\varpi \left( n_{2}\right) =u_{2}\right) \text{ \ \ \ or \ \ \ }\left(
\varpi \left( n_{1}\right) =u_{2}\ \text{\ \& \ }\varpi \left( n_{2}\right)
=u_{1}\right) \ . 
\]
In the first case we obtain 
\begin{equation}
\left( \varpi \left( n_{1}\right) ,\varpi \left( n_{2}\right) \right)
=\left( \varpi \left( \frak{y}_{1}\right) ,\varpi \left( \frak{y}_{2}\right)
\right) \ \left( 
\begin{array}{cc}
1 & p_{\sigma } \\ 
0 & q_{\sigma }
\end{array}
\right) =\left( \frak{y}_{1},\frak{y}_{2}\right) \ \left( 
\begin{array}{cc}
1 & p_{\tau } \\ 
0 & q_{\tau }
\end{array}
\right)  \label{BC1}
\end{equation}
In the second case: 
\begin{equation}
\left( \varpi \left( n_{1}\right) ,\varpi \left( n_{2}\right) \right)
=\left( \varpi \left( \frak{y}_{1}\right) ,\varpi \left( \frak{y}_{2}\right)
\right) \ \left( 
\begin{array}{cc}
1 & p_{\sigma } \\ 
0 & q_{\sigma }
\end{array}
\right) =\left( \frak{y}_{1},\frak{y}_{2}\right) \ \left( 
\begin{array}{cc}
p_{\tau } & 1 \\ 
q_{\tau } & 0
\end{array}
\right)  \label{BC2}
\end{equation}
Thus, by (\ref{BAS-CH}), (\ref{BC1}) and (\ref{BC2}) we get 
\[
\mathcal{A}=\left( 
\begin{array}{cc}
1 & \dfrac{p_{\tau }-p_{\sigma }}{q_{\sigma }}\smallskip \\ 
0 & \dfrac{q_{\tau }}{q_{\sigma }}
\end{array}
\right) \text{ \ \ \ \ and \ \ \ \ }\mathcal{A}=\allowbreak \allowbreak
\left( 
\begin{array}{cc}
p_{\tau } & \dfrac{1-p_{\sigma }\ p_{\tau }}{q_{\sigma }}\smallskip \\ 
q_{\tau } & -\dfrac{p_{\sigma }\ q_{\tau }}{q_{\sigma }}
\end{array}
\right) \ , 
\]
respectively. In the first case det$\left( \mathcal{A}\right) $ has to be
equal to $1$, which means that $q_{\sigma }=q_{\tau }$ and $p_{\tau
}-p_{\sigma }\equiv 0$ (mod $q_{\sigma }$), i.e., $p_{\tau }=p_{\sigma }$
(because $0\leq p_{\sigma },p_{\tau }\leq q_{\sigma }=q_{\tau }$). In the
second case, det$\left( \mathcal{A}\right) =-1$; hence, $q_{\sigma }=q_{\tau
}$ and $1-p_{\sigma }\ p_{\tau }\equiv 0$ (mod $q_{\sigma }$), i.e., $%
p_{\tau }=p_{\sigma }^{\prime }$.\smallskip\ \newline
(ii)$\Leftarrow $(iii): If $q_{\sigma }=q_{\tau }$ and $p_{\sigma }=p_{\tau
} $, we define $\varpi :=$ id$_{N_{\Bbb{R}}}$. Otherwise, $q_{\sigma
}=q_{\tau }$ and $p_{\tau }=p_{\sigma }^{\prime }$, and for an $\mathbf{x}%
\in N$ with $\mathbf{x}=\lambda _{1}\frak{y}_{1}+\lambda _{2}\frak{y}_{2}$, (%
$\lambda _{1} $, $\lambda _{2}\in \Bbb{Z}$), we set 
\[
\varpi \left( \mathbf{x}\right) :=\frac{1}{q_{\sigma }}\left( \lambda _{2}\
u_{1}+\left( \lambda _{1}q_{\sigma }-p_{\sigma }\lambda _{2}\right) \
u_{2}\right) \ . 
\]
Its scalar extension $\varpi =\varpi _{\Bbb{R}}:N_{\Bbb{R}}\rightarrow N_{%
\Bbb{R}}$ is the $\Bbb{R}$-vector space isomorphism with the desired
property. $_{\Box }$

\begin{remark}
\emph{Up to replacement of }$p$ \emph{by its socius }$p^{\prime }$ \emph{%
(which corresponds just to the interchange of the analytic coordinates),
these two numbers }$p$ \emph{and} $q$ \emph{parametrize uniquely the
analytic isomorphism class of the germ }$\left( U_{\sigma }\text{\emph{, orb}%
}\left( \sigma \right) \right) $. \emph{(In the terminology which will be
introduced in the next section, if }$q\neq 1$\emph{, then }$\left( U_{\sigma
}\text{\emph{, orb}}\left( \sigma \right) \right) $\emph{\ is nothing but a
cyclic quotient singularity of ``type'' }$\frac{1}{q}\left( q-p,1\right) $%
\emph{, cf. Oda \cite{Oda}, 1.24; prop. \ref{ISO2} can be therefore regarded
as a special case of the isomorphism criterion \ref{ISOCYC}).}
\end{remark}

\begin{lemma}
Let $N$ be a lattice of rank $2$, $M=$ \emph{Hom}$_{\Bbb{Z}}\left( N,\Bbb{Z}%
\right) $ its dual and $\sigma \subset N_{\Bbb{R}}$ a two-dimensional $%
\left( p,q\right) $-cone w.r.t. a $\Bbb{Z}$-basis $\left\{ \frak{y}_{1},%
\frak{y}_{2}\right\} $ of $N$. If we denote by $\left\{ \frak{m}_{1},\frak{m}%
_{2}\right\} $ the dual $\Bbb{Z}$-basis of $\left\{ \frak{y}_{1},\frak{y}%
_{2}\right\} $ in $M$, then the cone $\sigma ^{\vee }\subset M_{\Bbb{R}}$ is
a $\left( q-p,q\right) $-cone w.r.t. $\left\{ \frak{m}_{2},\frak{m}_{1}-%
\frak{m}_{2}\right\} $.
\end{lemma}

\noindent \textit{Proof. }Let Gen$\left( \sigma \right) =\left\{
n_{1},n_{2}\right\} $ be the two minimal generators of $\sigma $ with $n_{1}=%
\frak{y}_{1}$, $n_{2}=p\,\frak{y}_{1}+q\,\frak{y}_{2}$. Then 
\[
\sigma ^{\vee }=\text{ pos}\left( \left\{ \frak{m}_{2},q\frak{m}_{1}-p\frak{m%
}_{2}\right\} \right) =\text{ pos}\left( \left\{ \frak{m}_{2},\left(
q-p\right) \,\frak{m}_{2}+q\,\left( \frak{m}_{1}-\frak{m}_{2}\right)
\right\} \right) 
\]
and Gen$\left( \sigma ^{\vee }\right) =\left\{ \frak{m}_{2},\left(
q-p\right) \,\frak{m}_{2}+q\,\left( \frak{m}_{1}-\frak{m}_{2}\right)
\right\} $. Since $\left\{ \frak{m}_{2},\frak{m}_{1}-\frak{m}_{2}\right\} $
is a $\Bbb{Z}$-basis of $M$ and 
\[
0<q-p<q,\,\ \text{ gcd}\left( q-p,q\right) =1, 
\]
we are done. $_{\Box }\bigskip $

\noindent \textsf{(d) }From now on, and for the rest of the present section,
we fix a lattice $N$ of rank $2$, its dual $M$, a \ \textit{non-basic }%
two-dimensional $\left( p,q\right) $-cone $\sigma \subset N_{\Bbb{R}}$
w.r.t. a $\Bbb{Z}$-basis $\left\{ \frak{y}_{1},\frak{y}_{2}\right\} $ of $N$%
, the dual basis $\left\{ \frak{m}_{1},\frak{m}_{2}\right\} $ of $\left\{ 
\frak{y}_{1},\frak{y}_{2}\right\} $ in $M$, and the dual cone $\sigma ^{\vee
}\subset M_{\Bbb{R}}$ of $\sigma $. Moreover, we consider \textit{both} $%
\left( +\text{ \& }-\right) $-sign continued fraction expansions of both
rationals $\tfrac{q}{q-p}$ and $\tfrac{q}{p}$: 
\begin{equation}
\fbox{$
\begin{array}{ccc}
&  &  \\ 
& 
\begin{array}{c}
\ \ \ \ \ \ \ \ \ \ \ \ \ \ \dfrac{q}{q-p}=\left[ a_{1}^{\vee },a_{2}^{\vee
},\ldots ,a_{\nu }^{\vee }\right] =\left[ \!\left[ b_{1},b_{2},\ldots
,b_{\rho }\right] \!\right] \\ 
\  \\ 
\dfrac{q}{p}=\dfrac{q}{q-\left( q-p\right) }=\left[ a_{1},a_{2},\ldots ,a_{k}%
\right] =\left[ \!\left[ b_{1}^{\vee },b_{2}^{\vee },\ldots ,b_{t}^{\vee }%
\right] \!\right]
\end{array}
&  \\ 
&  & 
\end{array}
$}  \label{Q-P}
\end{equation}
and 
\[
\left( \frac{\mathsf{P}_{i}^{\vee }}{\mathsf{Q}_{i}^{\vee }}\right) _{-1\leq
i\leq \nu },\ \ \left( \frac{\mathsf{R}_{i}}{\mathsf{S}_{i}}\right) _{-1\leq
i\leq \rho },\ \ \left( \frac{\mathsf{P}_{i}}{\mathsf{Q}_{i}}\right)
_{-1\leq i\leq k},\ \ \left( \frac{\mathsf{R}_{i}^{\vee }}{\mathsf{S}%
_{i}^{\vee }}\right) _{-1\leq i\leq t} 
\]
the corresponding finite sequences of their convergents. It is well-known
(cf. \cite{Riem}, p. 223, \cite{Oda}, p. 29) that 
\[
\left( b_{1}+b_{2}+\cdots +b_{\rho }\right) -\rho =\left( b_{1}^{\vee
}+b_{2}^{\vee }+\cdots +b_{t}^{\vee }\right) -t=\rho +t-1\,. 
\]
On the other hand, examining only the two $\left( +\right) $-sign
expansions, the relationship between their entries can be written (strangely
enough) in a \textit{very} \textit{simple} form:

\begin{lemma}
For the ordered pair of the first two entries of the regular continued
fraction expansions \emph{(\ref{Q-P})} of $\frac{q}{q-p}$ and $\frac{q}{p}$
we obtain 
\begin{equation}
\left( a_{1}^{\vee },a_{1}\right) \notin \left\{ \left( 1,1\right) \right\}
\cup \left( \Bbb{Z}_{\geq 2}\right) ^{2}\   \label{MAKRIA}
\end{equation}
In particular, there are only two possibilities; namely either 
\begin{equation}
\left( a_{1}^{\vee }=1\ \ \text{\&\ \ }a_{1}\neq 1\right)
\Longleftrightarrow k=\nu -1  \label{PRW}
\end{equation}
or 
\begin{equation}
\left( a_{1}^{\vee }\neq 1\ \ \text{\&\ \ }a_{1}=1\right)
\Longleftrightarrow k=\nu +1  \label{DEU}
\end{equation}
In the first case, we have 
\begin{equation}
\fbox{$
\begin{array}{ccc}
&  &  \\ 
& a_{i}^{\vee }=\left\{ 
\begin{array}{ll}
a_{1}-1, & \text{\emph{for\ \ }}i=2 \\ 
a_{i-1}, & \text{\emph{for\ \ }}3\leq i\leq \nu \ \left( =k+1\right) 
\end{array}
\right.  &  \\ 
&  & 
\end{array}
$}  \label{A1}
\end{equation}
In the second case, 
\begin{equation}
\fbox{$
\begin{array}{ccc}
&  &  \\ 
& a_{i}^{\vee }=\left\{ 
\begin{array}{ll}
a_{2}+1, & \text{\emph{for\ \ }}i=1 \\ 
a_{i+1}, & \text{\emph{for\ \ }}2\leq i\leq \nu \ \left( =k-1\right) 
\end{array}
\right.  &  \\ 
&  & 
\end{array}
$}  \label{A2}
\end{equation}
Moreover, in case \emph{(\ref{PRW}) }we get\emph{:} 
\begin{equation}
\mathsf{P}_{i}=\mathsf{P}_{i+1}^{\vee },\ \ \mathsf{Q}_{i}=\mathsf{P}%
_{i+1}^{\vee }-\mathsf{Q}_{i+1}^{\vee },\ \forall i,\ 1\leq i\leq k\ \left(
=\nu -1\right)   \label{PQ1}
\end{equation}
while in case \emph{(\ref{DEU}):} 
\begin{equation}
\mathsf{P}_{i}=\mathsf{P}_{i-1}^{\vee },\ \ \mathsf{Q}_{i}=\mathsf{P}%
_{i-1}^{\vee }-\mathsf{Q}_{i-1}^{\vee },\ \forall i,\ 1\leq i\leq k\ \left(
=\nu +1\right)   \label{PQ2}
\end{equation}
\end{lemma}

\noindent \textit{Proof. }For $p=1$, we have 
\[
\dfrac{q}{q-1}=\ \left[ 1,q-1\right] =\stackunder{\left( q-1\right) \text{%
-times}}{\ \underbrace{\left[ \!\left[ 2,2,\ldots ,2,2\right] \!\right] }} 
\]
and the assertion is trivially true. From now on, suppose $p\geq 2$\textit{. 
}To pass from the ordered $t$-tuple $\left( b_{1}^{\vee },\ldots
,b_{t}^{\vee }\right) $ to the ordered $\rho $-tuple $\left( b_{1},\ldots
,b_{\rho }\right) $ we use the so-called Riemenschneider's \textit{point
diagram} (cf. \cite{Riem}, pp. 222-223)\textit{. }This is a plane
arrangement of points 
\[
\mathbf{PD}\subset \left\{ 1,\ldots ,t\right\} \times \left\{ 1,\ldots ,\rho
\right\} 
\]
with $t$ rows and $\rho $ columns, which is depicted as follows: 
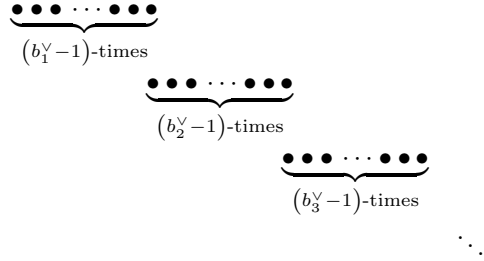
\[
\begin{array}{cc}
\stackunder{\left( b_{1}^{\vee }-1\right) \text{-times}}{\underbrace{\bullet
\bullet \bullet \,\cdots \bullet \bullet \,\bullet }} &  \\ 
\hspace{2.5cm}\hspace{0.5cm}\hspace{0.3cm}\hspace{0.3cm}\stackunder{\left(
b_{2}^{\vee }-1\right) \text{-times}}{\underbrace{\bullet \bullet \bullet
\,\cdots \bullet \bullet \,\bullet }} &  \\ 
\hspace{6cm}\hspace{0.5cm}\hspace{0.3cm}\hspace{0.4cm}\stackunder{\left(
b_{3}^{\vee }-1\right) \text{-times}}{\underbrace{\bullet \bullet \bullet
\,\cdots \bullet \bullet \,\bullet }} &  \\ 
& \mathbf{\ddots }
\end{array}
\]
and has the property that 
\[
\#\ \left\{ j\text{-th column points}\right\} =b_{j}-1,\ \ \forall j,\ \
1\leq j\leq \rho \ . 
\]
Separating the columns with exactly one entry, we consider the first column
whose number of points is $\geq 2$. (Such a column exists always for $p\geq
2 $). Obviously, this number of points has to be equal to 
\[
\#\ \left\{ 
\begin{array}{c}
\text{``twos'' occuring after the first position} \\ 
\text{in the ordered }t\text{-tuple\ }\left( b_{1}^{\vee },b_{2}^{\vee
},\ldots ,b_{t}^{\vee }\right)
\end{array}
\right\} =a_{2}+1\ . 
\]
Repeating the same procedure for the second, third, ... etc., column of $%
\mathbf{PD}$ whose numbers of points are $\geq 2$, and using the same
argument, we may rewrite the point diagram as: 
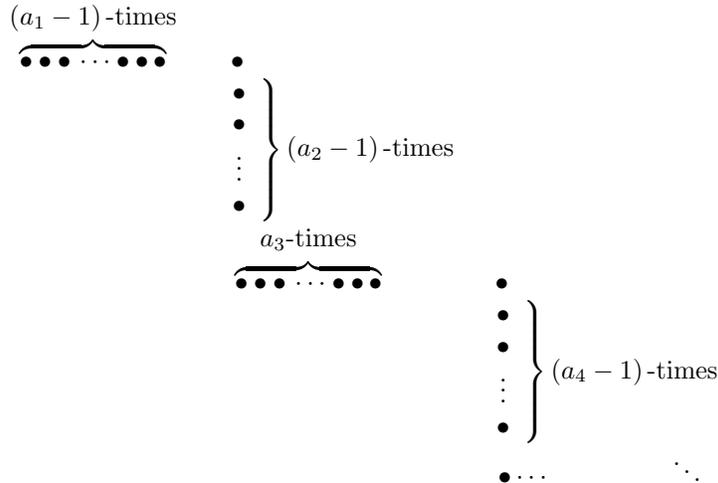
\[
\begin{array}[b]{lll}
\stackrel{
\begin{array}{c}
\left( a_{1}-1\right) \text{-times}
\end{array}
}{\overbrace{\bullet \bullet \bullet \,\cdots \bullet \bullet \,\bullet }} & 
\hspace{0.2cm}\bullet &  \\ 
& \left. 
\begin{array}{c}
\bullet \\ 
\bullet \\ 
\vdots \\ 
\bullet
\end{array}
\right\} \left( a_{2}-1\right) \text{-times} &  \\ 
& \hspace{0.25cm}\stackrel{
\begin{array}{c}
a_{3}\text{-times}
\end{array}
}{\overbrace{\bullet \bullet \bullet \,\cdots \bullet \bullet \,\bullet }} & 
\hspace{0.2cm}\bullet \\ 
&  & \left. 
\begin{array}{c}
\bullet \\ 
\bullet \\ 
\vdots \\ 
\bullet
\end{array}
\right\} \left( a_{4}-1\right) \text{-times} \\ 
&  & \hspace{0.2cm}\hspace{0.04cm}\bullet \cdots \hspace{0.8cm}\hspace{0.8cm}%
\ddots
\end{array}
\medskip 
\]
Inductively it is easy to show that 
\[
\left( b_{1},b_{2},\ldots ,b_{\rho }\right) = 
\]
\[
=\left\{ 
\begin{array}{ll}
\left( \stackunder{\left( a_{1}-1\right) \text{-times}}{\underbrace{2,...,2}}%
,a_{2}+2,\stackunder{\left( a_{3}-1\right) \text{-times}}{\underbrace{2,...,2%
}},a_{4}+2,\ldots ,a_{k-2}+2,\stackunder{\left( a_{k-1}-1\right) \text{-times%
}}{\underbrace{2,...,2}},a_{k}+1\right) , & \text{if\emph{\ \ }}k\ \text{
even} \\ 
\  &  \\ 
\left( \stackunder{\left( a_{1}-1\right) \text{-times}}{\underbrace{2,...,2}}%
,a_{2}+2,\stackunder{\left( a_{3}-1\right) \text{-times}}{\underbrace{2,...,2%
}},a_{4}+2,\ldots ,\stackunder{\left( a_{k-2}-1\right) \text{-times}}{%
\underbrace{2,...,2}},a_{k-1}+2,\stackunder{\left( a_{k}-1\right) \text{%
-times}}{\underbrace{2,...,2}}\right) , & \text{if\emph{\ \ }}k\ \text{ odd}
\end{array}
\right. \medskip \medskip 
\]
On the other hand,\medskip\ proposition \ref{PL-M} implies\ 
\[
\left( b_{1},b_{2},\ldots ,b_{\rho }\right) = 
\]
\[
=\left\{ 
\begin{array}{ll}
\left( a_{1}^{\vee }+1,\stackunder{\left( a_{2}^{\vee }-1\right) \text{-times%
}}{\underbrace{2,...,2}},a_{3}^{\vee }+2,\stackunder{\left( a_{4}^{\vee
}-1\right) \text{-times}}{\underbrace{2,...,2}},a_{5}^{\vee }+2,\ldots
,a_{\nu -1}^{\vee }+2,\stackunder{\left( a_{\nu }^{\vee }-1\right) \text{%
-times}}{\underbrace{2,...,2}}\right) , & \text{if\emph{\ \ }}\nu \ \text{
even} \\ 
\  &  \\ 
\left( a_{1}^{\vee }+1,\stackunder{\left( a_{2}^{\vee }-1\right) \text{-times%
}}{\underbrace{2,...,2}},a_{3}^{\vee }+2,\stackunder{\left( a_{4}^{\vee
}-1\right) \text{-times}}{\underbrace{2,...,2}},a_{5}^{\vee }+2,\ldots ,%
\stackunder{\left( a_{\nu -1}^{\vee }-1\right) \text{-times}}{\underbrace{%
2,...,2}}a_{\nu }^{\vee }+1\right) , & \text{if\ \emph{\ }}\nu \ \text{ odd}
\end{array}
\medskip \right. \smallskip 
\]
Since the ``twos'' must be placed at the same positions, we obtain (\ref
{MAKRIA}), (\ref{PRW}), (\ref{DEU}); moreover, since $a_{\nu }^{\vee }$ and $%
a_{k}$ are $\geq 2$, direct comparison gives (\ref{A1}), (\ref{A2}).
Finally, the equality (\ref{PQ1}) (resp. (\ref{PQ2})) is an immediate
consequence of (\ref{A1}) (resp. (\ref{A2})). $_{\Box }\bigskip $

\noindent We define now 
\[
d_{0}^{\vee }:=1,\ \text{ }d_{i}^{\vee }:=1+a_{2}^{\vee }+a_{4}^{\vee
}+\cdots +a_{2i}^{\vee },\ \ \forall i,\ \ 1\leq i\leq \QTOVERD\lfloor
\rfloor {\nu -1}{2}\ , 
\]
and 
\[
d_{0}:=1,\text{ \ }d_{i}:=1+a_{2}+a_{4}+\cdots +a_{2i},\ \ \forall i,\ \
1\leq i\leq \QTOVERD\lfloor \rfloor {k-1}{2}\ , 
\]
as well as 
\[
J:=\left\{ 
\begin{array}{cc}
\left\{ d_{1}^{\vee },d_{2}^{\vee },\ldots ,d_{\QOVERD\lfloor \rfloor {\nu
-1}{2}}^{\vee }\right\} \medskip , & \text{in case (\ref{PRW})} \\ 
\left\{ d_{0}^{\vee },d_{1}^{\vee },\ldots ,d_{\QOVERD\lfloor \rfloor {\nu
-1}{2}}^{\vee }\right\} , & \text{in case (\ref{DEU})}
\end{array}
\right. 
\]
and in the dual sense 
\[
J^{\vee }:=\left\{ 
\begin{array}{cc}
\left\{ d_{0},d_{2},\ldots ,d_{\QOVERD\lfloor \rfloor {k-1}{2}}\right\}
\medskip , & \text{in case (\ref{PRW})} \\ 
\left\{ d_{1},d_{1},\ldots ,d_{\QOVERD\lfloor \rfloor {k-1}{2}}\right\} , & 
\text{in case (\ref{DEU})}
\end{array}
\smallskip \right. 
\]
To formulate the main theorem of this section, let us further define 
\[
\Theta _{\sigma }:=\text{ conv}\left( \sigma \cap \left( N\Bbb{r}\left\{ 
\mathbf{0}\right\} \right) \right) \subset N_{\Bbb{R}},\ \ \ \text{resp.}\ \
\Theta _{\sigma ^{\vee }}:=\text{ conv}\left( \sigma ^{\vee }\cap \left( M%
\Bbb{r}\left\{ \mathbf{0}\right\} \right) \right) \subset M_{\Bbb{R}}, 
\]
denote by $\partial \Theta _{\sigma }^{\mathbf{cp}}$ (resp. $\partial \Theta
_{\sigma ^{\vee }}^{\mathbf{cp}}$) the part of the boundary polygon $%
\partial \Theta _{\sigma }$ (resp. $\partial \Theta _{\sigma ^{\vee }}$)
containing only its compact edges, and write vert$\left( \partial \Theta
_{\sigma }^{\mathbf{cp}}\right) $ (resp.$\ $vert$\left( \partial \Theta
_{\sigma ^{\vee }}^{\mathbf{cp}}\right) $) for the set of the vertices of $%
\partial \Theta _{\sigma }^{\mathbf{cp}}$ (resp. $\partial \Theta _{\sigma
^{\vee }}^{\mathbf{cp}}$).

\begin{theorem}[Determination of the vertices by Kleinian approximations]
\ \ \label{KAP}Consider lattice points $\left( \mathbf{u}_{j}\right) _{0\leq
j\leq \rho +1}$, $\left( \mathbf{v}_{i}\right) _{0\leq i\leq k+1}$ of $N$ $\ 
$and lattice points $\left( \mathbf{u}_{j}^{\vee }\right) _{0\leq j\leq t+1}$%
, $\left( \mathbf{v}_{i}^{\vee }\right) _{0\leq i\leq \nu +1}$ of $M$
defined by the vectorial recurrence relations\emph{:\smallskip } 
\begin{equation}
\left\{ 
\begin{array}{c}
\begin{array}{ll}
\mathbf{u}_{0}=n_{1}=\frak{y}_{1}, & \mathbf{u}_{1}=\frak{y}_{1}+\frak{y}%
_{2},
\end{array}
\medskip  \\ 
\mathbf{u}_{j}=b_{j-1}\ \mathbf{u}_{j-1}-\mathbf{u}_{i-2},\ \forall j,\ \
2\leq j\leq \rho +1;\ \left[ \mathbf{u}_{\rho +1}=n_{2}\right] 
\end{array}
\smallskip \right.   \label{VECTOR1}
\end{equation}
and\smallskip\ 
\begin{equation}
\left\{ 
\begin{array}{c}
\mathbf{u}_{0}^{\vee }=\frak{m}_{2},\ \ \ \mathbf{u}_{1}^{\vee }=\frak{m}%
_{1},\medskip  \\ 
\mathbf{u}_{j}^{\vee }=b_{j-1}^{\vee }\ \mathbf{u}_{j-1}^{\vee }-\mathbf{u}%
_{j-2}^{\vee },\ \forall j,\ \ 2\leq j\leq t+1;\ \left[ \mathbf{u}%
_{t+1}^{\vee }=\left( q-p\right) \frak{m}_{2}+q\left( \frak{m}_{1}-\frak{m}%
_{2}\right) \right] 
\end{array}
\right. \smallskip   \label{VECTOR2}
\end{equation}
and by the \emph{``Kleinian recurrence relations'' \cite{Klein2}, \cite
{Klein3}, \cite{Klein4}:\smallskip } 
\begin{equation}
\left\{ 
\begin{array}{ll}
\mathbf{v}_{i}=\mathsf{Q}_{i-1}\,\frak{y}_{1}+\mathsf{P}_{i-1}\,\frak{y}%
_{2},\  & \forall i,\ \ 0\leq i\leq k+1 \\ 
&  \\ 
\mathbf{v}_{i}^{\vee }=\mathsf{Q}_{i-1}^{\vee }\,\frak{m}_{2}+\mathsf{P}%
_{i-1}^{\vee }\,\left( \frak{m}_{1}-\frak{m}_{2}\right) ,\  & \forall i,\ \
0\leq i\leq \nu +1
\end{array}
\smallskip \right.   \label{Kleinian}
\end{equation}
respectively. Then 
\begin{equation}
\partial \Theta _{\sigma }^{\mathbf{cp}}\cap N=\left\{ \mathbf{u}_{j}\
\left| \ 0\leq j\leq \rho \right. +1\right\} ,\ \ \partial \Theta _{\sigma
^{\vee }}^{\mathbf{cp}}\cap M=\left\{ \mathbf{u}_{j}^{\vee }\ \left| \ 0\leq
j\leq t\right. +1\right\}   \label{TH1}
\end{equation}
and 
\begin{equation}
\text{\emph{vert}}\left( \partial \Theta _{\sigma }^{\mathbf{cp}}\right)
=\left\{ \mathbf{u}_{j}\ \left| \ j\in J\cup \left\{ 0,\rho +1\right\}
\right. \right\} ,\ \ \text{\emph{vert}}\left( \partial \Theta _{\sigma
^{\vee }}^{\mathbf{cp}}\right) =\left\{ \mathbf{u}_{j}^{\vee }\ \left| \
j\in J^{\vee }\cup \left\{ 0,t+1\right\} \right. \right\} \   \label{TH2}
\end{equation}
In particular, 
\begin{equation}
\fbox{$\ \forall i,\ 1\leq i\leq \QOVERD\lfloor \rfloor {k}{2}:\ \ \ \mathbf{%
v}_{2i}=\left\{ 
\begin{array}{cc}
\mathbf{u}_{d_{i}^{\vee }}\ , & \text{\emph{in case} \emph{(\ref{PRW})}} \\ 
\mathbf{u}_{d_{i-1}^{\vee }}\ , & \text{\emph{in case} \emph{(\ref{DEU})}}
\end{array}
\right. \ $}  \label{VVer}
\end{equation}
and 
\begin{equation}
\fbox{$\ \forall i,\ 1\leq i\leq \QOVERD\lfloor \rfloor {\nu }{2}:\ \ \ 
\mathbf{v}_{2i}^{\vee }=\left\{ 
\begin{array}{cc}
\mathbf{u}_{d_{i}}^{\vee }\ , & \text{\emph{in case (\ref{DEU})}} \\ 
\mathbf{u}_{d_{i-1}}^{\vee }\ , & \text{\emph{in case} \emph{(\ref{PRW})}}
\end{array}
\right. \ $}  \label{VVar}
\end{equation}
and therefore 
\begin{equation}
\fbox{$
\begin{array}{ccc}
&  &  \\ 
& \text{\emph{vert}}\left( \partial \Theta _{\sigma }^{\mathbf{cp}}\right)
=\left\{ \mathbf{v}_{0}\right\} \cup \left\{ \mathbf{v}_{2i}\ \left| \ 1\leq
i\leq \QOVERD\lfloor \rfloor {k}{2}\right. \right\} \cup \left\{ \mathbf{v}%
_{k+1}\right\}  &  \\ 
&  & 
\end{array}
$}  \label{VER1}
\end{equation}
and 
\begin{equation}
\fbox{$
\begin{array}{ccc}
&  &  \\ 
& \text{\emph{vert}}\left( \partial \Theta _{\sigma ^{\vee }}^{\mathbf{cp}%
}\right) =\left\{ \mathbf{v}_{0}^{\vee }\right\} \cup \left\{ \mathbf{v}%
_{2i}^{\vee }\ \left| \ 1\leq i\leq \QOVERD\lfloor \rfloor {\nu }{2}\right.
\right\} \cup \left\{ \mathbf{v}_{\nu +1}^{\vee }\right\}  &  \\ 
&  & 
\end{array}
$}\smallskip   \label{VER2}
\end{equation}
respectively. This means that the vertices of $\ \partial \Theta _{\sigma }^{%
\mathbf{cp}}$ and $\partial \Theta _{\sigma ^{\vee }}^{\mathbf{cp}}$ can be
read off just from the entries \emph{(}with \emph{even} and \emph{final }%
indices\emph{) }of the regular continued fraction expansions \emph{(\ref{Q-P}%
), }because \emph{(\ref{Kleinian}) }are uniquely determined by the vectorial
recurrence relations\emph{:\smallskip } 
\begin{equation}
\fbox{$
\begin{array}{c}
\\ 
\ \left\{ 
\begin{array}{l}
\medskip \mathbf{v}_{0}=n_{1}=\frak{y}_{1},\ \mathbf{v}_{1}=\frak{y}_{2},\ \
\ \ \ \,\mathbf{v}_{i}=a_{i-1}\ \mathbf{v}_{i-1}+\mathbf{v}_{i-2},\ \forall
i,\ \ 2\leq i\leq k+1 \\ 
\medskip \mathbf{v}_{0}^{\vee }=\frak{m}_{2},\ \mathbf{v}_{1}^{\vee }=\frak{m%
}_{1}-\frak{m}_{2},\ \ \mathbf{v}_{i}^{\vee }=a_{i-1}^{\vee }\ \mathbf{v}%
_{i-1}^{\vee }+\mathbf{v}_{i-2}^{\vee },\ \forall i,\ \ 2\leq i\leq \nu +1
\end{array}
\smallskip \right.  \\ 
\ \ 
\end{array}
$}  \label{Klein-rev}
\end{equation}
\end{theorem}

\noindent \textit{Proof}. At first notice that $\mathbf{v}_{0}=\mathbf{u}%
_{0}=n_{1}=\frak{y}_{1}$, $\mathbf{v}_{k+1}=\mathbf{u}_{\rho +1}=n_{2}$, and
that the lattice points $\mathbf{u}_{j}$ defined by (\ref{VECTOR1}) can be
determined by the vectorial matrix multiplication:\smallskip 
\[
\left( 
\begin{array}{c}
\mathbf{u}_{1} \\ 
\mathbf{u}_{2} \\ 
\vdots \\ 
\\ 
\vdots \\ 
\mathbf{u}_{\rho -1} \\ 
\mathbf{u}_{\rho }
\end{array}
\right) =\left( 
\begin{array}{cccccc}
b_{1} & -1 & 0 & \cdots & \cdots & 0 \\ 
-1 & b_{2} & -1 & \cdots & \cdots & 0 \\ 
0 & -1 & b_{3} & -1 & \cdots & 0 \\ 
\vdots & \vdots & \vdots & \vdots & \vdots & \vdots \\ 
&  &  &  &  &  \\ 
\vdots & \vdots & \vdots & \vdots & \vdots & \vdots \\ 
0 & \cdots & \cdots & 0 & -1 & b_{\rho }
\end{array}
\right) ^{-1}\ \left( 
\begin{array}{c}
n_{1} \\ 
0 \\ 
\vdots \\ 
\\ 
\vdots \\ 
0 \\ 
n_{2}
\end{array}
\right) \smallskip 
\]
(That this matrix is indeed invertible, is well-known from the theory of
continuants, cf. \cite{Perron}, I, \S 4. Furthermore, its determinant equals 
$q=\mathsf{R}_{\rho }$). Computing the corresponding adjoint matrix and
performing the multiplication, we obtain for all $j$, $0\leq j\leq \rho +1$, 
\begin{equation}
\mathbf{u}_{j}=\left( \frac{\widetilde{\mathsf{R}}_{j-1}}{q}\right)
\,n_{1}+\left( \frac{\mathsf{R}_{j-1}}{q}\right) \,n_{2}=\frac{1}{q}\left( 
\widetilde{\mathsf{R}}_{j-1}+p\,\mathsf{R}_{j-1}\right) \,\frak{y}_{1}+%
\mathsf{R}_{j-1}\,\frak{y}_{2}  \label{UJ1}
\end{equation}
where $\left( \widetilde{\mathsf{R}}_{j}\right) _{-1\leq j\leq \rho }$ is a
finite sequence of integers defined by the recurrence relations: 
\[
\widetilde{\mathsf{R}}_{-1}=q\text{, \ \ \ }\widetilde{\mathsf{R}}_{0}=q-p%
\text{, \ \ \ }\widetilde{\mathsf{R}}_{j}=b_{j}\,\widetilde{\mathsf{R}}%
_{j-1}-\widetilde{\mathsf{R}}_{j-2},\ \ \forall j,\ \ 1\leq j\leq \rho 
\]
(with $\widetilde{\mathsf{R}}_{\rho -2}=b_{\rho }$, $\widetilde{\mathsf{R}}%
_{\rho -1}=1$, $\widetilde{\mathsf{R}}_{\rho }=0$). If we define, in
addition, another finite sequence $\left( \widetilde{\widetilde{\mathsf{R}}}%
_{j}\right) _{-1\leq j\leq \rho }$ of integers via 
\[
\widetilde{\widetilde{\mathsf{R}}}_{-1}=\widetilde{\widetilde{\mathsf{R}}}%
_{0}=1\text{, \ \ \ }\widetilde{\widetilde{\mathsf{R}}}_{1}=b_{1}-1\text{, \
\ \ }\widetilde{\widetilde{\mathsf{R}}}_{j}=b_{j}\,\widetilde{\widetilde{%
\mathsf{R}}}_{j-1}-\widetilde{\widetilde{\mathsf{R}}}_{j-2},\ \ \forall j,\
\ 1\leq j\leq \rho 
\]
(with $\widetilde{\widetilde{\mathsf{R}}}_{\rho }=p$), then for all $j$, $%
0\leq j\leq \rho $, we get the equalities: 
\begin{equation}
\left\{ 
\begin{array}{l}
\widetilde{\mathsf{R}}_{j-1}\,\mathsf{R}_{j}-\widetilde{\mathsf{R}}_{j}\,%
\mathsf{R}_{j-1}=q\smallskip \\ 
\mathsf{R}_{j}\,\widetilde{\widetilde{\mathsf{R}}}_{j-1}-\mathsf{R}_{j-1}\,%
\widetilde{\widetilde{\mathsf{R}}}_{j}=1\smallskip \\ 
\widetilde{\mathsf{R}}_{j}+p\,\mathsf{R}_{j}-q\,\widetilde{\widetilde{%
\mathsf{R}}}_{j}=0
\end{array}
\right.  \label{EQUALITIES}
\end{equation}
On the other hand, using the notation introduced in the proof of \ref{PL-M},
we deduce: 
\[
\left( 
\begin{array}{ll}
\mathsf{R}_{j} & -\mathsf{R}_{j-1} \\ 
\widetilde{\widetilde{\mathsf{R}}}_{j} & -\widetilde{\widetilde{\mathsf{R}}}%
_{j-1}
\end{array}
\right) =\mathcal{W}\cdot \left( \dprod\limits_{s=1}^{j}\ \mathcal{M}%
^{-}\left( b_{s}\right) \right) \ \stackrel{\text{by (\ref{nsuc})}}{=}\ 
\mathcal{W}\cdot \left( 
\begin{array}{ll}
\mathsf{R}_{j} & -\mathsf{R}_{j-1} \\ 
\mathsf{S}_{j} & -\mathsf{S}_{j-1}
\end{array}
\right) \ . 
\]
Hence, $\widetilde{\widetilde{\mathsf{R}}}_{j}=\mathsf{R}_{j}-\mathsf{S}_{j}$%
, and for all $j$, $0\leq j\leq \rho +1$, (\ref{UJ1}), (\ref{EQUALITIES})
give: 
\begin{equation}
\mathbf{u}_{j}=\widetilde{\widetilde{\mathsf{R}}}_{j-1}\,\frak{y}_{1}+%
\mathsf{R}_{j-1}\,\frak{y}_{2}=\left( \mathsf{R}_{j-1}-\mathsf{S}%
_{j-1}\right) \,\frak{y}_{1}+\mathsf{R}_{j-1}\,\frak{y}_{2}  \label{UJ2}
\end{equation}
Clearly, $\mathbf{u}_{j}$'s belong to $\Theta _{\sigma }\cap N$ (cf. (\ref
{sdec})). That these points are exactly the lattice points of $\partial
\Theta _{\sigma }^{\mathbf{cp}}\cap N$ was proved by Oda in \cite{Oda},
lemma 1.20, pp. 25-26, by making use of the fact that $\left\{ \mathbf{u}%
_{0},\mathbf{u}_{1}\right\} $ form a $\Bbb{Z}$-basis of $N$ and induction on 
$j$. Since for all $j$, $2\leq j\leq \rho +1$, 
\[
\left( \mathbf{u}_{j-2},\ \mathbf{u}_{j-1},\ \mathbf{u}_{j}\ \ \ \text{are
collinear}\right) \Longleftrightarrow b_{j-1}=2, 
\]
the first equality in (\ref{TH2}) follows from proposition \ref{PL-M}. To
prove (\ref{VVer}), note that in case (\ref{PRW}) the equality (\ref{PQ1})
and the matrix-transfer rule of rem. \ref{matr-transf} imply 
\[
\mathsf{P}_{2i-1}=\mathsf{P}_{2i}^{\vee }=\mathsf{R}_{d_{i}^{\vee }-1}\text{%
, \ \ }\mathsf{Q}_{2i-1}=\mathsf{Q}_{2i}^{\vee }=\mathsf{R}_{d_{i}^{\vee
}-1}-\mathsf{S}_{d_{i}^{\vee }-1}\text{ } 
\]
for all $i$, $1\leq i\leq \QOVERD\lfloor \rfloor {k}{2}$. This means that 
\[
\mathbf{v}_{2i}=\mathsf{Q}_{2i-1}\,\frak{y}_{1}+\mathsf{P}_{2i-1}\,\frak{y}%
_{2}=\mathsf{Q}_{2i}^{\vee }\,\frak{y}_{1}+\mathsf{P}_{2i}^{\vee }\,\frak{y}%
_{2}=\left( \mathsf{R}_{d_{i}^{\vee }-1}-\mathsf{S}_{d_{i}^{\vee }-1}\right)
\,\frak{y}_{1}+\mathsf{R}_{d_{i}^{\vee }-1}\,\frak{y}_{2}\ \stackrel{\text{%
by (\ref{UJ2})}}{=}\ \mathbf{u}_{d_{i}^{\vee }}\ \text{.} 
\]
Case (\ref{DEU}) can be treated similarly. Finally, all assertions (\ref{TH1}%
), (\ref{TH2}), (\ref{VVar}), (\ref{VER2}) concerning the dual cone $\sigma
^{\vee }$ are proved by the same method (i.e., just by replacing $k$ and $%
\rho $ by $\nu $ and $t$, respectively, interchanging the ``checks'' in the
dual sense, and using the same arguments for Gen$\left( \sigma ^{\vee
}\right) $). $_{\Box }$

\begin{remark}
\label{AAP}\emph{(i) As in the case of }$\mathbf{u}_{j}$\emph{'s, one may
determine }$\mathbf{v}_{i}$\emph{'s in (\ref{Klein-rev})} \emph{by the
following vectorial matrix multiplication:\smallskip } 
\[
\left( 
\begin{array}{c}
\mathbf{v}_{1} \\ 
\mathbf{v}_{2} \\ 
\vdots  \\ 
\\ 
\vdots  \\ 
\mathbf{v}_{k-1} \\ 
\mathbf{v}_{k}
\end{array}
\right) =\left( 
\begin{array}{cccccc}
a_{1} & -1 & 0 & \cdots  & \cdots  & 0 \\ 
1 & a_{2} & -1 & \cdots  & \cdots  & 0 \\ 
0 & 1 & a_{3} & -1 & \cdots  & 0 \\ 
\vdots  & \vdots  & \vdots  & \vdots  & \vdots  & \vdots  \\ 
&  &  &  &  &  \\ 
\vdots  & \vdots  & \vdots  & \vdots  & \vdots  & \vdots  \\ 
0 & \cdots  & \cdots  & 0 & 1 & a_{k}
\end{array}
\right) ^{-1}\ \left( 
\begin{array}{c}
-n_{1} \\ 
0 \\ 
\vdots  \\ 
\\ 
\vdots  \\ 
0 \\ 
n_{2}
\end{array}
\right) \smallskip \smallskip 
\]
\emph{(ii) It should be particularly mentioned that} 
\[
\ \#\left( \text{\emph{vert}}\left( \partial \Theta _{\sigma }^{\mathbf{cp}%
}\right) \right) -\#\left( \text{\emph{vert}}\left( \partial \Theta _{\sigma
^{\vee }}^{\mathbf{cp}}\right) \right) =\QOVERD\lfloor \rfloor
{k}{2}-\QOVERD\lfloor \rfloor {\nu }{2}=\left\{ 
\begin{array}{lll}
1, & \text{\emph{if}} & \text{\ }\nu \notin 2\,\Bbb{Z}\text{ \ \emph{\&} \ }%
k=\nu +1 \\ 
\, &  &  \\ 
0, & \text{\emph{if}} & \left\{ 
\begin{array}{ll}
\text{\emph{either}} & \nu \notin 2\,\Bbb{Z}\text{ \ \emph{\&} \ }k=\nu -1
\\ 
\text{\emph{or}} & \nu \in 2\,\Bbb{Z}\text{ \ \emph{\&} \ }k=\nu +1
\end{array}
\right.  \\ 
\, &  &  \\ 
-1, & \text{\emph{if}} & \nu \in 2\,\Bbb{Z}\text{ \ \emph{\&} \ }k=\nu -1
\end{array}
\medskip \right. 
\]
\emph{(iii)} \emph{The }$\mathbf{v}_{i}$\emph{'s with odd indices (and final
indices }$1,k+1$\emph{) are the vertices of }$\partial \Theta _{\tau }^{%
\mathbf{cp}}$\emph{, i.e.,} 
\[
\text{\emph{vert}}\left( \partial \Theta _{\tau }^{\mathbf{cp}}\right)
=\left\{ \mathbf{v}_{1}\right\} \cup \left\{ \mathbf{v}_{2i+1}\ \left| \
1\leq i\leq \QOVERD\lfloor \rfloor {k-1}{2}\right. \right\} \cup \left\{ 
\mathbf{v}_{k+1}\right\} \ ,
\]
\emph{where }$\tau :=\emph{\ pos}\left( \frak{n}_{2},n_{2}\right) \subset N_{%
\Bbb{R}}$\emph{.\ Geometrically, the fact that }$\mathbf{v}_{2i+1}$\emph{'s
belong to }$\tau $\emph{\ follows directly from (\ref{Kleinian}) and the
inequalities (\ref{de-in}). Moreover, }$\tau $\emph{\ is a }$\left( \left[ q%
\right] _{p},p\right) $\emph{-cone w.r.t. }$\left\{ \frak{n}_{1},\frak{n}%
_{2}\right\} $\emph{. (See figure \textbf{1}}, \emph{where the illustrated
cone }$\sigma $\emph{\ is a }$\left( 4,7\right) $\emph{-cone and }$\tau $ 
\emph{is a} $\left( 3,4\right) $\emph{-cone).}\newline
\begin{figure}[htbp]
\begin{center}
\input{fig1.pstex_t}
\center{Figure \textbf{1}}
\end{center}
\end{figure}
\end{remark}

\section{The underlying spaces of abelian quotient singularities\\
as affine toric varieties\label{USP}}

\noindent Abelian quotient singularities can be directly investigated by
means of the theory of toric varieties. If $G$ is a finite subgroup of GL$%
\left( r,\Bbb{C}\right) $, then $\left( \Bbb{C}^{*}\right) ^r/G$ is
automatically an algebraic torus embedded in $\Bbb{C}^r/G$.\medskip

\noindent \textsf{(a) }Let $G$ be a finite subgroup of GL$\left( r,\Bbb{C}%
\right) $ which is \textit{small}, i.e. with no pseudoreflections, acting
linearly on $\Bbb{C}^{r}$, and $\mathsf{p}:$ $\Bbb{C}^{r}\rightarrow \Bbb{C}%
^{r}/G$ the quotient map. Denote by $\left( \Bbb{C}^{r}/G,\left[ \mathbf{0}%
\right] \right) $ the (germ of the) corresponding quotient singularity with $%
\left[ \mathbf{0}\right] :=\mathsf{p}\left( \mathbf{0}\right) $.

\begin{proposition}[Singular locus]
\label{SLOC}If $G$ is a small finite subgroup of \emph{GL}$\left( r,\Bbb{C}%
\right) $, then 
\[
\text{\emph{Sing}}\left( \Bbb{C}^{r}/G\right) =\mathsf{p}\left( \left\{ 
\mathbf{z}\in \Bbb{C}^{r}\ \left| \ G_{\mathbf{z}}\neq \left\{ \text{\emph{Id%
}}\right\} \right. \right\} \right) 
\]
where $G_{\mathbf{z}}:=\left\{ g\in G\ \left| \ g\cdot \mathbf{z=z}\right.
\right\} $ is the isotropy group of $\mathbf{z=}\left( z_{1},\ldots
,z_{r}\right) \in \Bbb{C}^{r}$.
\end{proposition}

\begin{theorem}[Prill's group-theoretic isomorphism criterion]
\label{prill-th}Let $G_{1}$, $G_{2}$ be two small finite subgroups of \emph{%
GL}$\left( r,\Bbb{C}\right) $. Then there exists an analytic isomorphism 
\[
\left( \Bbb{C}^{r}/G_{1},\left[ \mathbf{0}\right] \right) \cong \left( \Bbb{C%
}^{r}/G_{2},\left[ \mathbf{0}\right] \right) 
\]
if and only if $G_{1}$ and $G_{2}$ are conjugate to each other within \emph{%
GL}$\left( r,\Bbb{C}\right) $.
\end{theorem}

\noindent \textit{Proof. }See Prill \cite{Prill}, thm. 2, p. 382. $_{\Box }$
\medskip \smallskip

\noindent \textsf{(b) } Let $G$ be a finite, small, \textit{abelian}
subgroup of GL$\left( r,\Bbb{C}\right) $, $r\geq 2$, having order $l=\left|
G\right| \geq 2$, and let 
\[
\left\{ e_{1}=\left( 1,0,\ldots ,0,0\right) ^{\intercal },\ldots
,e_{r}=\left( 0,0,\ldots ,0,1\right) ^{\intercal }\right\} 
\]
denote the standard basis of $\Bbb{Z}^{r}$, $N_{0}:=\sum_{i=1}^{r}\Bbb{Z}%
e_{i}$ the standard rectangular lattice, $M_{0}$ its dual, and\smallskip 
\[
T_{N_{0}}:=\text{Max-Spec}\left( \Bbb{C}\left[ \frak{x}_{1}^{\pm 1},\ldots ,%
\frak{x}_{r}^{\pm 1}\right] \right) =\left( \Bbb{C}^{\ast }\right) ^{r}\,\,. 
\]
Clearly, 
\[
T_{N_{G}}:=\text{Max-Spec}\left( \Bbb{C}\left[ \frak{x}_{1}^{\pm 1},\ldots ,%
\frak{x}_{r}^{\pm 1}\right] ^{G}\right) =\left( \Bbb{C}^{\ast }\right)
^{r}/G 
\]
is an $r$-dimensional algebraic torus with $1$-parameter group $N_{G}$ and
with group of characters $M_{G}$. Using the exponential map 
\[
\left( N_{0}\right) _{\Bbb{R}}\ni \left( y_{1},\ldots ,y_{r}\right)
^{\intercal }=\mathbf{y\longmapsto }\text{ exp}\left( \mathbf{y}\right)
:=\left( e^{\left( 2\pi \sqrt{-1}\right) y_{1}},\ldots ,e^{\left( 2\pi \sqrt{%
-1}\right) y_{r}}\right) ^{\intercal }\in T_{N_{0}} 
\]
and the injection $\mathbf{\iota }:T_{N_{0}}\hookrightarrow $ GL$\left( r,%
\Bbb{C}\right) $ defined by 
\[
T_{N_{0}}\ni \left( t_{1},\ldots ,t_{r}\right) ^{\intercal }=\mathbf{t}%
\hookrightarrow \text{ }\mathbf{\iota }\left( \mathbf{t}\right) :=\text{diag}%
\left( t_{1},\ldots ,t_{r}\right) \in \text{ GL}\left( r,\Bbb{C}\right) \ , 
\]
we have obviously 
\[
N_{G}=\left( \mathbf{\iota \circ }\text{ exp}\right) ^{-1}\left( G\right) 
\text{ \ \ \ \ \ \ (\ and determinant \ det}\left( N_{G}\right) =\dfrac{1}{l}%
\ \text{)} 
\]
(as long as we make a concrete choice of eigencoordinates to diagonalize the
action of the elements of $G$ on $\Bbb{C}^{r}$) with\smallskip 
\[
M_{G}=\left\{ m\in M_{0}\,\left| 
\begin{array}{c}
\frak{x}^{m}=\,\frak{x}_{1}^{\mu _{1}}\,\cdots \,\frak{x}_{r}^{\mu _{r}}%
\text{ \ \ is a }G\text{-invariant } \\ 
\text{Laurent monomial (}m=\left( \mu _{1},\ldots ,\mu _{r}\right) \text{)}
\end{array}
\right. \right\} \,\,\text{(and det}\left( M_{G}\right) =l\text{)}. 
\]
\newline
$\bullet $ If we define 
\[
\fbox{$
\begin{array}{ccc}
& \sigma _{0}:=\ \text{pos}\left( \left\{ e_{1},..,e_{r}\right\} \right) & 
\end{array}
$} 
\]
to be the $r$-dimensional positive orthant, and $\Delta _{G}$ to be the
fan\smallskip 
\[
\Delta _{G}:=\left\{ \sigma _{0}\text{ together with its faces}\right\} 
\]
then by the exact sequence $0\rightarrow G\cong N_{G}/N_{0}\rightarrow
T_{N_{0}}\rightarrow T_{N_{G}}\rightarrow 0$ induced by the canonical
duality pairing 
\[
M_{0}/M_{G}\times N_{G}/N_{0}\rightarrow \Bbb{Q}/\Bbb{Z\hookrightarrow C}%
^{\ast }\smallskip 
\]
(cf. \cite{Fulton}, p. 34, and \cite{Oda}, pp. 22-23), we get as projection
map:\smallskip\ $\Bbb{C}^{r}=X\left( N_{0},\Delta _{G}\right) \rightarrow
X\left( N_{G},\Delta _{G}\right) $, where 
\[
\fbox{$X\left( N_{G},\Delta _{G}\right) =U_{\sigma _{0}}=\Bbb{C}^{r}/G=$
Max-Spec$\left( \Bbb{C}\left[ \frak{x}_{1},\ldots ,\frak{x}_{r}\right]
^{G}\right) \hookleftarrow T_{N_{G}}$} 
\]
$\bullet $ Formally, we identify $\left[ \mathbf{0}\right] $ with orb$\left(
\sigma _{0}\right) $. Moreover, using the notation introduced in \S \ref
{Prel}\textsf{(g)}, the singular locus of $X\left( N_{G},\Delta _{G}\right) $
can be written (by \ref{SLOC} and \ref{SMCR}) as the union 
\[
\text{Sing}\left( X\left( N_{G},\Delta _{G}\right) \right) =\text{ orb}%
\left( \sigma _{0}\right) \cup \left( \bigcup \left\{ \text{Max-Spec}\left( 
\Bbb{C}\left[ \overline{\sigma }_{0}^{\vee }\cap M_{G}\left( \tau \right) 
\right] \right) \,\left| 
\begin{array}{c}
\,\tau \precneqq \sigma _{0}\text{,\thinspace dim}\left( \tau \right) \geq
2\smallskip \text{ \thinspace } \\ 
\text{and \thinspace mult}\left( \tau ;N_{G}\right) \geq 2
\end{array}
\right. \right\} \right) \ . 
\]
$\bullet $ In particular, if the acting group $G$ is \textit{cyclic}, then,
fixing diagonalization of the action on $\Bbb{C}^{r}$, we may assume that $G$
is generated by the element 
\[
\text{diag}\left( \zeta _{l}^{\alpha _{1}}\,,\ldots ,\,\zeta _{l}^{\alpha
_{r}}\right) 
\]
(with $\zeta _{l}:=e^{\frac{2\pi \sqrt{-1}}{l}}$) for $r$ integers $\alpha
_{1},\ldots ,\alpha _{r}\in \left\{ 0,1,\ldots ,l-1\right\} $, at least $2$
of which are $\neq 0$. This $r$-tuple $\left( \alpha _{1},\ldots ,\alpha
_{r}\right) $ of \textit{weights} is unique only up to the usual conjugacy
relations (see \ref{ISOCYC} below), and $N_{G}$ is to be identified with the
so-called \textit{lattice of weights} 
\[
N_{G}=N_{0}+\Bbb{Z\ }\left( \frac{1}{l}\left( \alpha _{1},\ldots ,\alpha
_{r}\right) ^{\intercal }\ \right) \ \ \ \ 
\]
containing all lattice points representing the elements of 
\[
G=\left\{ \text{diag}\left( \zeta _{l}^{\left[ \lambda \alpha _{1}\right]
_{l}},\ldots ,\zeta _{l}^{\left[ \lambda \alpha _{r}\right] _{l}}\right) \
\left| \ \ \ \lambda \in \Bbb{Z},\ \ 0\leq \lambda \leq l-1\right. \right\}
\ . 
\]

\begin{definition}
\label{type}\emph{Under these conditions, we say that the quotient
singularity} $\left( X\left( N_{G},\Delta _{G}\right) ,\text{\emph{orb}}%
\left( \sigma _{0}\right) \right) $\emph{\ is} \textit{of type} 
\begin{equation}
\fbox{$
\begin{array}{ccc}
& \dfrac{1}{l}\left( \alpha _{1},\ldots ,\alpha _{r}\right)  & 
\end{array}
$}  \label{ttype}
\end{equation}
\end{definition}

\noindent $\bullet $ Note that, since $G$ is small, gcd$\left( l,\alpha
_{1},\ldots ,\widehat{\alpha _{i}},\ldots ,\alpha _{r}\right) =1$, for all $%
i $, $1\leq i\leq r$. (The symbol $\widehat{\alpha _{i}}$ means here that $%
\alpha _{i}$ is omitted.)

\begin{definition}
\emph{The} \textit{splitting codimension} \emph{of the underlying space }$%
U_{\sigma _{0}}=\Bbb{C}^{r}/G$\emph{\ of an abelian quotient singularity is
defined to be the number} 
\[
\text{\emph{splcod}}\left( U_{\sigma _{0}}\right) :=\text{\emph{max\ }}%
\left\{ \varkappa \in \left\{ 2,\ldots ,r\right\} \ \left| \ 
\begin{array}{c}
U_{\sigma _{0}}\cong U_{\tau }\times \Bbb{C}^{r-\varkappa },\ \text{\emph{%
s.t.}} \\ 
\tau \precneqq \sigma _{0}\text{\emph{, dim}}\left( \tau \right) =\varkappa 
\\ 
\text{\emph{and Sing}}\left( U_{\tau }\right) \neq \varnothing 
\end{array}
\right. \right\} 
\]
\emph{If splcod}$\left( U_{\sigma _{0}}\right) =r$,\emph{\ then orb}$\left(
\sigma _{0}\right) $\emph{\ is called an }\textit{msc-singularity}\emph{,
i.e., a singularity having the maximum splitting codimension.\smallskip\ }
\end{definition}

\begin{lemma}
\label{ISOL}\emph{(i)} A cyclic quotient singularity of type $\emph{(\ref
{ttype})}$ has splitting codimension $\varkappa \in \left\{ 2,\ldots
,r-1\right\} $ if and only if there exists an index-subset $\left\{ \nu
_{1},\nu _{2},\ldots ,\nu _{r-\varkappa }\right\} \subset \left\{ 1,\ldots
,r\right\} $, such that 
\[
\alpha _{\nu _{1}}=\alpha _{\nu _{2}}=\cdots =\alpha _{\nu _{r-\varkappa
}}=0\ ,
\]
which is, in addition, \emph{maximal} w.r.t. this property.\smallskip 
\newline
\emph{(ii)} A cyclic quotient msc-singularity of type $\emph{(\ref{ttype})}$
is isolated if and only if 
\[
\emph{gcd}\left( \alpha _{i},l\right) =1,\ \forall i,\ \ 1\leq i\leq r\ .
\]
\end{lemma}

\noindent \textit{Proof. }It is immediate by the way we let $G$ act on $\Bbb{%
C}^r$. $_{\Box }\medskip \smallskip $

\noindent \textsf{(c) }For two integers $l,\,r\geq 2$, we define 
\[
\Lambda \left( l;r\right) :=\left\{ \left( \alpha _1,..,\alpha _r\right) \in
\left\{ 0,1,2,..,l-1\right\} ^r\left| 
\begin{array}{c}
\ \text{gcd}\left( l,\alpha _1,..,\widehat{\alpha _i},..,\alpha _r\right) =1,
\\ 
\text{for all }i,\ \ 1\leq i\leq r
\end{array}
\right. \right\} 
\]
and for $\left( \left( \alpha _1,\ldots ,\alpha _r\right) ,\left( \alpha
_1^{\prime },\ldots ,\alpha _r^{\prime }\right) \right) \in \Lambda \left(
l;r\right) \times \Lambda \left( l;r\right) $ the relation\smallskip 
\[
\left( \alpha _1,\ldots ,\alpha _r\right) \backsim \left( \alpha _1^{\prime
},\ldots ,\alpha _r^{\prime }\right) :\Longleftrightarrow \left\{ 
\begin{array}{l}
\text{there exists a permutation } \\ 
\theta :\left\{ 1,\ldots ,r\right\} \rightarrow \left\{ 1,\ldots ,r\right\}
\\ 
\text{and an integer }\lambda \text{, }1\leq \lambda \leq l-1\text{, } \\ 
\text{with gcd}\left( \lambda ,l\right) =1\text{, such that} \\ 
\alpha _{\theta \left( i\right) }^{\prime }=\left[ \lambda \cdot \alpha _i%
\right] _l\text{, }\forall i,\ 1\leq i\leq r
\end{array}
\right\} \smallskip 
\]
It is easy to see that $\backsim $ is an equivalence relation on $\Lambda
\left( l;r\right) \times \Lambda \left( l;r\right) $.

\begin{corollary}[Isomorphism criterion for cyclic acting groups]
\label{ISOCYC} \ \smallskip \newline
Let $G$, $G^{\prime }$ be two small, cyclic finite subgroups of \emph{GL}$%
\left( r,\Bbb{C}\right) $ acting on $\Bbb{C}^{r}$, and let the corresponding
quotient singularities be of type $\frac{1}{l}\,\left( \alpha _{1},\ldots
,\alpha _{r}\right) $ and $\frac{1}{l^{\prime }}\,\left( \alpha _{1}^{\prime
},\ldots ,\alpha _{r}^{\prime }\right) $ respectively. Then there exists an
analytic \emph{(}torus-equivariant\emph{) }isomorphism 
\[
\left( X\left( N_{G},\Delta _{G}\right) ,\text{\emph{orb}}\left( \sigma
_{0}\right) \right) \cong \left( X\left( N_{G^{\prime }},\Delta _{G^{\prime
}}\right) ,\text{\emph{orb}}\left( \sigma _{0}\right) \right) 
\]
if and only if $l=l^{\prime }$ and $\left( \alpha _{1},\ldots ,\alpha
_{r}\right) \backsim \left( \alpha _{1}^{\prime },\ldots ,\alpha
_{r}^{\prime }\right) $ within $\Lambda \left( l;r\right) $.
\end{corollary}

\noindent \textit{Proof. }It follows from \ref{prill-th} (cf. Fujiki \cite
{Fujiki}, lemma 2, p. 296). $_{\Box }$

\begin{proposition}[Gorenstein-condition]
\label{Gor-prop} \ \smallskip \newline
Let $\left( \Bbb{C}^{r}/G,\left[ \mathbf{0}\right] \right) =\left( X\left(
N_{G},\Delta _{G}\right) ,\text{\emph{orb}}\left( \sigma _{0}\right) \right) 
$ be an abelian quotient singularity. Then the following conditions are
equivalent \emph{:\smallskip } \newline
\emph{(i)} $X\left( N_{G},\Delta _{G}\right) =U_{\sigma _{0}}=\Bbb{C}^{r}/G$
is Gorenstein,\smallskip \newline
\emph{(ii) }$G\subset $ \emph{SL}$\left( r,\Bbb{C}\right) $,\smallskip 
\newline
\emph{(iii)} $\left\langle \left( 1,1,\ldots .1,1\right) ,n\right\rangle
\geq 1$, for all $n,\,n\in \sigma _{0}\cap \left( N_{G}\smallsetminus
\left\{ \mathbf{0}\right\} \right) $,\smallskip \newline
\emph{(iv) }$\left( X\left( N_{G},\Delta _{G}\right) ,\text{\emph{orb}}%
\left( \sigma _{0}\right) \right) $ is a canonical singularity of index $1$%
.\smallskip \newline
In particular, if $\left( \Bbb{C}^{r}/G,\left[ \mathbf{0}\right] \right) $
is cyclic of type $\frac{1}{l}\,\left( \alpha _{1},\ldots ,\alpha
_{r}\right) $, then \emph{(i)-(iv)} are equivalent to 
\[
\sum_{j=1}^{r}\alpha _{j}\equiv 0\ \text{\emph{(mod} }l\text{\emph{)}}
\]
\end{proposition}

\noindent \textit{Proof. }See e.g. Reid \cite{Reid1}, thm. 3.1. $_{\Box
}\bigskip $

\noindent $\bullet $ If $X\left( N_{G},\Delta _{G}\right) $ is Gorenstein,
then the cone $\sigma _{0}=$ pos$\,\left( \frak{s}_{G}\right) $ is supported
by the so-called \textit{junior lattice simplex\smallskip } 
\[
\fbox{$
\begin{array}{ccc}
& \frak{s}_{G}=\text{conv}\left( \left\{ e_{1},..,e_{r}\right\} \right) & 
\end{array}
$} 
\]
(w.r.t. $N_{G}$, cf. \cite{Ito-Reid}, \cite{BD}). Note that up to $\mathbf{0}
$ there is no other lattice point of $\sigma _{0}\cap N_{G}$ lying ``under''
the affine hyperplane of $\Bbb{R}^{r}$ containing $\frak{s}_{G}$. Moreover,
the lattice points representing the $l-1$ non-trivial group elements are
exactly those belonging to the intersection of a dilation $\lambda \,\frak{s}%
_{G}$ of $\frak{s}_{G}$ with $\mathbf{Par}\left( \sigma _{0}\right) $, for
some integer $\lambda $, $1\leq \lambda \leq r-1$.

\section{Lattice triangulations, crepant projective resolutions and main
theorems\label{LATR}}

\noindent In this section we first briefly recall some general theorems
concerning the projective, crepant resolutions of Gorenstein abelian
quotient singularities in terms of appropriate lattice triangulations of the
junior simplex. (For detailed expositions we refer to \cite{DH}, \cite{DHZ1}%
, \cite{DHZ2}). After that we formulate our main theorems.\bigskip

\noindent \textsf{(a) }By vert$\left( \mathcal{S}\right) $ we denote the set
of vertices of a polyhedral complex $\mathcal{S}$. By a \textit{%
triangulation }$\mathcal{T}$ of a polyhedral complex $\mathcal{S}$ we mean a
geometric simplicial subdivision of $\mathcal{S}$ with vert$\left( \mathcal{S%
}\right) \subset $ vert$\left( \mathcal{T}\right) $. A polytope $P$ will be,
as usual, identified with the polyhedral complex consisting of $P$ itself
together with all its faces. If $\mathcal{S}_{1}$, $\mathcal{S}_{2}$ are two 
\textit{simplicial} complexes, then we denote by $\mathcal{S}_{1}\ast 
\mathcal{S}_{2}$ their \textit{join}.\bigskip\ 

\noindent \textsf{(b) }A triangulation $\mathcal{T}$ of an $r$-dimensional
polyhedral complex $\mathcal{S}$ is called \textit{coherent }(or \textit{%
regular}) if there exists a strictly upper convex $\mathcal{T}$-support
function $\psi :\left| \mathcal{T}\right| \rightarrow \Bbb{R}$, i.e., a
piecewise-linear real function defined on the underlying space $\left| 
\mathcal{T}\right| $ of $\mathcal{T}$, for which 
\[
\psi \left( \delta \ \mathbf{x}+\left( 1-\delta \right) \mathbf{\ y}\right)
\geq \delta \ \psi \left( \mathbf{x}\right) +\left( 1-\delta \right) \ \psi
\left( \mathbf{y}\right) ,\text{ for all \ }\mathbf{x},\mathbf{y}\in \left| 
\mathcal{T}\right| ,\ \text{and\ }\emph{\ }\delta \in \left[ 0,1\right] \ , 
\]
so that for every maximal simplex $\mathbf{s}$ of $\mathcal{T}$, there is a
linear function $\eta _{\mathbf{s}}:\left| \mathbf{s}\right| \mathbf{%
\rightarrow }\Bbb{R}$ satisfying $\psi \left( \mathbf{x}\right) \leq \eta _{%
\mathbf{s}}\left( \mathbf{x}\right) $, for all $\mathbf{x}\in $\emph{\ }$%
\left| \mathcal{T}\right| $, with equality being valid only for those $%
\mathbf{x}$ belonging to $\mathbf{s}$. The set of all strictly upper convex $%
\mathcal{T}$-support functions will be denoted by SUCSF$_{\Bbb{R}}\left( 
\mathcal{T}\right) $. \bigskip

\noindent \textsf{(c) }Let $N$ denote an $r$-dimensional lattice. By a 
\textit{lattice polytope }(w.r.t. $N$)\textit{\ }is meant a polytope in $N_{%
\Bbb{R}}\Bbb{\cong R}^r$ with vertices belonging to $N$. If $\left\{
n_0,n_1,\ldots ,n_k\right\} $ is a set of $k\leq r$ affinely independent
lattice points, $\mathbf{s}$ the lattice\emph{\ }$k$-dimensional simplex $%
\mathbf{s}=$ conv$\left( \left\{ n_0,n_1,n_2,\ldots ,n_k\right\} \right) $,
and $N_{\mathbf{s}}:=$ lin$\left( \left\{ n_1-n_0,\ldots ,n_k-n_0\right\}
\right) \cap N$, then\smallskip \newline
$\bullet $ we say that $\mathbf{s}$ is an \textit{elementary simplex}\emph{\ 
}if 
\[
\left\{ \mathbf{y}-n_0\ \left| \ \mathbf{y}\in \mathbf{s}\right. \right\}
\cap N_{\mathbf{s}}=\left\{ \mathbf{0},\,n_1-n_0,\ldots ,\,n_k-n_0\right\} . 
\]
\newline
$\bullet $ $\mathbf{s}$ is \textit{basic} if it has anyone of the following
equivalent properties:\smallskip \newline
(i) $\left\{ n_1-n_0,n_2-n_0,\ldots ,n_k-n_0\right\} $ is a $\Bbb{Z}$-basis
of $N_{\mathbf{s}}$,\smallskip \newline
(ii) $\mathbf{s\ }$has relative volume\ Vol$\left( \mathbf{s};N_{\mathbf{s}%
}\right) =\dfrac{\text{Vol}\left( \mathbf{s}\right) }{\text{det}\left( N_{%
\mathbf{s}}\right) }=\dfrac 1{k!}\ \ $(w.r.t. $N_{\mathbf{s}}$) .

\begin{lemma}
\label{EL-BA}\emph{(i)} Every basic lattice simplex is elementary.\newline
\emph{(ii)} Elementary lattice simplices of dimension $\leq 2$ are basic.
\end{lemma}

\noindent \textit{Proof. }See \cite{DH}, lemma 6.2. $_{\Box }$

\begin{example}
\label{Counter}\emph{The lattice }$r$\emph{-simplex} 
\[
\mathbf{s}=\text{\emph{\ conv}}\left( \left\{ \mathbf{0},e_{1},e_{2},\ldots
,e_{r-2},e_{r-1},\left( 1,1,\ldots ,1,1,k\right) ^{\intercal }\right\}
\right) \subset \Bbb{R}^{r}\text{\emph{,\thinspace } }r\geq 3\text{\emph{, }}%
k\geq 2\text{\emph{,}}
\]
\emph{(w.r.t. }$\Bbb{Z}^{r}$\emph{)} \emph{serves as example of an
elementary but non-basic simplex because }$\mathbf{s}\cap \Bbb{Z}^{r}=$ 
\emph{vert}$\left( \mathbf{s}\right) $ \emph{and} 
\[
r!\,\text{\emph{Vol}}\left( \mathbf{s};\Bbb{Z}^{r}\right) =\left| \text{%
\emph{det}}\left( e_{1},\ldots ,e_{r-1},\left( 1,1,\ldots ,1,1,k\right)
^{\intercal }\right) \right| =k\neq 1\ .
\]
\end{example}

\begin{definition}
\emph{A triangulation} $\mathcal{T}$ \emph{of a lattice polytope} $P\subset
N_{\Bbb{R}}\cong \Bbb{R}^{r}$ \emph{(w.r.t.} $N$\emph{)} \emph{is called} 
\textit{lattice triangulation }\emph{if vert}$\left( P\right) \subset $ 
\emph{vert}$\left( \mathcal{T}\right) \subset N$\emph{. The set of all
lattice triangulations of a lattice polytope }$P$ \emph{(w.r.t.} $N$\emph{)
will be denoted by} $\mathbf{LTR}_{N}\left( P\right) $.
\end{definition}

\begin{definition}
\emph{A lattice triangulation} $\mathcal{T}$ \emph{of }$P\subset N_{\Bbb{R}%
}\cong \Bbb{R}^{r}$ \emph{(w.r.t.} $N$\emph{)} \emph{is called }\textit{%
maximal triangulation }\emph{if vert}$\left( \mathcal{T}\right) =N\cap P$%
\emph{. A lattice triangulation} $\mathcal{T}$ \emph{of }$P$ \emph{is
obviously maximal if and only if each simplex }$\mathbf{s}$ \emph{of }$%
\mathcal{T}$ \emph{is elementary. A lattice triangulation} $\mathcal{T}$ 
\emph{of }$P$ \emph{is said to be }\textit{basic}\emph{\ if }$\mathcal{T}$ \ 
\emph{consists of exclusively basic simplices. We define :} 
\[
\begin{array}{l}
\mathbf{LTR}_{N}^{\emph{\scriptsize max}}\,\left( P\right) :=\left\{ 
\mathcal{T}\in \mathbf{LTR}_{N}\left( P\right) \ \left| \ \mathcal{T}\text{%
\emph{\ \ is a maximal triangulation of \ }}P\right. \right\} \smallskip \
,\smallskip  \\ 
\mathbf{LTR}_{N}^{\emph{\scriptsize basic}}\left( P\right) :=\left\{ 
\mathcal{T}\in \mathbf{LTR}_{N}^{\emph{\scriptsize max}}\,\left( P\right) \
\left| \ \mathcal{T}\text{\emph{\ \ is a basic triangulation of \ }}P\right.
\right\} \ .
\end{array}
\]
\emph{(Moreover, adding the prefix} $\mathbf{Coh}$\emph{-} \emph{to anyone
of the above sets, we shall mean the subsets of their elements which are
coherent). The }\textit{hierarchy of lattice triangulations }\emph{of a }$P$%
\emph{\ (as above) is given by the following inclusion-diagram:}
\end{definition}

\[
\begin{array}{ccccc}
\mathbf{LTR}_N^{\text{{\scriptsize basic}}}\left( P\right) & \subset & 
\mathbf{LTR}_N^{\text{{\scriptsize max}}}\,\left( P\right) & \subset & 
\mathbf{LTR}_N\left( P\right) \smallskip \\ 
\bigcup &  & \bigcup &  & \bigcup \smallskip \\ 
\mathbf{Coh}\text{-}\mathbf{LTR}_N^{\text{{\scriptsize basic}}}\left(
P\right) & \subset & \mathbf{Coh}\text{-}\mathbf{LTR}_N^{\text{{\scriptsize %
max}}}\,\left( P\right) & \subset & \mathbf{Coh}\text{-}\mathbf{LTR}_N\left(
P\right)
\end{array}
\]

\begin{proposition}
\label{NON-EM}For any lattice polytope $P\subset N_{\Bbb{R}}\cong \Bbb{R}^{r}
$ \emph{(}w.r.t. $N$\emph{)} the set of maximal coherent triangulations\emph{%
\ }$
\begin{array}{ccc}
\mathbf{Coh}\text{-}\mathbf{LTR}_{N}^{\emph{max}}\emph{\,}\left( P\right)  & 
\subset  & \mathbf{Coh}\text{-}\mathbf{LTR}_{N}\left( P\right) 
\end{array}
$ of $P$ is non-empty.
\end{proposition}

\noindent \textit{Proof. }Consider the s.c.p.cone supported by $P$ in $N_{%
\Bbb{R}}\oplus \Bbb{R\cong R}^{r+1}$, and then use \cite{OP}, cor. 3.8, p.
394. $_{\Box }$\emph{\ }

\begin{remark}[``Pathologies'']
\emph{(i) Already for }$r=2$\emph{\ there exist lots of examples of }$P$%
\emph{'s admitting basic, non-coherent triangulations (``whirlpool
phenomenon''). \smallskip \newline
(ii) A more remarkable pathological counterexample which was constructed
recently by Hibi and Ohsugi \cite{Hi-O} is a }$9$\emph{-dimensional }$0/1$%
\emph{-polytope (with }$15$ \emph{vertices) which possesses basic
triangulations, but none of whose coherent triangulations is basic. Hence,
for high-dimensional }$P$\emph{'s, }$\mathbf{LTR}_{N}^{\emph{basic}}\emph{\,}%
\left( P\right) \neq \varnothing $ \emph{does not necessarily imply }$%
\mathbf{Coh}$-$\mathbf{LTR}_{N}^{\emph{basic}}\emph{\,}\left( P\right) \neq
\varnothing $\emph{. (To the best of our knowledge, it is not as yet clear
if there exists any counterexample of this kind when we restrict ourselves
to the class of }\textit{lattice}\emph{\ }\textit{simplices }\emph{or
not.)\medskip }
\end{remark}

\noindent \textsf{(d) }To pass from triangulations to desingularizations we
need to introduce some extra notation.\smallskip

\begin{definition}
\emph{Let }$\left( X\left( N_{G},\Delta _{G}\right) ,\text{\emph{orb}}\left(
\sigma _{0}\right) \right) $ \emph{be an} $r$\emph{-dimensional abelian
Gorenstein quotient singularity} \emph{(}$r\geq 2$\emph{), and }$\frak{s}_{G}
$ \emph{the }$\left( r-1\right) $\emph{-dimensional junior simplex. For any
simplex }$\mathbf{s}$ \emph{of a lattice triangulation $\mathcal{T}$ of }$%
\frak{s}_{G}$ \emph{let }$\sigma _{\mathbf{s}}$ \emph{denote the s.c.p. cone}
\[
\sigma _{\mathbf{s}}:=\left\{ \lambda \,\mathbf{y}\in \left( N_{G}\right) _{%
\Bbb{R}}\ \left| \ \lambda \in \Bbb{R}_{\geq 0}\,,\ \mathbf{y}\in \mathbf{s}%
\right. \right\} \ \,\left( \,=\text{\emph{pos}}\left( \mathbf{s}\right) \,\,%
\text{\emph{within} \thinspace }\left( N_{G}\right) _{\Bbb{R}}\right) \ 
\]
\emph{supporting} $\mathbf{s}$. \emph{We define the fan } 
\[
\widehat{\Delta _{G}}\left( \mathcal{T}\right) :=\left\{ \sigma _{\mathbf{s}%
}\ \left| \ \mathbf{s}\in \mathcal{T}\right. \right\} 
\]
\emph{of s.c.p. cones in }$\left( N_{G}\right) _{\Bbb{R}}\cong \Bbb{R}^{r}$%
\emph{,} \emph{and } 
\[
\begin{array}{ll}
\mathbf{PCDES}\left( X\left( N_{G},\Delta _{G}\right) \right) := & \left\{ 
\begin{array}{c}
\text{\textit{partial crepant} }T_{N_{G}}\text{\emph{-equivariant }} \\ 
\text{\emph{desingularizations} \emph{of} \ }X\left( N_{G},\Delta
_{G}\right)  \\ 
\begin{array}{c}
\ \text{\emph{with overlying spaces having }} \\ 
\text{\emph{\ at most (}}\Bbb{Q}\text{\emph{-factorial) }\textit{canonical} }
\\ 
\text{\emph{singularities (of index }}1\text{\emph{)}}
\end{array}
\end{array}
\right\} \smallskip \ ,\medskip \ \smallskip \medskip  \\ 
\mathbf{PCDES}^{\emph{max}}\left( X\left( N_{G},\Delta _{G}\right) \right) :=
& \left\{ 
\begin{array}{c}
\text{\textit{partial crepant} }T_{N_{G}}\text{\emph{-equivariant }} \\ 
\text{\emph{desingularizations} \emph{of} \ }X\left( N_{G},\Delta
_{G}\right)  \\ 
\begin{array}{c}
\ \text{\emph{with overlying spaces having }} \\ 
\text{\emph{\ at most (}}\Bbb{Q}\text{\emph{-factorial) }\textit{terminal} }
\\ 
\text{\emph{singularities (of index }}1\text{\emph{)}}
\end{array}
\end{array}
\right\} \smallskip \ ,\medskip  \\ 
\mathbf{CDES}\left( X\left( N_{G},\Delta _{G}\right) \right) := & \left\{ 
\begin{array}{c}
\text{\textit{crepant} }T_{N_{G}}\text{\emph{-equivariant (full)}} \\ 
\text{\emph{\ desingularizations of} }X\left( N_{G},\Delta _{G}\right) 
\end{array}
\right\} \ .\medskip \smallskip 
\end{array}
\]
\emph{(Setting the prefix} $\mathbf{QP}$\emph{- in the front of anyone of
them, we shall mean the corresponding subsets of them consisting of those} 
\emph{desingularizations whose overlying spaces are }\textit{quasiprojective}%
\emph{.) }
\end{definition}

\begin{theorem}[Desingularizing by triangulations]
\label{TRCOR} \ \smallskip \newline
Let $\left( X\left( N_{G},\Delta _{G}\right) ,\text{\emph{orb}}\left( \sigma
_{0}\right) \right) $ be an $r$-dimensional abelian Gorenstein quotient
singularity \emph{(}$r\geq 2$\emph{).} Then there exist one-to-one
correspondences \emph{:\smallskip } 
\[
\fbox{$
\begin{array}{ccc}
&  &  \\ 
& 
\begin{array}{ccc}
\left( \mathbf{Coh}\text{-}\right) \mathbf{LTR}_{N_{G}}^{\emph{basic}}\left( 
\frak{s}_{G}\right) \smallskip  & \stackrel{\text{\emph{1:1}{\scriptsize %
\smallskip }}}{\longleftrightarrow } & \left( \mathbf{QP}\text{-}\right) 
\mathbf{CDES}\left( X\left( N_{G},\Delta _{G}\right) \right)  \\ 
\cap  &  & \cap  \\ 
\left( \mathbf{Coh}\text{-}\right) \mathbf{LTR}_{N_{G}}^{\emph{max}}\,\left( 
\frak{s}_{G}\right) \smallskip  & \stackrel{\text{\emph{1:1}{\scriptsize %
\smallskip }}}{\longleftrightarrow } & \left( \mathbf{QP}\text{-}\right) 
\mathbf{PCDES}^{\emph{max}}\left( X\left( N_{G},\Delta _{G}\right) \right) 
\\ 
\cap  &  & \cap  \\ 
\left( \mathbf{Coh}\text{-}\right) \mathbf{LTR}_{N_{G}}\left( \frak{s}%
_{G}\right)  & \stackrel{\text{\emph{1:1}{\scriptsize \smallskip }}}{%
\longleftrightarrow } & \left( \mathbf{QP}\text{-}\right) \mathbf{PCDES}%
\left( X\left( N_{G},\Delta _{G}\right) \right) 
\end{array}
&  \\ 
&  & 
\end{array}
$}\smallskip 
\]
which are realized by crepant $T_{N_{G}}$-equivariant birational morphism of
the form\smallskip \smallskip 
\begin{equation}
f_{\mathcal{T}}=\text{\emph{id}}_{\ast }:X\left( N_{G},\text{ }\widehat{%
\Delta _{G}}\left( \mathcal{T}\right) \right) \longrightarrow X\left(
N_{G},\Delta _{G}\right)   \label{DESING}
\end{equation}
induced by mapping 
\[
\mathcal{T}\longmapsto \widehat{\Delta _{G}}\left( \mathcal{T}\right) ,\ \ \
\ \ \widehat{\Delta _{G}}\left( \mathcal{T}\right) \longmapsto X\left( N_{G},%
\text{ }\widehat{\Delta _{G}}\left( \mathcal{T}\right) \right) \ .
\]
\end{theorem}

\noindent \textit{Proof. }See \cite{DHZ1}, \S\ 4 and \cite{DH}, thm. 6.9. $%
_{\Box }\bigskip $ \ \smallskip

\noindent \textsf{(e) }It is now clear by theorem \ref{TRCOR} that the main
question (formulated in \S \ref{INTRO}), restricted to the category of
torus-equivariant desingularizations of $X\left( N_{G},\Delta _{G}\right) $%
's, is equivalent to the following:\smallskip \newline
$\bullet $ \underline{\textbf{Question}}\textbf{: For a Gorenstein abelian
quotient singularity }$\left( X\left( N_{G},\Delta _{G}\right) ,\text{orb}%
\left( \sigma _{0}\right) \right) $\textbf{\ with junior simplex} $\frak{s}%
_{G}$\textbf{, under which conditions do we} \textbf{have }$\mathbf{Coh}$-$%
\mathbf{LTR}_{N_{G}}\left( \frak{s}_{G}\right) \neq \varnothing $ \textbf{%
?\smallskip }\newline
We shall answer this question in full generality for $2$-parameter
Gorenstein cyclic quotient singularitites below in thm. \ref{MT1}. (Further
techniques being applicable to other families of singularities will be
discussed in \cite{DHZ2}). Our starting-point is the following simple, but
very useful necessary criterion for the existence of basic triangulations of 
$\frak{s}_{G}$.

\begin{theorem}[Necessary Existence-Criterion]
\label{KILLER} \ \smallskip \newline
Let\emph{\ }$\left( X\left( N_{G},\Delta _{G}\right) ,\emph{orb}\left(
\sigma _{0}\right) \right) $ be a Gorenstein abelian quotient singularity.
If $\frak{s}_{G}$ admits a basic triangulation $\mathcal{T}$, then 
\begin{equation}
\fbox{$
\begin{array}{ccc}
&  &  \\ 
& \mathbf{Hlb}_{N_{G}}\left( \sigma _{0}\right) =\frak{s}_{G}\cap N_{G} & 
\\ 
&  & 
\end{array}
$}  \label{HILBCON}
\end{equation}
\end{theorem}

\noindent \textit{Proof}. See \cite{DH} thm. 6.15. $_{\Box }$

\begin{remark}
\emph{\label{FIRLA}(i) The well-known counterexamples (due to Bouvier and
Gonzalez-Sprinberg \cite{B-GS2}) of the two Gorenstein cyclic singularities
proving that theorem \ref{Seb} cannot be in general true for }$r\geq 4$ 
\emph{(i.e., in our terminology, those with types }$\frac{1}{7}\,\left(
1,3,4,6\right) $\emph{\ and }$\frac{1}{16}\,\left( 1,7,11,13\right) $\emph{%
), indicated the first technical difficulties for desigularizing in higher
dimensions (i.e., by }\textit{only}\emph{\ using the members of the Hilbert
basis of the corresponding monoid as the set of the minimal
cone-generators). Nevertheless, since they are both }\textit{terminal
singularities}\emph{, i.e., they do not possess junior elements up to the
vertices of }$\frak{s}_{G}$\emph{,} \emph{they} \textit{cannot} \emph{serve
as counterexamples to the inverse implication of thm. \ref{KILLER} (because }%
$\frak{s}_{G}$ \emph{does not admit basic triangulations by construction and}
$\frak{s}_{G}\cap N_{G}\subsetneqq \mathbf{Hlb}_{N_{G}}\left( \sigma
_{0}\right) $\emph{); in fact,} \emph{for a long time it was completely
unknown if condition (\ref{HILBCON}) might be sufficient or not for the
existence of a basic triangulation of }$\frak{s}_{G}$.\emph{\ Only recently
Firla and Ziegler (\cite{Firla}, \S 4.2 \& \cite{Fir-Zi}) discovered (by
computer testing) }$10$\emph{\ appropriate }\textit{counterexamples in
dimension} $4!$ \emph{Among them, the counterexample of the Gorenstein
cyclic quotient singularity with the smallest possible acting group-order,
fulfilling property (\ref{HILBCON}) and admitting no crepant,
torus-equivariant resolutions, is that of type }$\frac{1}{39}\,\left(
1,5,8,25\right) $\emph{.\smallskip }\newline
\emph{(ii) Our intention in the present paper is to show that the property
of all the above mentioned counterexamples to have }$3$ \emph{or more
``parameters'' (i.e., freely chosen weights) is not a pure chance! As we
shall see in theorems \ref{1PAR}, \ref{MT1} and \ref{MT2}, for }$1$\emph{-
and }$2$\emph{-parameter series of Gorenstein cyclic quotient singularities
condition (\ref{HILBCON}) \textbf{is}} \emph{indeed sufficient for the
existence of crepant, }$T_{N_{G}}$\emph{-equivariant, full resolutions in 
\textbf{all} dimensions. }
\end{remark}

\noindent Let us first recall what happens in the $1$-parameter case. (By
lemma \ref{ISOL} (i) it is enough to consider only msc-singularities.
Otherwise the problem can be reduced to a lower-dimensional one).

\begin{theorem}[On $1$-parameter singularity series]
\label{1PAR} \ \ \newline
Let $\left( X\left( N_{G},\Delta _{G}\right) ,\emph{orb}\left( \sigma
_{0}\right) \right) $ be an $r$-dimensional Gorenstein cyclic quotient
msc-singularity \emph{(}with $l=\left| G\right| \geq r\geq 4$\emph{)}.
Suppose that its type contains at least $r-1$ equal weights. Then $\left(
X\left( N_{G},\Delta _{G}\right) ,\emph{orb}\left( \sigma _{0}\right)
\right) $ is analytically isomorphic to Gorenstein cyclic quotient $1$%
-parameter singularity of type\smallskip 
\begin{equation}
\dfrac{1}{l}\ \left( \stackunder{\left( r-1\right) \text{\emph{-times}}}{%
\underbrace{1,1,\ldots ,1,1}},\ l-\left( r-1\right) \right)   \label{type1p}
\end{equation}
and $X\left( N_{G},\Delta _{G}\right) $ admits a \textbf{unique} crepant,
projective, $T_{N_{G}}$-equivariant, full desingularization 
\[
f=\emph{id}_{\ast }:X\left( N_{G},\widehat{\Delta _{G}}\right) \rightarrow
X\left( N_{G},\Delta _{G}\right) 
\]
\textbf{iff }\emph{either }$l\equiv 0$\emph{\ mod}$\left( r-1\right) $\emph{%
\ or }$l\equiv 1$\emph{\ mod}$\left( r-1\right) $. \emph{(}These conditions
are actually equivalent to condition \emph{(\ref{HILBCON}) }of \emph{thm. 
\ref{KILLER}).\smallskip }\newline
Moreover, the dimensions of the non-trivial cohomology groups of $X\left(
N_{G},\widehat{\Delta }_{G}\right) $ are given by the formulae\emph{:} 
\begin{equation}
\text{\emph{dim}}_{\Bbb{Q}}H^{2i}\left( X\left( N_{G},\widehat{\Delta }%
_{G}\right) ;\Bbb{Q}\right) =\left\{ 
\begin{array}{llll}
1 & , & \text{\emph{for }} & i=0 \\ 
\  &  &  &  \\ 
\QTOVERD\lfloor \rfloor {l}{r-1} & , & \text{\emph{for}} & i\in \left\{
1,2,\ldots ,r-2\right\}  \\ 
\  &  &  &  \\ 
\QTOVERD\lfloor \rfloor {l-1}{r-1} & , & \text{\emph{for}} & i=r-1
\end{array}
\right.   \label{cohd1}
\end{equation}
\end{theorem}

\noindent \textit{Proof. }See Dais-Henk \cite{DH}, thm. 8.2. $_{\Box }$

\begin{remark}
\emph{In fact, there are exactly }$\QTOVERD\lfloor \rfloor {l}{r-1}$\emph{\
exceptional prime divisors supported by the preimage of Sing}$\left( X\left(
N_{G},\Delta _{G}\right) \right) $ \emph{via} $f$\emph{; }$\QTOVERD\lfloor
\rfloor {l}{r-1}-1$\emph{\ of them are analytically isomorphic to the
projectivization of certain decomposable bundles over }$\Bbb{P}_{\Bbb{C}%
}^{r-2}$\emph{\ (having only twisted hyperplane bundles as summands), and
the last one is isomorphic }$\Bbb{P}_{\Bbb{C}}^{r-1}$\emph{\ (resp. to }$%
\Bbb{P}_{\Bbb{C}}^{r-2}\times \Bbb{C}$\emph{) for }$l\equiv 1$\emph{\ mod}$%
\left( r-1\right) $\emph{\ (resp. }$l\equiv 0$\emph{\ mod}$\left( r-1\right) 
$\emph{) [cf. \cite{DH}, thm. 8.4].}
\end{remark}

\begin{theorem}[Main Theorem I: On $2$-parameter singularity series]
\label{MT1}\smallskip\ \ \ \newline
Let $\left( X\left( N_{G},\Delta _{G}\right) ,\emph{orb}\left( \sigma
_{0}\right) \right) $ be a Gorenstein cyclic quotient msc-singularity of
type $\dfrac{1}{l}\left( \alpha _{1},\ldots ,\alpha _{r}\right) $ with $%
l=\left| G\right| \geq r\geq 4$\emph{, } for which at least $r-2$ of its
defining weights are equal. Then $X\left( N_{G},\Delta _{G}\right) $ admits
crepant, $T_{N_{G}}$-equivariant, full desingularizations \textbf{if and
only if} condition \emph{(\ref{HILBCON}) }is satisfied. Moreover, at least
one of these desingularizations is projective.
\end{theorem}

\begin{remark}
\emph{To examine the validity of condition (\ref{HILBCON}) in practice, one
has first to determine all the elements of the Hilbert basis }$\mathbf{Hlb}%
_{N_{G}}\left( \sigma _{0}\right) $\emph{\ and then to test if all of them
belong to the junior simplex or not.\smallskip\ }\allowbreak \newline
\emph{On the other hand, there is another, more direct method for working
with condition (\ref{HILBCON}); namely to translate the geometric properties
of }$\frak{s}_{G}$ \emph{into number-theoretic conditions fulfilled by the
weights of the defining ``type'' of the singularity. In the next theorem we
apply this method in the case in which }$r-2$ \emph{weights} \emph{are equal
to }$1$\emph{\ and the sum of all weights equals }$l$\emph{. As it turns
out, since all ``pure'' junior elements of }$s_{G}$\emph{\ have to be
coplanar, these number-theoretic conditions involve only linear congruences,
restrictions for certain gcd's and regular continued fraction expansions.
Hence, we have just to perform the euclidean division algorithm, 
                                %finitely many time 
 which is a polynomial-time-procedure.\smallskip\ }%
\newline
$\bullet $ \emph{Up to a very first step (involving the determination of
suitable Hermitian normal forms), and up to some extra numerological
conditions, we do not lose so much in generality (at least from the
algorithmic point of view) by restricting ourselves to the study of cyclic
singularities of type (\ref{typab}) (see rem. \ref{REDUCTION}
below).\smallskip }\newline
$\bullet $ \emph{Note that the problem of \ computing the elements of \ the
Hilbert basis of a }\textit{general}\emph{\ pointed rational cone is ``%
\textit{NP-hard}'' (cf. Henk-Weismantel \cite{HEW}, \S\ 3).\smallskip }%
\newline
$\bullet $ \emph{The formulae giving the dimensions of the non-trivial
cohomology groups of the fully resolvable }$2$\emph{-parameter singularities
are much more complicated than (\ref{cohd1}) and will be treated separately
by means of Ehrhart polynomials in section \ref{COHOM}.}
\end{remark}

\begin{theorem}[Main Theorem II]
\label{MT2}Let $\left( X\left( N_{G},\Delta _{G}\right) ,\emph{orb}\left(
\sigma _{0}\right) \right) $ be a Gorenstein cyclic quotient singularity of
type\smallskip\ 
\begin{equation}
\fbox{$
\begin{array}{ccc}
&  &  \\ 
& \dfrac{1}{l}\ \left( \stackunder{\left( r-2\right) \text{\emph{-times}}}{%
\underbrace{1,1,\ldots ,1,1}},\ \alpha ,\beta \right) ,\ \ \ l\geq r\geq 4,\
\ \alpha ,\beta \in \Bbb{N},\ \ \alpha +\beta =l-\left( r-2\right)  &  \\ 
&  & 
\end{array}
$}\smallskip   \label{typab}
\end{equation}
$\bullet $\emph{\ }Define 
\begin{equation}
\frak{t}_{1}:=\text{\emph{\ gcd}}\left( \alpha ,l\right) ,\ \ \frak{t}_{2}:=%
\text{\emph{gcd}}\left( \beta ,l\right) =\text{\emph{\ gcd}}\left( \alpha
+\left( r-2\right) ,l\right)   \label{MA1}
\end{equation}
and 
\begin{equation}
\frak{z}_{1}:=\frac{l}{\frak{t}_{2}},\ \ \ \ \frak{z}_{2}:=\frac{\alpha
+\left( r-2\right) }{\frak{t}_{2}}\smallskip   \label{MA2}
\end{equation}
$\bullet $ After that, if $\frak{z}_{2}\neq 1$, express\emph{\ }$\frak{z}%
_{1}/\frak{z}_{2}$\emph{\ }as regular continued fraction\smallskip 
\begin{equation}
\fbox{$
\begin{array}{ccc}
&  &  \\ 
& \dfrac{\frak{z}_{1}}{\frak{z}_{2}}=\left[ a_{1},a_{2},\ldots ,a_{\nu
-1},a_{\nu }\right]  &  \\ 
&  & 
\end{array}
$}\smallskip   \label{MA3}
\end{equation}
and define\emph{\ }$\frak{c}_{1}\in \Bbb{Z}_{<0}$\emph{, }$\frak{c}_{2}\in 
\Bbb{N}$\emph{, }by\medskip\ the formulae\emph{:} 
\[
\frak{c}_{1}=\left\{ 
\begin{array}{ll}
\dfrac{-\frak{z}_{2}}{\left[ a_{\nu },a_{\nu -1},\ldots ,a_{3},a_{2}\right] }
& \text{\emph{if} \ }a_{2}\geq 2\text{ \ \emph{and \ }}\nu \text{\emph{\ odd}%
} \\ 
\  &  \\ 
\dfrac{-\frak{z}_{2}}{\left[ a_{\nu },a_{\nu -1},\ldots ,a_{4},a_{3}+1\right]
} & \text{\emph{if} \ }a_{2}=1\text{ \ \emph{and \ }}\nu \text{\emph{\ odd}}
\\ 
&  \\ 
\left( \dfrac{1}{\left[ a_{\nu },a_{\nu -1},\ldots ,a_{3},a_{2}\right] }%
-1\right) \frak{z}_{2} & \text{\emph{if} \ }a_{2}\geq 2\text{ \ \emph{and \ }%
}\nu \text{\emph{\ even}} \\ 
&  \\ 
\left( \dfrac{1}{\left[ a_{\nu },a_{\nu -1},\ldots ,a_{4},a_{3}+1\right] }%
-1\right) \frak{z}_{2} & \text{\emph{if} \ }a_{2}=1\text{ \ \emph{and \ }}%
\nu \text{\emph{\ even}}
\end{array}
\right. 
\]
and 
\[
\frak{c}_{2}=\left\{ 
\begin{array}{ll}
\dfrac{\frak{z}_{1}}{\left[ a_{\nu },a_{\nu -1},\ldots ,a_{2},a_{1}\right] }
& \text{\emph{if} \ }a_{1}\geq 2\text{ \ \emph{and \ }}\nu \text{\emph{\ odd}%
} \\ 
\  &  \\ 
\dfrac{\frak{z}_{1}}{\left[ a_{\nu },a_{\nu -1},\ldots ,a_{3},a_{2}+1\right] 
} & \text{\emph{if} \ }a_{1}=1\text{ \ \emph{and \ }}\nu \text{\emph{\ odd}}
\\ 
&  \\ 
\left( 1-\dfrac{1}{\left[ a_{\nu },a_{\nu -1},\ldots ,a_{2},a_{1}\right] }%
\right) \frak{z}_{1} & \text{\emph{if} \ }a_{1}\geq 2\text{ \ \emph{and \ }}%
\nu \text{\emph{\ even}} \\ 
&  \\ 
\left( 1-\dfrac{1}{\left[ a_{\nu },a_{\nu -1},\ldots ,a_{3},a_{2}+1\right] }%
\right) \frak{z}_{1} & \text{\emph{if} \ }a_{1}=1\text{ \ \emph{and \ }}\nu 
\text{\emph{\ even}}
\end{array}
\medskip \right. \medskip 
\]
\emph{(}For $\frak{z}_{2}=1$, set $\frak{c}_{1}:=-1,\frak{c}_{2}:=\frak{z}%
_{1}+1$\emph{).\medskip }\newline
$\bullet $ Finally define 
\begin{equation}
\fbox{$
\begin{array}{ccc}
&  &  \\ 
& \breve{p}:=\dfrac{\frak{c}_{1}\cdot l+\frak{c}_{2}\cdot \alpha }{\frak{t}%
_{1}},\ \ \ \ \ q:=\dfrac{l}{\frak{t}_{1}\cdot \frak{t}_{2}},\ \ \ \ \ p:=%
\left[ \,\breve{p}\,\right] _{q}\  &  \\ 
&  & 
\end{array}
$}\medskip   \label{MA4}
\end{equation}
and if $p\neq 0$, write $q/p$ as regular continued fraction\medskip\ 
\begin{equation}
\fbox{$
\begin{array}{ccc}
&  &  \\ 
& \dfrac{q}{p}=\left[ \lambda _{1},\lambda _{2},\ldots ,\lambda _{\kappa
-1},\lambda _{\kappa }\right]  &  \\ 
&  & 
\end{array}
$}  \label{MA5}
\end{equation}
$\blacktriangleright $ Then $X\left( N_{G},\Delta _{G}\right) $ admits
crepant, $T_{N_{G}}$-equivariant, full desingularizations \emph{(}i.e., 
\emph{(\ref{HILBCON}) }is satisfied, as in \emph{thm. \ref{MT1}),} and at
least one of them is projective, \textbf{if and only if} one of the
following $\emph{(}$mutually exclusive\emph{)} conditions \emph{(i), (ii)}
is fulfilled\emph{:\medskip }\newline
\emph{(i) }The greatest common divisor of\emph{\ }$\alpha ,\beta $ and\emph{%
\ }$l$\emph{\ }equals\smallskip 
\begin{equation}
\fbox{$
\begin{array}{ccc}
& \emph{gcd}\left( \alpha ,\beta ,l\right) =r-2 & 
\end{array}
$}\smallskip   \label{CON1}
\end{equation}
\emph{(ii)} The greatest common divisor of\emph{\ }$\alpha ,\beta $ and\emph{%
\ }$l$ equals $1$, $\left[ \frak{t}_{1}\right] _{r-2}=\left[ \frak{t}_{2}%
\right] _{r-2}=1$, and \textbf{either} \ $p=0$ \emph{(}and consequently $q=1$%
\emph{)}\textbf{\ or } the above defined characteristic numbers satisfy the
following relations\emph{:\medskip } \emph{\ } 
\begin{equation}
\fbox{$
\begin{array}{ccc}
&  &  \\ 
& \left\{ 
\begin{array}{l}
\dfrac{\breve{p}-p}{q}\equiv 0\ \emph{mod}\left( r-2\right) , \\ 
\\ 
\lambda _{2j}\equiv 0\ \text{\emph{mod}}\left( r-2\right) ,\ \forall j,\ \
j\in \left\{ 1,2,\ldots ,\QDOVERD\lfloor \rfloor {\kappa -1}{2}\right\}
,\smallskip  \\ 
\text{\emph{whenever }}\kappa \geq 3,\text{\emph{and}} \\ 
\\ 
\lambda _{\kappa }\equiv 1\ \text{\emph{mod}}\left( r-2\right) \smallskip 
\\ 
\text{\emph{in the case in which the length }}\kappa \,\,\left( \geq
2\right) \ \text{\emph{is \textbf{even}}}.
\end{array}
\right.  &  \\ 
&  & 
\end{array}
$}\ \   \label{CON2}
\end{equation}
\end{theorem}

\noindent Though conditions (i), (ii) of theorem \ref{MT2} are fairly
restrictive, it is remarkable that they are fulfilled by several subseries
of $2$-parameter Gorenstein cyclic singularities having \textit{infinitely
many} members in \textit{each dimension}.

\begin{example}
\emph{The subseries of non-isolated singularities with defining types } 
\[
\dfrac{1}{\left( \xi +\xi ^{\prime }+1\right) \cdot \left( r-2\right) }%
\,\left( \stackunder{\left( r-2\right) \text{\emph{-times}}}{\underbrace{%
1,1,\ldots ,1,1}},\xi \cdot \left( r-2\right) ,\xi ^{\prime }\cdot \left(
r-2\right) \right) \smallskip 
\]
\emph{and }$\xi $\emph{, }$\xi ^{\prime }$\emph{\ }$\in \Bbb{N}$\emph{, gcd}$%
\left( \xi ,\xi ^{\prime }\right) =1$\emph{, }$r\geq 4$\emph{, satisfies
obviously (\ref{CON1}).}
\end{example}

\begin{example}
\emph{The subseries of isolated singularities with defining types } 
\[
\dfrac{1}{2\left( r-1\right) ^{i}+r-2}\,\left( \stackunder{\left( r-2\right) 
\text{\emph{-times}}}{\underbrace{1,1,\ldots ,1,1}},\left( r-1\right)
^{i},\left( r-1\right) ^{i}\right) \,
\]
\emph{and }$i\in \Bbb{N}$\emph{, }$r\geq 4$\emph{, satisfies (\ref{CON2})
because\ }$\breve{p}=2\left( r-1\right) ^{i}+\left( r-1\right) -2$\emph{\ \
and } 
\[
q=2\left( r-1\right) ^{i}+\left( r-1\right) -1,\ \ \ p=q-1,\ \ \ \frac{q}{p}=%
\left[ 1,q-1\right] \ .
\]
\end{example}

\begin{example}
\label{Mohris}\emph{The example of }$4$\emph{-dimensional subseries due to
Mohri \cite{Mohri}:} 
\[
\dfrac{1}{4\,\xi }\,\left( 1,1,2\,\xi -1,2\,\xi -1\right) ,\ \ \xi \in \Bbb{%
N\ },
\]
\emph{satisfies (\ref{CON2}) and contains only isolated singularities
because gcd}$\left( 4\xi ,2\xi -1\right) =1$ \emph{and} 
\[
\frac{\frak{z}_{1}}{\frak{z}_{2}}=\frac{4\,\xi }{2\xi +1}=\left[ 1,1,\xi -1,2%
\right] ,\ \ \frak{c}_{1}=-\left( \xi +1\right) ,\ \ \frak{c}_{2}=2\xi +1,
\]
\emph{i.e., } 
\[
\breve{p}=-\left( 4\xi +1\right) ,\ q=4\,\xi ,\ p=4\,\xi -1,\ \frac{q}{p}=%
\left[ 1,4\,\xi -1\right] ,\text{\emph{\ with\ }}4\,\xi -1\equiv 1\text{ }%
\left( \text{\emph{mod }}2\right) .
\]
\emph{Note that also the single suitably resolvable cyclic singularity} $%
\dfrac{1}{11}\left( 1,1,3,6\right) \ $\emph{\ found in \cite{Mohri} belongs
to the subseries of isolated cyclic quotient singularities with type} 
\[
\dfrac{1}{4r-5}\left( \stackunder{\left( r-2\right) \text{\emph{-times}}}{%
\underbrace{1,1,\ldots ,1,1}},r-1,2r-2\right) 
\]
\emph{satisfying (\ref{CON2}). Moreover, there are examples like }$%
1/28\left( 1,1,1,4,21\right) $ \emph{for which }$p=0,q=1.$
\end{example}

\begin{remark}
\emph{For those readers who would like to test rapidly if one of the above
conditions (i), (ii) of theorem \ref{MT2} is fulfilled for a concrete }$2$%
\emph{-parameter cyclic quotient (just by giving }$\alpha $\emph{, }$\beta $%
\emph{\ and }$r$\emph{\ as input), we refer to the www-page \cite{Haus} of
the second author.}
\end{remark}

\section{Proof of main theorems I, II\label{MAIN}}

\noindent In this section we prove theorems \ref{MT1} and \ref{MT2}.\bigskip

\noindent \textsf{(a)} Let $\left( X\left( N_{G},\Delta _{G}\right) ,\text{%
orb}\left( \sigma _{0}\right) \right) $ be an $r$-dimensional Gorenstein
cyclic quotient msc-singularity (with $l=\left| G\right| \geq r\geq 4$).%
\emph{\ }Suppose that its type contains at least $r-2$ equal weights.
Without loss of generality (i.e., up to analytic isomorphism, cf. cor. \ref
{ISOCYC} and prop. \ref{Gor-prop}), we may assume that $\left( X\left(
N_{G},\Delta _{G}\right) ,\text{orb}\left( \sigma _{0}\right) \right) $ is
of type\smallskip 
\begin{equation}
\fbox{$
\begin{array}{ccc}
&  &  \\ 
& \dfrac{1}{l}\ \left( \stackunder{\left( r-2\right) \text{-times}}{%
\underbrace{k,\ldots ,k}},\alpha ,\beta \right) ,\,\,\,\ \text{with\ }%
\left\{ 
\begin{array}{c}
\alpha +\beta +k\,\left( r-2\right) \equiv 0\,\left( \text{mod }l\right) ,
\\ 
\text{gcd}\left( k,\alpha ,l\right) =\text{gcd}\left( k,\beta ,l\right) =1
\end{array}
\right. \text{ } &  \\ 
&  & 
\end{array}
$}\   \label{gentype}
\end{equation}
$(\alpha ,\beta ,k\in \Bbb{N}$). Obviously,\smallskip 
\[
\mathbf{Hlb}_{N_{G}}\left( \sigma _{0}\right) \subset \left\{ e_{1},\ldots
,e_{r}\right\} \cup \left\{ \left. \frac{1}{l}\left( \stackunder{\left(
r-2\right) \text{-times}}{\underbrace{\left[ j\cdot k\right] _{l},\ldots ,%
\left[ j\cdot k\right] _{l}}},\left[ j\cdot \alpha \right] _{l},\left[
j\cdot \beta \right] _{l}\right) ^{\intercal }\ \right| \ 1\leq j\leq
l-1\right\} \ .
\]
(by prop. \ref{MINGS}). Define the hyperplane 
\[
\mathcal{H}:=\left\{ \mathbf{x=}\left( x_{1},\ldots ,x_{r}\right)
^{\intercal }\in \left( N_{G}\right) _{\Bbb{R}}\ \left| \ \sum_{i=1}^{r}\
x_{i}=1\right. \right\} \text{,}
\]
and the $3$-dimensional linear $\frak{L}$ subspace of $\left( N_{G}\right) _{%
\Bbb{R}}$%
\[
\frak{L}:=\left\{ \mathbf{x=}\left( x_{1},\ldots ,x_{r}\right) ^{\intercal
}\in \left( N_{G}\right) _{\Bbb{R}}\ \left| \ x_{1}-x_{i}=0,\ \forall i,\ \
2\leq i\leq r-2\right. \right\} \ .
\]
Next consider the $3$-dimensional s.c.p.cone 
\[
\overline{\sigma _{0}}:=\sigma _{0}\cap \frak{L}=\text{ pos}\left( \left\{ 
\frac{1}{\left( r-2\right) }\ \dsum\limits_{i=1}^{r-2}\ e_{i},\ e_{r-1},\
e_{r}\right\} \right) \subset \left( \overline{N_{G}}\right) _{\Bbb{R}%
}\subset \left( N_{G}\right) _{\Bbb{R}}
\]
supporting the triangle 
\[
\overline{\frak{s}_{G}}:=\frak{s}_{G}\cap \frak{L}=\sigma _{0}\cap \mathcal{H%
}\cap \frak{L}=\overline{\sigma _{0}}\cap \mathcal{H}=\text{ conv}\left(
\left\{ \frac{1}{\left( r-2\right) }\ \dsum\limits_{i=1}^{r-2}\ e_{i},\
e_{r-1},\ e_{r}\right\} \right) \ ,
\]
where $\frak{s}_{G}$ denotes, as usual, the corresponding junior lattice
simplex (w.r.t. $N_{G}$), and 
\[
\overline{N_{G}}:=\left\{ 
\begin{array}{c}
\text{ the sublattice of \ }N_{G}\,\text{ generated} \\ 
\text{ (as subgroup) by }N_{G}\cap \frak{L}
\end{array}
\right\} \frak{\ }.
\]
Note that if 
\[
\frak{n}_{G}:=n\left( \Bbb{R}_{\geq 0}\left( \frac{1}{\left( r-2\right) }%
\sum_{i=1}^{r-2}e_{i}\right) \right) 
\]
is the first primitive lattice point $\frak{n}_{G}$ of $\overline{N_{G}}\Bbb{%
r}\left\{ \mathbf{0}\right\} $ belonging to the ray which is defined by $%
\frac{1}{\left( r-2\right) }\sum_{i=1}^{r-2}e_{i}$ (cf. \S \ref{Prel} 
\textsf{(e)}) and\smallskip 
\[
\mu _{G}:=\text{min}\left\{ \varkappa \in \Bbb{Q}_{>0}\Bbb{\ }\left| \
\varkappa \cdot \left( \frac{1}{\left( r-2\right) }\sum_{i=1}^{r-2}e_{i}%
\right) \in \left( \overline{N_{G}}\Bbb{r}\left\{ \mathbf{0}\right\} \right)
\right. \right\} ,
\]
then 
\[
\text{Gen}\left( \overline{\sigma _{0}}\right) =\left\{ \frak{n}_{G},\
e_{r-1},\ e_{r}\right\} 
\]
and $\overline{\frak{s}_{G}}$ is a \textit{lattice triangle} w.r.t. $%
\overline{N_{G}}$ if and only if $\mu _{G}=1$(!).\medskip 

\noindent Now define the lattice polygon 
\[
\frak{Q}_{G}:=\text{ conv}\left( \frak{s}_{G}\cap \overline{N_{G}}\right)
\subset \left( \overline{N_{G}}\right) _{\Bbb{R}}. 
\]
It should be mentioned that if $\mu _{G}=1$, then $\frak{Q}_{G}=\overline{%
\frak{s}_{G}}$ and if $\mu _{G}\neq 1$, we have $\frak{Q}_{G}\subsetneqq 
\overline{\frak{s}_{G}}$.

\begin{definition}
\label{ROS}\emph{In case in which }$\mu _{G}\neq 1$\emph{\ we denote by }$%
\frak{w}$\emph{\ (resp. }$\frak{w}^{\prime }$\emph{) the unique lattice
point belonging to} 
\[
\frak{Q}_{G}\cap \text{\emph{conv}}\left( \frac{\frak{n}_{G}}{\mu _{G}}%
,e_{r-1}\right) \cap \overline{N_{G}}\text{ \thinspace \thinspace \emph{%
\thinspace \thinspace (resp. to} }\frak{Q}_{G}\cap \text{\emph{conv}}\left( 
\frac{\frak{n}_{G}}{\mu _{G}},e_{r}\right) \cap \overline{N_{G}}\text{\emph{)%
}}
\]
\emph{so that} 
\[
\text{\emph{conv}}\left( \frac{\frak{n}_{G}}{\mu _{G}},\frak{w}\right) \cap 
\overline{N_{G}}=\left\{ \frak{w}\right\} \,\,\ \text{\emph{(resp.} \emph{%
conv}}\left( \frac{\frak{n}_{G}}{\mu _{G}},\frak{w}^{\prime }\right) \cap 
\overline{N_{G}}=\left\{ \frak{w}^{\prime }\right\} \text{\emph{).}}
\]
\emph{Moreover, if conv}$\left( e_{r-1},e_{r}\right) \subsetneqq \frak{Q}_{G}
$\emph{, we fix the (clockwise ordered, uniquely determined) enumeration} 
\[
\frak{w}_{0}=\frak{w},\ \frak{w}_{1},\ \frak{w}_{2},\ \ldots \ ,\frak{w}%
_{\rho },\frak{w}_{\rho +1}=\frak{w}^{\prime }
\]
\emph{of all lattice points of }$\overline{N_{G}}$\emph{\ lying on} 
\[
\left( \left( \frak{Q}_{G}\Bbb{r}\text{\emph{conv}}\left( \left\{ \frak{w},%
\frak{w}^{\prime },e_{r-1},e_{r}\right\} \right) \right) \cap \partial \frak{%
Q}_{G}\right) \cup \left\{ \frak{w},\frak{w}^{\prime }\right\} \ .
\]
\emph{(For }$r=4$, $\mu _{G}\neq 1$\emph{,} \emph{figure\textbf{\ 2 }
illustrates these lattice points within the junior tetrahedron. Note that
the singularity is isolated if and only if }$\frak{w}=e_{r-1}$\emph{\ and }$%
\frak{w}^{\prime }=e_{r}$\emph{).}\newline
\begin{figure}[htbp]
\begin{center}
\input{fig2.pstex_t}
\center{Figure \textbf{2}}
\end{center}
\end{figure}
\end{definition}

\begin{definition}
\emph{For the given }$r\geq 4$\emph{\ we define }$\Xi _{r}$\emph{\ to be the
set} 
\[
\Xi _{r}:=\left\{ \left( \xi _{1},\xi _{2},\ldots ,\xi _{r-3}\right) \in
\left( \left\{ 1,2,\ldots ,r-2\right\} \right) ^{r-3}\ \left| \ \right.
1\leq \xi _{1}<\xi _{2}<\cdots \ \cdots <\xi _{r-3}\leq r-2\right\} .
\]
\end{definition}

\begin{definition}[Maximal triangulations constructed by ``joins'']
\label{JOI}\ \smallskip\ \ \newline
\emph{An auxiliary subclass of maximal lattice triangulations of the junior
simplex }$s_{G}$\emph{\ which can be described easily and used efficiently
for several geometric arguments is that consisting of tringulations of the
form }$\mathcal{T}$ $=\mathcal{T}\left[ \frak{T}\right] $\emph{, with} 
\[
\mathcal{T}\left[ \frak{T}\right] :=\left\{ 
\begin{array}{lll}
\frak{E}_{\frak{T}}\frak{\,}\cup \left\{ 
\begin{array}{c}
\text{\emph{conv}}\left( \left\{ \frak{w}_{i},\frak{w}_{i+1},e_{1},e_{2},%
\ldots ,e_{r-2}\right\} \right) ,\,0\leq i\leq \rho ,\smallskip  \\ 
\text{\emph{together with their faces}}
\end{array}
\right\}  & , & \text{\emph{if }\ }\mu _{G}\neq 1 \\ 
& \  &  \\ 
\frak{E}_{\frak{T}} & , & \text{\emph{if }\ }\mu _{G}=1
\end{array}
\right. 
\]
\emph{where} 
\[
\frak{T}\in \mathbf{LTR}_{\,\overline{N_{G}}}^{\emph{\scriptsize max}%
}\,\left( \frak{Q}_{G}\right) \ \left( =\mathbf{LTR}_{\,\overline{N_{G}}}^{%
\text{{\scriptsize basic}}}\,\left( \frak{Q}_{G}\right) \right) 
\]
\emph{and} 
\[
\frak{E}_{\frak{T}}:=\left\{ 
\begin{array}{c}
\left\{ \text{\emph{conv}}\left( \left\{ n_{1},n_{2},n_{3},e_{\xi
_{1}},e_{\xi _{2}},\ldots ,e_{\xi _{r-3}}\right\} \right) \ \left| 
\begin{array}{c}
\ \text{\emph{for \ all\ }}\left( \xi _{1},\xi _{2},\ldots ,\xi
_{r-3}\right) \in \Xi _{r}\text{ \ \emph{and}} \\ 
\text{\emph{all triangles conv}}\left( \left\{ n_{1},n_{2},n_{3}\right\}
\right) \ \text{\emph{of} \thinspace }\frak{T}
\end{array}
\right. \right\} \smallskip  \\ 
\text{\emph{together with their faces}}
\end{array}
\right\} \ .\medskip 
\]
$\bullet $ \emph{One easily verifies that all simplices of the above
constructed triangulations }$\mathcal{T}\left[ \frak{T}\right] $\emph{\ are
in fact representable as joins of smaller simplices; in particular we have} 
\[
\text{\emph{conv}}\left( \left\{ n_{1},n_{2},n_{3},e_{\xi _{1}},e_{\xi
_{2}},\ldots ,e_{\xi _{r-3}}\right\} \right) =\text{\emph{conv}}\left(
\left\{ n_{1},n_{2},n_{3}\right\} \right) \ast \text{\emph{conv}}\left(
\left\{ e_{\xi _{1}},e_{\xi _{2}},\ldots ,e_{\xi _{r-3}}\right\} \right) 
\]
\emph{and } 
\[
\text{\emph{conv}}\left( \left\{ \frak{w}_{i},\frak{w}_{i+1},e_{1},e_{2},%
\ldots ,e_{r-2}\right\} \right) =\text{\emph{conv}}\left( \left\{ \frak{w}%
_{i},\frak{w}_{i+1}\right\} \right) \ast \text{\emph{conv}}\left( \left\{
e_{1},e_{2},\ldots ,e_{r-2}\right\} \right) ,\,\forall i,\ 0\leq i\leq \rho ,
\]
\emph{respectively. Hence, we may alternatively describe }$\frak{E}_{\frak{T}%
}$ \emph{and }$\mathcal{T}\left[ \frak{T}\right] $ \emph{as\smallskip\ } 
\[
\frak{E}_{\frak{T}}=\left\{ \frak{T}\ast \text{\emph{conv}}\left( \left\{
e_{\xi _{1}},e_{\xi _{2}},\ldots ,e_{\xi _{r-3}}\right\} \right) \ \left| \
\right. \ \text{\emph{for \ all\ \ }}\left( \xi _{1},\xi _{2},\ldots ,\xi
_{r-3}\right) \in \Xi _{r}\right\} 
\]
\emph{and} 
\[
\mathcal{T}\left[ \frak{T}\right] =\left\{ 
\begin{array}{lll}
\frak{E}_{\frak{T}}\frak{\,}\cup \left[ \left( \dbigcup\limits_{i=0}^{\rho }%
\text{\emph{conv}}\left( \left\{ \frak{w}_{i},\frak{w}_{i+1}\right\} \right)
\right) \ast \text{\emph{conv}}\left( \left\{ e_{1},e_{2},\ldots
,e_{r-2}\right\} \right) \right]  & , & \text{\emph{if }\ }\mu _{G}\neq 1 \\ 
& \  &  \\ 
\frak{E}_{\frak{T}} & , & \text{\emph{if }\ }\mu _{G}=1
\end{array}
\right. 
\]
\end{definition}

\begin{remark}
\emph{(i) If \ }$\left( X\left( N_{G},\Delta _{G}\right) ,\emph{orb}\left(
\sigma _{0}\right) \right) $ \emph{as in thm. \ref{MT1} is a }\textit{pure}%
\emph{\ }$2$\emph{-parameter singularity, i.e., if not more than }$r-2$\emph{%
\ of the weights of its defining type are equal, then there is no }\textit{%
unique}\emph{\ full (resp. maximal partial) crepant resolution of }$X\left(
N_{G},\Delta _{G}\right) $\emph{, but at least such resolutions due to
triangulations of the form }$\mathcal{T}\left[ \frak{T}\right] $ \emph{are
induced by maximal (and therefore basic) triangulations }$\frak{T}$ \emph{of
the lattice polygon }$\frak{Q}_{G}$\emph{. Hence, the }\textit{flops }\emph{%
connecting any two of them are induced by classical ``elementary
transformations'' (cf. Oda \cite{Oda}, prop. 1.30, p. 49).\medskip }\newline
\emph{(ii) If \ a singularity }$\left( X\left( N_{G},\Delta _{G}\right) ,%
\emph{orb}\left( \sigma _{0}\right) \right) $ \emph{as in thm. \ref{MT1} has
a crepant, full resolution coming from a triangulation of the form
\thinspace }$\mathcal{T}\left[ \frak{T}\right] $\emph{, then applying
techniques similar to those of \cite{DH} (by considering the stars of the
vertices }$n$\emph{\ of $\mathcal{T}$ } \emph{and the corresponding closures 
}$V\left( \Bbb{R}_{\geq 0}n\right) $\emph{, cf. \S \ref{Prel}\textsf{(g)})
it is possible to specify the structure of the exceptional prime divisors up
to analytic isomorphism. For instance, all compactly supported exceptional
prime divisors w.r.t. }$f_{\mathcal{T}\left[ \frak{T}\right] }$ $\ $\ \emph{%
are the total spaces of fibrations having basis }$\Bbb{P}_{\Bbb{C}}^{r-3}$%
\emph{\ and typical fiber isomorphic either to }$\Bbb{P}_{\Bbb{C}}^{1}$\emph{%
\ or to a }$2$\emph{-dimensional compact toric variety (i.e. to a }$\Bbb{P}_{%
\Bbb{C}}^{2}$\emph{\ or }$\Bbb{F}_{\varkappa }=\Bbb{P}\left( \mathcal{O}_{%
\Bbb{P}_{\Bbb{C}}^{1}}\oplus \mathcal{O}_{\Bbb{P}_{\Bbb{C}}^{1}}\left(
\varkappa \right) \right) $\emph{, probably blown up at finitely many
points, cf. \cite{Oda}, thm. 1.28, p. 42).\medskip }
\end{remark}

\noindent

\begin{proposition}
\label{COHERENCY}The set $\mathbf{Coh}$-$\mathbf{LTR}_{\,\overline{N_{G}}}^{%
\emph{max}}\,\left( \frak{Q}_{G}\right) $ is non-empty. Moreover, if $\frak{T%
}\in \mathbf{Coh}$-$\mathbf{LTR}_{\,\overline{N_{G}}}^{\emph{max}}\,\left( 
\frak{Q}_{G}\right) $, then 
\[
\mathcal{T}\left[ \frak{T}\right] \in \mathbf{Coh}\text{-}\mathbf{LTR}%
_{N_{G}}^{\emph{\scriptsize max}}\,\left( \frak{s}_{G}\right) \ .
\]
\end{proposition}

\noindent \textit{Proof. }$\mathbf{Coh}$-$\mathbf{LTR}_{\,\overline{N_{G}}}^{%
\text{max}}\,\left( \frak{Q}_{G}\right) \neq \varnothing $ follows from
prop. \ref{NON-EM}. We shall henceforth fix a coherent, strictly upper
convex support funtion $\theta :\left| \frak{T}\right| \rightarrow \Bbb{R}$.
Now since $X\left( N_{G},\Delta _{G}\right) $ itself is an affine toric
variety, it is quasiprojective. Hence, besides $\theta $, there is also
another coherent, strictly upper convex support function $\phi :\frak{s}%
_{G}\rightarrow \Bbb{R}$. For any index-subset $\left\{ i_{1},\ldots
,i_{\kappa }\right\} \subset \left\{ 1,2,\ldots ,r\right\} $ let $\left.
\phi \right| _{i_{1},\ldots ,i_{\kappa }}$ denote the restriction of $\phi $
onto conv$\left( e_{i_{1}},e_{i_{2}},\ldots ,e_{i_{\kappa }}\right) $. We
define a support function 
\[
\psi :\left| \,\mathcal{T}\left[ \frak{T}\right] \,\right| \longrightarrow 
\Bbb{R} 
\]
by setting 
\[
\psi \left( \delta \cdot \mathbf{x+}\left( 1-\delta \right) \cdot \mathbf{y}%
\right) :=\delta \cdot \theta \left( \mathbf{x}\right) +\left( 1-\delta
\right) \cdot \left. \phi \right| _{\xi _{1},\ldots ,\xi _{r-3}}\left( 
\mathbf{y}\right) 
\]
for all $\mathbf{x\in }\left| \frak{T}\right| $, $\mathbf{y}\in $ conv$%
\left( \left\{ e_{\xi _{1}},\ldots ,e_{\xi _{r-3}}\right\} \right) $, $%
\delta \in \left[ 0,1\right] $ and all $\left( \xi _{1},\xi _{2},\ldots ,\xi
_{r-3}\right) \in \Xi _{r}$, and by 
\[
\psi \left( \delta \cdot \mathbf{x+}\left( 1-\delta \right) \cdot \mathbf{y}%
\right) :=\delta \cdot \theta \left( \mathbf{x}\right) +\left( 1-\delta
\right) \cdot \left. \phi \right| _{1,2,\ldots ,r-2}\left( \mathbf{y}\right) 
\]
for all $\mathbf{x\in }$ conv$\left( \left\{ \frak{w}_{i},\frak{w}%
_{i+1}\right\} \right) $, $\mathbf{y}\in $ conv$\left( \left\{ e_{1},\ldots
,e_{r-2}\right\} \right) $, $\delta \in \left[ 0,1\right] $ and for all $i$, 
$0\leq i\leq \rho $ (whenever $\mu _{G}\neq 1$). It is easy to verify that $%
\psi \in $ SUCSF$_{\Bbb{R}}\left( \mathcal{T}\left[ \frak{T}\right] \right) $%
. $_{\Box }$

\begin{lemma}
\label{SNN}Let $n_{1},n_{2},n_{3}$ be three lattice points of $\overline{%
N_{G}}\Bbb{r}\left\{ \mathbf{0}\right\} $, such that \emph{pos}$\left(
\left\{ n_{1},n_{2},n_{3}\right\} \right) $ is a $3$-dimensional basic cone
w.r.t. $\overline{N_{G}}$. Then the cone 
\[
\emph{pos}\left( \left\{ n_{1},n_{2},n_{3},e_{\xi _{1}},e_{\xi _{2}},\ldots
,e_{\xi _{r-3}}\right\} \right) \subset \left( \overline{N_{G}}\right) _{%
\Bbb{R}}
\]
is basic w.r.t. the lattice $N_{G}$ \ for every $\left( r-3\right) $-tuple $%
\left( \xi _{1},\xi _{2},\ldots ,\xi _{r-3}\right) \in \Xi _{r}$.
\end{lemma}

\noindent \textit{Proof. }For a $u\in N_{G}$ there exist $\gamma _{1},\gamma
_{2},\ldots ,\gamma _{r},\gamma _{r+1}\in \Bbb{Z}$, such that 
\[
u=\sum_{i=1}^{r}\ \gamma _{i}\ e_{i}+\gamma _{r+1}\left( \frac{1}{l}\left(
k,\ldots ,k,\alpha ,\beta \right) ^{\intercal }\right) \ .
\]
On the other hand, there exist $\kappa _{1},\kappa _{2},\ldots ,\kappa
_{r-3}\in \Bbb{Z}$, such that 
\[
u-\sum_{i=1}^{r-3}\ \kappa _{i}\ e_{\xi _{i}}\in \overline{N_{G}}\text{,\ \
\ \ }\overline{N_{G}}=\Bbb{Z\,}n_{1}\oplus \Bbb{Z\,}n_{2}\oplus \Bbb{Z\,}%
n_{3}\ ,
\] 
which means that the above cone is indeed basic. $_{\Box }$
%This means that 
%\[
%\text{mult}\left( \text{pos}\left( \left\{ n_{1},n_{2},n_{3},e_{\xi
%_{1}},\ldots ,e_{\xi _{r-3}}\right\} \right) ;N_{G}\right) =\text{mult}%
%\left( \text{pos}\left( \left\{ n_{1},n_{2},n_{3}\right\} \right) ;\overline{%
%N_{G}}\right) =1\text{, }
%\]
%i.e., that the above cone is indeed basic w.r.t. $N_{G}$. 

\begin{proposition}
\label{CENTRAL}For a Gorenstein cyclic quotient singularity $\left( X\left(
N_{G},\Delta _{G}\right) ,\emph{orb}\left( \sigma _{0}\right) \right) $ of
type \emph{(\ref{gentype})} the following conditions are equivalent\emph{%
:\smallskip }\newline
\emph{(i) } There exists a crepant, $T_{N_{G}}$-equivariant, full
desingularization 
\[
f=\emph{id}_{\ast }:X\left( N_{G},\widehat{\Delta _{G}}\right)
\longrightarrow X\left( N_{G},\Delta _{G}\right) 
\]
of the quotient space $X\left( N_{G},\Delta _{G}\right) $.\smallskip \newline
\emph{(ii)} There exists a refinement $\widehat{\Delta _{G}}$ of the fan $%
\Delta _{G}$ consisting of basic cones.\smallskip \newline
\emph{(iii) }$\mathbf{Hlb}_{N_{G}}\left( \sigma _{0}\right) =\frak{s}%
_{G}\cap N_{G}.\smallskip \smallskip $\newline
\emph{(iv) }$\overline{\frak{s}_{G}}\cap \overline{N_{G}}=\left\{ 
\begin{array}{lll}
\mathbf{Hlb}_{\,\overline{N_{G}}}\left( \overline{\sigma _{0}}\right) \Bbb{r}%
\left\{ \frak{n}_{G}\right\}  & , & \text{if \ }\mu _{G}\neq 1 \\ 
& \  &  \\ 
\mathbf{Hlb}_{\,\overline{N_{G}}}\left( \overline{\sigma _{0}}\right)  & , & 
\text{if \ }\mu _{G}=1\ .
\end{array}
\right. $
\end{proposition}

\noindent \textit{Proof. }(i)$\Leftrightarrow $(ii) follows from thm. \ref
{TRCOR}. The implication (i)$\Rightarrow $(iii) follows from thm. \ref
{KILLER}. (iii)$\Leftrightarrow $(iv) is obvious.\smallskip

\noindent (iii)$\Rightarrow $(ii): Note that $\mathbf{Hlb}_{N_{G}}\left(
\sigma _{0}\right) \Bbb{r}\left\{ e_{1},\ldots ,e_{r}\right\} \neq
\varnothing $, because otherwise $\sigma _{0}$ would be basic w.r.t. $N_{G}$
(which is impossible). More precisely, 
\[
\mathbf{Hlb}_{N_{G}}\left( \sigma _{0}\right) \Bbb{r}\left\{
e_{1},e_{2},\ldots ,e_{r-2}\right\} =\left\{ 
\begin{array}{lll}
\mathbf{Hlb}_{\,\overline{N_{G}}}\left( \overline{\sigma _{0}}\right) \Bbb{r}%
\left\{ \frak{n}_{G}\right\} & , & \text{if \ }\mu _{G}\neq 1 \\ 
& \  &  \\ 
\mathbf{Hlb}_{\overline{\,N_{G}}}\left( \overline{\sigma _{0}}\right) & , & 
\text{if \ }\mu _{G}=1
\end{array}
\right. 
\]
Now by Szeb\"{o}'s theorem \ref{Seb} there exists a proper subdivision of $%
\overline{\sigma _{0}}$ into basic subcones w.r.t. $\,\overline{N_{G}}$, 
\[
\overline{\sigma _{0}}=\dbigcup\limits_{j\in J}\,\text{pos}\left(
n_{1}^{\left( j\right) },n_{2}^{\left( j\right) },n_{3}^{\left( j\right)
}\right) , 
\]
such that Gen$\left( \overline{\sigma _{0}}\right) =\mathbf{Hlb}_{\,%
\overline{N_{G}}}\left( \overline{\sigma _{0}}\right) $. This subdivision
induces a refinement 
\[
\widehat{\Delta _{G}}=\left\{ \text{pos}\left( \left\{ n_{1}^{\left(
j\right) },n_{2}^{\left( j\right) },n_{3}^{\left( j\right) },e_{\xi
_{1}},e_{\xi _{2}},\ldots ,e_{\xi _{r-3}}\right\} \right) \ \left| \ \text{%
for \ all\emph{\ }}\left( \xi _{1},\xi _{2},\ldots ,\xi _{r-3}\right) \in
\Xi _{r}\text{ \ and\ }j\in J\right. \right\} 
\]
of $\Delta _{G}$ with Gen$\left( \widehat{\Delta _{G}}\right) \subset 
\mathbf{Hlb}_{N_{G}}\left( \sigma _{0}\right) $ ($=\frak{s}_{G}\cap N_{G}$).
In particular, if for some index $j\in J$ one of the lattice points $%
n_{1}^{\left( j\right) }$, $n_{2}^{\left( j\right) }$, $n_{3}^{\left(
j\right) }$, say $n_{3}^{\left( j\right) }$, equals $\frak{n}_{G}$, then 
\[
\text{pos}\left( \left\{ n_{1}^{\left( j\right) },n_{2}^{\left( j\right)
},n_{3}^{\left( j\right) },e_{\xi _{1}},e_{\xi _{2}},\ldots ,e_{\xi
_{r-3}}\right\} \right) =\text{pos}\left( \left\{ n_{1}^{\left( j\right)
},n_{2}^{\left( j\right) },e_{1},e_{2},\ldots ,e_{r-3},e_{r-2}\right\}
\right) \text{.} 
\]
Applying lemma \ref{SNN} we deduce that all cones of $\widehat{\Delta _{G}}$
are basic w.r.t. $N_{G}$. $_{\Box }\medskip $

\noindent To win projectivity we shall use a consequence of this proposition
based on the following lemma.

\begin{lemma}
\label{BOH}Let $N$ be a lattice of rank $3,$ $\sigma $ a s.c.p. cone $%
\subset N_{\Bbb{R}}$, and let $F_{1},F_{2},\ldots ,F_{\varkappa }$ denote
the compact facets of the lower convex hull \emph{conv}$\left( \sigma \cap
\left( N\Bbb{r}\left\{ \mathbf{0}\right\} \right) \right) $. Suppose that
the Hilbert basis of $\sigma $ w.r.t. $N$ equals 
\begin{equation}
\mathbf{Hlb}_{N}\left( \sigma \right) =\left( \bigcup_{i=1}^{\varkappa
}\,F_{i}\right) \cap N\,\,.  \label{HILBER}
\end{equation}
If $\ \left\{ \mathbf{s}_{i,1},\mathbf{s}_{i,2},\ldots ,\mathbf{s}_{i,\pi
_{i}}\right\} $ is an arbitrary triangulation of $F_{i}$ into elementary 
\emph{(}and therefore basic\emph{) }lattice triangles, then 
\[
\bigcup_{i=1}^{\varkappa }\,\bigcup_{j=1}^{\pi _{i}}\,\text{\emph{pos}}%
\left( \text{\emph{\thinspace }}\mathbf{s}_{i,j}\right) 
\]
constitutes a subdivision of $\sigma $ into basic subcones w.r.t. $N$, such
that for all $j$, $1\leq j\leq \pi _{i}$ and all $i$, $1\leq i\leq \varkappa 
$, 
\[
\text{\emph{Gen}}\left( \text{\emph{pos}}\left( \text{\emph{\thinspace }}%
\mathbf{s}_{i,j}\right) \right) \subset \mathbf{Hlb}_{N}\left( \sigma
\right) \,.
\]
\end{lemma}

\noindent \textit{Proof. }It suffices to show that all $\,$pos$\left( \text{%
\emph{\thinspace }}\mathbf{s}_{i,j}\right) $ are basic w.r.t. $N$. If there
were indices $j=j_{\bullet }\in \left\{ 1,..,\pi _{i}\right\} $ and $%
i=i_{\bullet }\in $ $\left\{ 1,..,\varkappa \right\} $, such that vert$%
\left( \mathbf{s}_{i_{\bullet },j_{\bullet }}\right) =\left\{ \mathbf{n}_{1},%
\mathbf{n}_{2},\mathbf{n}_{3}\right\} $ and pos$\left( \mathbf{s}%
_{i_{\bullet },j_{\bullet }}\right) $ a non-basic cone, then we would find a
point $\mathbf{n}_{\bullet }\in \mathbf{Hlb}_{N}\left( \text{pos}\left( 
\mathbf{s}_{i_{\bullet },j_{\bullet }}\right) \right) \Bbb{r}\left\{ \mathbf{%
n}_{1},\mathbf{n}_{2},\mathbf{n}_{3}\right\} $. Writing the affine hull of $%
F_{i_{\bullet }}$ as a hyperplane 
\[
\text{aff}\left( F_{i_{\bullet }}\right) =\left\{ \mathbf{x\in }N_{\Bbb{R}}\
\left| \ \left\langle \mathbf{m}_{\bullet }\mathbf{,x}\right\rangle =\gamma
\right. \right\} ,\text{ with\ }\gamma \in \Bbb{N},\ \ \mathbf{m}_{\bullet
}\in \text{Hom}_{\Bbb{Z}}\left( N,\Bbb{Z}\right) ,
\]
we would have $\left\langle \mathbf{m}_{\bullet }\mathbf{,n}%
_{1}\right\rangle =$ $\left\langle \mathbf{m}_{\bullet }\mathbf{,n}%
_{2}\right\rangle =\left\langle \mathbf{m}_{\bullet }\mathbf{,n}%
_{3}\right\rangle =\gamma $ on the one hand, and 
\begin{equation}
\ \left\langle \mathbf{m}_{\bullet }\mathbf{,x}\right\rangle \geq \gamma ,\
\ \forall \mathbf{x},\ \ \ \mathbf{x}\in \sigma \cap \left( N\Bbb{r}\left\{ 
\mathbf{0}\right\} \right)   \label{ANISOS}
\end{equation}
on the other. Expressing $\mathbf{n}_{\bullet }$ as linear combination of
the form (cf. (\ref{HilbP}))
\[
\mathbf{n}_{\bullet }=\delta _{1}\,\mathbf{n}_{1}+\delta _{2}\,\mathbf{n}%
_{2}+\delta _{3}\,\mathbf{n}_{3},\text{ \ \ with \thinspace \thinspace\ }%
0\leq \delta _{1},\delta _{2},\delta _{3}<1,
\]
(where at least two of $\delta _{1},\delta _{2},\delta _{3}$ are $\neq 0$), and assuming that $\left\langle \mathbf{m}_{\bullet }%
\mathbf{,n}_{\bullet }\right\rangle \geq 2\,\gamma $, we would obtain 
\[
\left\langle \mathbf{m}_{\bullet }\mathbf{,\mathbf{n}_{1}+\mathbf{n}_{2}+%
\mathbf{n}_{3}-n}_{\bullet }\right\rangle \leq \,\gamma \Longrightarrow 
\mathbf{\mathbf{n}_{1}+\mathbf{n}_{2}+\mathbf{n}_{3}-n}_{\bullet }\in \left( 
\mathbf{s}_{i_{\bullet },j_{\bullet }}\Bbb{r}\left\{ \mathbf{n}_{1},\mathbf{n%
}_{2},\mathbf{n}_{3}\right\} \right) 
\]
contradicting the fact that $\mathbf{s}_{i_{\bullet },j_{\bullet }}$ itself
is basic by construction. Thus, $\left\langle \mathbf{m}_{\bullet }\mathbf{,n%
}_{\bullet }\right\rangle <2\,\gamma $. If there were two lattice points $%
\mathbf{n}_{\bullet }^{\prime },\mathbf{n}_{\bullet }^{\prime \prime }$
belonging to $\sigma \cap \left( N\Bbb{r}\left\{ \mathbf{0}\right\} \right) $%
, such that $\mathbf{n}_{\bullet }=$ $\mathbf{n}_{\bullet }^{\prime }+%
\mathbf{n}_{\bullet }^{\prime \prime }$, then by (\ref{ANISOS}): 
\[
2\,\gamma >\left\langle \mathbf{m}_{\bullet }\mathbf{,n}_{\bullet
}\right\rangle =\left\langle \mathbf{m}_{\bullet }\mathbf{,n}_{\bullet
}^{\prime }\right\rangle +\left\langle \mathbf{m}_{\bullet }\mathbf{,n}%
_{\bullet }^{\prime \prime }\right\rangle \geq 2\,\gamma ,
\]
which is impossible. This means that $\mathbf{n}_{\bullet }$ necessarily
belongs to $\mathbf{Hlb}_{N}\left( \sigma \right) $ and by (\ref{HILBER}): 
\[
\mathbf{n}_{\bullet }\in \left( \mathbf{s}_{i_{\bullet },j_{\bullet }}\Bbb{r}%
\left\{ \mathbf{n}_{1},\mathbf{n}_{2},\mathbf{n}_{3}\right\} \right) ,
\]
which again contradicts our hypothesis. $_{\Box }$

\begin{proposition}
\label{COBAS}If the conditions of prop. \emph{\ref{CENTRAL}} for the
quotient singularity $\left( X\left( N_{G},\Delta _{G}\right) ,\emph{orb}%
\left( \sigma _{0}\right) \right) $ of type \emph{(\ref{gentype}) }are
satisfied, then all partial crepant, $T_{N_{G}}$-equivariant,
desingularizations 
\[
f=\text{\emph{id}}_{\ast }:X\left( N_{G},\widehat{\Delta _{G}}\left( 
\mathcal{T}\left[ \frak{T}\right] \right) \right) \longrightarrow X\left(
N_{G},\Delta _{G}\right) 
\]
of $X\left( N_{G},\Delta _{G}\right) $ induced by maximal triangulations of
the junior simplex $\frak{s}_{G}$ of the form $\mathcal{T}\left[ \frak{T}%
\right] $ \emph{(}as in \emph{\ref{JOI}) }are in particular \emph{``}full%
\emph{'' (}and hence all $\mathcal{T}\left[ \frak{T}\right] $'s basic w.r.t. 
$N_{G}$\emph{). }
\end{proposition}

\noindent \textit{Proof. }\ For the $3$-dimensional s.c.p. cone $\overline{%
\sigma _{0}}=$ pos$\left( \left\{ \frak{n}_{G},e_{r-1},e_{r}\right\} \right) 
$ condition (iv) of proposition \ref{CENTRAL} implies: 
\[
\mathbf{Hlb}_{\overline{N_{G}}}\left( \overline{\sigma _{0}}\right) =\left\{ 
\frak{n}_{G}\right\} \cup \left( \overline{\frak{s}_{G}}\cap \overline{N_{G}}%
\right) =\left\{ \frak{n}_{G}\right\} \cup \left( \frak{Q}_{G}\cap \overline{%
N_{G}}\right) \ .
\]
Hence, 
\[
\left\{ 
\begin{array}{c}
\text{compact facets of the lower } \\ 
\text{convex hull conv}\left( \overline{\sigma _{0}}\cap \left( \overline{%
N_{G}}\Bbb{r}\left\{ \mathbf{0}\right\} \right) \right) 
\end{array}
\right\} =\left\{ 
\begin{array}{ll}
\frak{Q}_{G}, & \text{for }\mu _{G}=1 \\ 
&  \\ 
\frak{Q}_{G}\cup \tbigcup\limits_{j=0}^{\rho }\text{conv}\left( \left\{ 
\frak{n}_{G},\frak{w}_{j},\frak{w}_{j+1}\right\} \right) , & \text{for }\mu
_{G}\neq 1.
\end{array}
\right. 
\]
Applying lemma \ref{BOH} for any $\frak{T}\in \mathbf{LTR}_{\,\overline{N_{G}%
}}^{\text{basic}}\,\left( \frak{Q}_{G}\right) $ $\ $(for $\mu _{G}=1$, resp. 
$\frak{T}\cup \tbigcup\limits_{j=0}^{\rho }$ conv$\left( \left\{ \frak{w}%
_{j},\frak{w}_{j+1},\frak{n}_{G}\right\} \right) $, for $\mu _{G}\neq 1$) we
obtain a subdivision of $\overline{\sigma _{0}}$ into basic subcones w.r.t. $%
\overline{N_{G}}$ whose set of minimal generators belongs to the Hilbert
basis of $\overline{\sigma _{0}}$. By lemma \ref{SNN} all triangulations $%
\mathcal{T}\left[ \frak{T}\right] $ have to be basic too (w.r.t. $N_{G}$). $%
_{\Box }$

\begin{corollary}
\label{CORPROJ}If the conditions of prop. \emph{\ref{CENTRAL}} are
satisfied, then there exists a crepant, $T_{N_{G}}$-equivariant, full, 
\textbf{projective} desingularization of the quotient space $X\left(
N_{G},\Delta _{G}\right) $.
\end{corollary}

\noindent \noindent \textit{Proof. }It follows from \ref{COHERENCY} and \ref
{COBAS}. $_{\Box }\medskip \smallskip $

\noindent $\bullet $ \textbf{Proof of thm. \ref{MT1}: }It follows
straightforwardly from \ref{KILLER}, \ref{CENTRAL} and\emph{\ }\ref{CORPROJ}%
. $_{\Box }\bigskip \smallskip $

\noindent \textsf{(b) }Our strategy to prove thm. \ref{MT2} is based on the
reduction of the problem by lattice transformations to a $2$-dimensional
one, and on the application of the techniques of \S \ref{BRUCH}. Hereafter
let $\left( X\left( N_{G},\Delta _{G}\right) ,\text{orb}\left( \sigma
_{0}\right) \right) $ denote a Gorenstein cyclic quotient singularity of
type\smallskip\ (\ref{typab}). Using the above introduced notation we obtain
\begin{equation}
\overline{N_{G}}=\Bbb{Z\ }\left( \dfrac{1}{l}\ \left( 1,\ldots ,1,\alpha
,\beta \right) ^{\intercal }\right) +\Bbb{Z\ }e_{r-1}+\Bbb{Z\ }e_{r}
\label{EIDIKO}
\end{equation}
and the corresponding $\frak{n}_{G},\mu _{G}$ are given as follows:

\begin{lemma}
The first primitive lattice point $\frak{n}_{G}$ of $\overline{N_{G}}\Bbb{r}%
\left\{ \mathbf{0}\right\} $ belonging to the ray which is defined by $\frac{%
1}{\left( r-2\right) }\sum_{i=1}^{r-2}e_{i}$ equals 
\[
\frak{n}_{G}=\frac{1}{\text{\emph{gcd}}\left( \alpha ,\beta ,l\right) }\
\sum_{i=1}^{r-2}e_{i},
\]
i.e., 
\[
\emph{Gen}\left( \overline{\sigma _{0}}\right) =\left\{ \frak{n}_{G},\
e_{r-1},\ e_{r}\right\} ,
\]
because 
\[
\mu _{G}:=\text{\emph{min}}\left\{ \varkappa \in \Bbb{Q}_{>0}\Bbb{\ }\left|
\ \varkappa \cdot \left( \frac{1}{\left( r-2\right) }\sum_{i=1}^{r-2}e_{i}%
\right) \in \left( \overline{N_{G}}\Bbb{r}\left\{ \mathbf{0}\right\} \right)
\right. \right\} =\frac{r-2}{\text{\emph{gcd}}\left( \alpha ,\beta ,l\right) 
}\ .\medskip 
\]
\end{lemma}

\noindent \textit{Proof. }\ By (\ref{EIDIKO}) every lattice point $n\in 
\overline{N_{G}}$ can be written as a linear combination 
\[
n=\xi _{1}\ \left( \frac{1}{l}\,\left( \stackunder{\left( r-2\right) \text{%
-times}}{\underbrace{1,1,\ldots ,1,1}},\alpha ,\beta \right) ^{\intercal
}\right) +\xi _{2}\ e_{r-1}+\xi _{3}\ e_{r},\ \ \ \xi _{1},\xi _{2},\xi
_{3}\in \Bbb{Z\ }.
\]
We define the set 
\[
\frak{A}_{G}:=\left\{ 
\begin{array}{c}
\text{all lattice vectors of \ }\overline{N_{G}}\text{ \ whose } \\ 
\text{last two coordinates are }=0
\end{array}
\right\} \ .
\]
For an $n\in \overline{N_{G}}$ to belong to $\frak{A}_{G}$ means that $\xi
_{1}\,\alpha +\xi _{2}\,l=\xi _{1}\,\beta +\xi _{3}\,l=0$, i.e., 
\[
\xi _{2}\in \left( \frac{-\alpha }{\text{gcd}\left( \alpha ,l\right) }%
\right) \,\Bbb{Z}\text{, \ \ \ }\xi _{3}\in \left( \frac{-\beta }{\text{gcd}%
\left( \alpha ,l\right) }\right) \,\Bbb{Z},
\]
and 
\[
\xi _{1}\in \left( \frac{l}{\text{gcd}\left( \alpha ,l\right) }\right) \,%
\Bbb{Z}\cap \left( \frac{l}{\text{gcd}\left( \alpha ,l\right) }\right) \,%
\Bbb{Z=}\left( \frac{l}{\text{gcd}\left( \alpha ,\beta ,l\right) }\right) \,%
\Bbb{Z\ }.
\]
Hence, $\frak{A}_{G}\subset \overline{N_{G}}$ can be alternatively expressed
as 
\[
\frak{A}_{G}=\Bbb{Z\ }\left( \frac{1}{\text{gcd}\left( \alpha ,\beta
,l\right) }\,\left( \stackunder{\left( r-2\right) \text{-times}}{\underbrace{%
1,1,\ldots ,1,1}},0,0\right) ^{\intercal }\right) \ .
\]
Since the last two coordinates of the vector $\frac{1}{\left( r-2\right) }%
\sum_{i=1}^{r-2}e_{i}$ are $=0$, we obtain 
\[
\mu _{G}=\text{min}\left\{ \varkappa \in \Bbb{Q}_{>0}\Bbb{\ }\left| \
\varkappa \cdot \left( \frac{1}{\left( r-2\right) }\sum_{i=1}^{r-2}e_{i}%
\right) \in \frak{A}_{G}\right. \right\} =\frac{r-2}{\text{gcd}\left( \alpha
,\beta ,l\right) }
\]
as we asserted. $_{\Box }$

\begin{remark}
\label{MI1}\emph{In case }$\mu _{G}=1$\emph{, i.e., if gcd}$\left( \alpha
,\beta ,l\right) =r-2$\emph{, it is clear by \ref{COBAS} and \ref{CORPROJ}
that there will be a crepant (resp. crepant and projective), }$T_{N_{G}}$%
\emph{-equivariant, full desingularization of the quotient space }$X\left(
N_{G},\Delta _{G}\right) $\emph{. This is the reason for which, from now on,
we shall focus our attention to the case }$\mu _{G}\neq 1$\emph{, and find
out which direct arithmetical conditions are equivalent to the geometric
ones of prop. \ref{CENTRAL} (and involve exclusively the two given
parameters }$\alpha $\emph{\ and }$\beta $\emph{). This will be done in five
steps and requires several lemmas.\smallskip }
\end{remark}

\noindent $\bullet $\textbf{\ First step}. Define the linear
transformation\smallskip 
\[
\Phi :\left( N_{G}\right) _{\Bbb{R}}\longrightarrow \left( N_{G}\right) _{%
\Bbb{R}},\ \ \ \Phi \left( \mathbf{x}\right) =A\cdot \mathbf{x,} 
\]
with 
\[
A:=\left( 
\begin{array}{ccccccccc}
0 & 0 & 0 & \cdots & 0 & 0 & -\beta & 0 & 1 \\ 
0 & 0 & 0 & \cdots & 0 & 0 & -\alpha & 1 & 0 \\ 
0 & 0 & 0 & \cdots & 0 & 0 & l & 0 & 0 \\ 
0 & 0 & 0 & \cdots & 0 & 1 & -1 & 0 & 0 \\ 
0 & 0 & 0 & \cdots & 1 & 0 & -1 & 0 & 0 \\ 
\vdots & \vdots & \vdots & \ddots & \vdots & \vdots & \vdots & \vdots & 
\vdots \\ 
0 & 0 & 1 & \cdots & 0 & 0 & -1 & 0 & 0 \\ 
0 & 1 & 0 & \cdots & 0 & 0 & -1 & 0 & 0 \\ 
1 & 0 & 0 & \cdots & 0 & 0 & \stackunder{\left( r-2\right) \text{-pos.}}{%
\underbrace{-1}} & 0 & 0
\end{array}
\right) \in \text{ GL}\left( r,\Bbb{Q}\right) 
\]
Since 
\[
N_{G}=\Bbb{Z\ }e_{1}+\Bbb{Z\ }e_{2}+\cdots ++\Bbb{Z\ }e_{r-3}+\Bbb{Z\ }%
\left( \frac{1}{l}\left( 1,\ldots ,1,\alpha ,\beta \right) ^{\intercal
}\right) +\Bbb{Z\ }e_{r-1}+\Bbb{Z\ }e_{r} 
\]
and $\Phi \left( e_{i}\right) =e_{r-i+1}$, $\forall i$, $i\in \left\{
1,\ldots ,r-3,r-1,r\right\} $, $\Phi \left( \frac{1}{l}\left( 1,\ldots
,1,\alpha ,\beta \right) ^{\intercal }\right) =e_{3}$, we have\smallskip 
\[
\Lambda _{G}:=\Phi \left( N_{G}\right) =\Bbb{Z\,}e_{1}+\Bbb{Z\,}e_{2}+\cdots
+\Bbb{Z\,}e_{r},\ \ \ \overline{\Lambda _{G}}:=\Phi \left( \overline{N_{G}}%
\right) =\Bbb{Z\,}e_{1}+\Bbb{Z\,}e_{2}+\Bbb{Z\,}e_{3}. 
\]
We can therefore work within $\left( \ \overline{\Lambda _{G}}\right) _{\Bbb{%
R}}\cong \Bbb{R}^{3}$ with the unit vectors 
\[
\text{e}_{1}=\left( 1,0,0\right) ^{\intercal },\text{e}_{2}=\left(
0,1,0\right) ^{\intercal },\text{e}_{3}=\left( 0,0,1\right) ^{\intercal } 
\]
and write\smallskip 
\[
\Phi \left( \overline{\sigma _{0}}\right) =\text{ pos}\left( \left\{ \text{e}%
_{2},\text{e}_{1},\eta _{G}\right\} \right) =\text{ pos}\left( \left\{ \text{%
e}_{2},\text{e}_{1},\frac{\eta _{G}}{\mu _{G}}\right\} \right) ,\ \ \Phi
\left( \overline{\frak{s}_{G}}\right) =\text{ conv}\left( \left( \left\{ 
\text{e}_{2},\text{e}_{1},\frac{\eta _{G}}{\mu _{G}}\right\} \right) \right) 
\]
where 
\[
\eta _{G}:=\Phi \left( \frak{n}_{G}\right) =\frac{1}{\text{gcd}\left( \alpha
,\beta ,l\right) }\ \left( 
\begin{array}{c}
-\beta \\ 
-\alpha \\ 
l
\end{array}
\right) \ \ . 
\]
Furthermore, setting 
\[
\frak{u}_{0}:=\Phi \left( \frak{w}\right) ,\ \frak{u}_{i}:=\Phi \left( \frak{%
w}_{i}\right) ,\ \forall i,\ \ 1\leq i\leq \rho ,\ \ \text{and}\ \ \frak{u}%
_{\rho +1}:=\Phi \left( \frak{w}^{\prime }\right) , 
\]
we obtain 
\[
\Phi \left( \frak{Q}_{G}\right) =\text{conv}\left( \Phi \left( \overline{%
\frak{s}_{G}}\right) \cap \overline{N_{G}}\right) =\text{conv}\left( \left\{ 
\text{e}_{2},\frak{u}_{0},\frak{u}_{1},\ldots ,\frak{u}_{\rho },\frak{u}%
_{\rho +1},\text{e}_{1}\right\} \right) \subset \left( \ \overline{\Lambda
_{G}}\right) _{\Bbb{R}} 
\]
(see figure $\mathbf{3}$).

\begin{lemma}
\label{hilmil}If $\mu _{G}\neq 1$, then each of the conditions of prop. 
\emph{\ref{CENTRAL}} is equivalent to 
\[
\mathbf{Hlb}_{\,\overline{\Lambda _{G}}}\left( \Phi \left( \overline{\sigma
_{0}}\right) \right) \Bbb{r}\left\{ \eta _{G}\right\} =\Phi \left( \overline{%
\frak{s}_{G}}\right) \cap \overline{\Lambda _{G}}\ .
\]
\end{lemma}

\noindent\textit{Proof.} Since $\Phi \left( \mathbf{Hlb}_{\,\overline{N_{G}}%
}\left( \overline{\sigma _{0}}\right) \right) =\mathbf{Hlb}_{\,\overline{%
\Lambda _{G}}}\left( \Phi \left( \overline{\sigma _{0}}\right) \right) $,
this follows from \ref{CENTRAL} (iv). $_{\Box }$

\begin{lemma}
\label{UNIM1}If $\mu _{G}\neq 1$, then each of the conditions of prop. \emph{%
\ref{CENTRAL}} is equivalent to the following: 
\[
\text{\emph{pos}}\left( \left\{ \frak{u}_{i-1},\frak{u}_{i},\eta
_{G}\right\} \right) \text{ \ \emph{is basic w.r.t. } }\overline{\Lambda _{G}%
}\text{\emph{,} \ \emph{for all }}\emph{i}\text{\emph{,} \ }1\leq i\leq \rho
+1.
\]
\end{lemma}

\noindent \textit{Proof. }``$\Rightarrow $'':  This implication follows from
lemmas \ref{BOH} and\ \ref{hilmil} applied to the $3$-dimensional cone $\Phi
\left( \overline{\sigma _{0}}\right) $.\smallskip \newline
``$\Leftarrow $'': Consider a maximal triangulation $\left\{ \mathbf{s}_{j}\
\left| \ j\in J\right. \right\} $ of the lattice polygon $\Phi \left( \frak{Q%
}_{G}\right) $. By lemma \ref{EL-BA}(i) this triangulation has to be basic
w.r.t. the sublattice of $\overline{\Lambda _{G}}$ generated (as subgroup)
by aff$\left( \frak{Q}_{G}\right) $. Consequently 
\[
\left\{ \text{pos}\left( \mathbf{s}_{j}\right) \ \left| \ j\in J\right.
\right\} \cup \left\{ \text{pos}\left( \left\{ \frak{u}_{i-1},\frak{u}%
_{i},\eta _{G}\right\} \right) \ \left| \ 1\leq i\leq \rho +1\right.
\right\} 
\]
constitutes a subdivision of the entire $\Phi \left( \overline{\sigma _{0}}%
\right) $ into basic cones w.r.t. $\overline{\Lambda _{G}}$ whose set of
minimal generators is exactly $\left( \Phi \left( \overline{\frak{s}_{G}}%
\right) \cap \overline{\Lambda _{G}}\right) \cup \left\{ \eta _{G}\right\} $%
. It remains to use the inverse implication in the statement of lemma \ref
{hilmil}. $_{\Box }\bigskip $

\noindent $\bullet $\textbf{\ Second step}. We define the unimodular
transformation\smallskip 
\[
\Psi :\left( \overline{\Lambda _{G}}\right) _{\Bbb{R}}\longrightarrow \left( 
\overline{\Lambda _{G}}\right) _{\Bbb{R}},\ \ \ \Psi \left( \left( 
\begin{array}{c}
x_{1} \\ 
x_{2} \\ 
x_{3}
\end{array}
\right) \right) =\left( 
\begin{array}{ccc}
1 & 0 & 1 \\ 
1 & 1 & 1 \\ 
0 & 0 & 1
\end{array}
\right) \ \left( 
\begin{array}{c}
x_{1} \\ 
x_{2} \\ 
x_{3}
\end{array}
\right) \ .
\]
Obviously, 
\[
\Psi \left( \text{e}_{1}\right) =\left( 
\begin{array}{c}
1 \\ 
1 \\ 
0
\end{array}
\right) ,\ \Psi \left( \text{e}_{2}\right) =\left( 
\begin{array}{c}
0 \\ 
1 \\ 
0
\end{array}
\right) ,\ \Psi \left( \frac{\eta _{G}}{\mu _{G}}\right) =\frac{1}{r-2}\
\left( 
\begin{array}{c}
-\beta +l \\ 
r-2 \\ 
l
\end{array}
\right) =\frac{1}{r-2}\ \left( 
\begin{array}{c}
\alpha +\left( r-2\right)  \\ 
r-2 \\ 
l
\end{array}
\right) \ .
\]
Using the embedding 
\[
\Bbb{R}^{2}\ni \left( 
\begin{array}{c}
x_{1} \\ 
x_{2}
\end{array}
\right) \stackunder{\iota }{\hookrightarrow }\left( 
\begin{array}{c}
x_{1} \\ 
1 \\ 
x_{2}
\end{array}
\right) \in \left( \overline{\Lambda _{G}}\right) _{\Bbb{R}}\cong \Bbb{R}^{3}
\]
we may work with the lattice $\widetilde{\Lambda _{G}}$ of rank 2 defined
by\smallskip 
\[
\widetilde{\Lambda _{G}}:=\iota ^{-1}\left( \left\{ \left. \left( 
\begin{array}{c}
x_{1} \\ 
x_{2} \\ 
x_{3}
\end{array}
\right) \in \left( \overline{\Lambda _{G}}\right) _{\Bbb{R}}\cong \Bbb{R}%
^{3}\ \right| \ x_{2}=1\right\} \cap \overline{\Lambda _{G}}\right) 
\]
and having $\left\{ \frak{e}_{1},\frak{e}_{2}\right\} $ as a $\Bbb{Z}$-basis
(where now $\frak{e}_{1},\frak{e}_{2}$ denote here the unit vectors $\left(
1,0\right) ^{\intercal }$and $\left( 0,1\right) ^{\intercal }$ of $\left( 
\widetilde{\Lambda _{G}}\right) _{\Bbb{R}}\cong \Bbb{R}^{2}$ (!)
respectively). In particular, 
\[
\iota ^{-1}\left( \Psi \left( \Phi \left( \overline{\frak{s}_{G}}\right)
\right) \right) =\text{ conv}\left( \left( \left\{ \mathbf{0},\frak{e}_{1},%
\frak{v}_{G}\right\} \right) \right) 
\]
where $\iota ^{-1}\left( \Psi \left( \text{e}_{1}\right) \right) =\frak{e}%
_{1}$, $\ \iota ^{-1}\left( \Psi \left( \text{e}_{2}\right) \right) =\mathbf{%
0}$, and 
\[
\frak{v}_{G}:=\iota ^{-1}\left( \Psi \left( \frac{\eta _{G}}{\mu _{G}}%
\right) \right) =\frac{1}{r-2}\ \left( 
\begin{array}{c}
a+\left( r-2\right)  \\ 
l
\end{array}
\right) \ .
\]
Furthermore, setting 
\[
\widetilde{\frak{u}}_{i}:=\iota ^{-1}\left( \Psi \left( \frak{u}_{i}\right)
\right) ,\ \forall i,\ 0\leq i\leq \rho +1,
\]
we transform $\frak{Q}_{G}$ onto the polygon 
\[
\widetilde{\frak{Q}}_{G}:=\iota ^{-1}\left( \Psi \left( \Phi \left( \frak{Q}%
_{G}\right) \right) \right) =\text{conv}\left( \left\{ \mathbf{0,}\widetilde{%
\frak{u}}_{0},\widetilde{\frak{u}}_{1},\widetilde{\frak{u}}_{2},\ldots ,%
\widetilde{\frak{u}}_{\rho },\widetilde{\frak{u}}_{\rho +1},\frak{e}%
_{1}\right\} \right) \subset \left( \ \widetilde{\Lambda _{G}}\right) _{\Bbb{%
R}}
\]
(see figure $\mathbf{3}$).

\begin{lemma}
\label{UNIM2}If $\mu _{G}\neq 1$, then each of the conditions of prop. \emph{%
\ref{CENTRAL}} is equivalent to each of the following\emph{:\smallskip }%
\newline
\emph{(i) }The volumes of the triangles \emph{conv}$\left( \left\{ 
\widetilde{\frak{u}}_{i-1},\widetilde{\frak{u}}_{i},\frak{v}_{G}\right\}
\right) $ are equal to\smallskip\ 
\[
\text{\emph{Vol}}\left( \text{\emph{conv}}\left( \left\{ \widetilde{\frak{u}}%
_{i-1},\widetilde{\frak{u}}_{i},\frak{v}_{G}\right\} \right) \right) =\frac{1%
}{2\,\mu _{G}}=\frac{\text{\emph{gcd}}\left( \alpha ,\beta ,l\right) }{%
2\,\left( r-2\right) }\text{\emph{, \ \ }for\ all}\emph{\ }i\emph{,}\ 1\leq
i\leq \rho +1.
\]
\emph{(ii) }For\ all$\emph{\ }i\emph{,}\ 1\leq i\leq \rho +1$\emph{, the
triangle conv}$\left( \left\{ \widetilde{\frak{u}}_{i-1},\widetilde{\frak{u}}%
_{i},\frak{v}_{G}\right\} \right) $ is\emph{\ }basic\emph{\ }w.r.t. the
extended lattice 
\[
\widetilde{\Lambda _{G}^{\emph{ext}}}:=\widetilde{\Lambda _{G}}+\Bbb{Z\,}%
\frak{v}_{G}\ .
\]
\end{lemma}

\noindent \textit{Proof. }That (i) is equivalent to the conditions of \ref
{CENTRAL} follows from the equalities 
\[
\left| \text{det}\left( \frak{u}_{i-1},\frak{u}_{i},\eta _{G}\right) \right|
=\mu _{G}\ \left| \text{det}\left( \frak{u}_{i-1},\frak{u}_{i},\frac{\eta
_{G}}{\mu _{G}}\right) \right| =\mu _{G}\ \left( 3!\right) \ \text{Vol}%
\left( \text{conv}\left( \left\{ \mathbf{0},\frak{u}_{i-1},\frak{u}_{i},%
\frac{\eta _{G}}{\mu _{G}}\right\} \right) \right) =
\]
\[
=\mu _{G}\ \left( 3!\right) \ \text{Vol}\left( \text{conv}\left( \left\{ 
\mathbf{0},\Psi \left( \frak{u}_{i-1}\right) ,\Psi \left( \frak{u}%
_{i}\right) ,\Psi \left( \frac{\eta _{G}}{\mu _{G}}\right) \right\} \right)
\right) =\mu _{G}\ \left( 3!\right) \ \frac{1}{3}\ \text{Vol}\left( \text{%
conv}\left( \left\{ \widetilde{\frak{u}}_{i-1},\widetilde{\frak{u}}_{i},%
\frak{v}_{G}\right\} \right) \right) 
\]
and lemma \ref{UNIM1}. Now since 
\[
\text{det}\left( \widetilde{\Lambda _{G}^{\text{{\small ext}}}}\right) =%
\frac{1}{\#\left( \left\{ \left( \left[ j\cdot \alpha \right] _{\left(
r-2\right) },\left[ j\cdot l\right] _{\left( r-2\right) }\right) \ \left| \
1\leq j\leq r-2\right. \right\} \right) }=\frac{1}{\mu _{G}}
\]
the equivalence (i)$\Leftrightarrow $(ii) becomes obvious. $_{\Box }\bigskip 
$

\noindent $\bullet $\textbf{\ Third step}. We define the affine integral
transformation (w.r.t. $\widetilde{\Lambda _{G}^{\text{{\small ext}}}}$) : 
\[
\Upsilon :\left( \widetilde{\Lambda _{G}^{\text{{\small ext}}}}\right) _{%
\Bbb{R}}\longrightarrow \left( \ \widetilde{\Lambda _{G}^{\text{{\small ext}}%
}}\right) _{\Bbb{R}},\ \ \ \Upsilon \left( \left( 
\begin{array}{c}
x_{1} \\ 
x_{2}
\end{array}
\right) \right) =\left( 
\begin{array}{cc}
-1 & 0 \\ 
0 & -1
\end{array}
\right) \ \left( 
\begin{array}{c}
x_{1} \\ 
x_{2}
\end{array}
\right) +\frak{v}_{G}
\]
being the composite of a unimodular transformation (w.r.t. $\ \widetilde{%
\Lambda _{G}}$) and a lattice translation (w.r.t. $\widetilde{\Lambda _{G}^{%
\text{{\small ext}}}}$). Obviously, 
\[
\Upsilon \left( \text{conv}\left( \left\{ \widetilde{\frak{u}}_{i-1},%
\widetilde{\frak{u}}_{i},\frak{v}_{G}\right\} \right) \right) =\text{conv}%
\left( \left\{ \mathbf{0},\frak{v}_{G}-\widetilde{\frak{u}}_{i-1},\frak{v}%
_{G}-\widetilde{\frak{u}}_{i}\right\} \right) \ .
\]
Next let us consider the $2$-dimensional s.c.p. cone 
\[
\mathbf{\tau }_{G}:=\text{ pos}\left( \left\{ \frak{v}_{G}-\widetilde{\frak{u%
}}_{0},\frak{v}_{G}-\widetilde{\frak{u}}_{\rho +1}\right\} \right) \subset
\left( \Upsilon \left( \ \widetilde{\Lambda _{G}^{\text{{\small ext}}}}%
\right) \right) _{\Bbb{R}}=\left( \widetilde{\Lambda _{G}^{\text{{\small ext}%
}}}\right) _{\Bbb{R}}\cong \Bbb{R}^{2}
\]
and, as in \S \ref{BRUCH}, define 
\[
\Theta _{\mathbf{\tau }_{G}}:=\text{ conv}\left( \mathbf{\tau }_{G}\cap
\left( \left( \widetilde{\Lambda _{G}^{\text{{\small ext}}}}\right) \Bbb{r}%
\left\{ \mathbf{0}\right\} \right) \right) \subset \left( \widetilde{\Lambda
_{G}^{\text{{\small ext}}}}\right) _{\Bbb{R}}\cong \Bbb{R}^{2}
\]
and denote by $\partial \Theta _{\mathbf{\tau }_{G}}^{\mathbf{cp}}$ the part
of the boundary $\partial \Theta _{\mathbf{\tau }_{G}}$ of $\Theta _{\mathbf{%
\tau }_{G}}$ containing only its compact edges (see fig. $\mathbf{3}$).%
%\newline
\begin{figure}[htbp]
\begin{center}\vspace{-1cm}
\input{fig3.pstex_t}
\center{Figure \textbf{3}}
\end{center}
\end{figure}
Furthermore, let us denote by 
\[
\partial \Theta _{\mathbf{\tau }_{G}}^{\mathbf{cp}}\cap \widetilde{\Lambda
_{G}^{\text{{\small ext}}}}=\left\{ \frak{l}_{0}=\frak{v}_{G}-\widetilde{%
\frak{u}}_{0},\frak{l}_{1},\ldots ,\frak{l}_{\rho ^{\prime }},\frak{l}_{\rho
^{\prime }+1}=\frak{v}_{G}-\widetilde{\frak{u}}_{\rho +1}\right\} 
\]
the (clockwise ordered, uniquely determined) enumeration of the lattice
points of $\partial \Theta _{\mathbf{\tau }_{G}}^{\mathbf{cp}}$ belonging to
the extended lattice $\widetilde{\Lambda _{G}^{\text{{\small ext}}}}$. Since
we work with two \textit{different} lattices of rank $2$ it might happen
that $\rho \neq \rho ^{\prime }$ or even that 
\[
\left( \partial \Theta _{\mathbf{\tau }_{G}}^{\mathbf{cp}}\cap \widetilde{%
\Lambda _{G}^{\text{{\small ext}}}}\right) \cap \left\{ \frak{v}_{G}-%
\widetilde{\frak{u}}_{1},\frak{v}_{G}-\widetilde{\frak{u}}_{2},\ldots ,\frak{%
v}_{G}-\widetilde{\frak{u}}_{\rho }\right\} =\varnothing \ .
\]
Nevertheless, in our particular situation we have:

\begin{lemma}
\label{UNIM3}If $\mu _{G}\neq 1$, then the conditions of lemma \emph{\ref
{UNIM2}} are equivalent to each of the following\emph{:\medskip }\newline
\emph{(i) \ \ conv}$\left( \left\{ \mathbf{0},\frak{v}_{G}-\widetilde{\frak{u%
}}_{i-1},\frak{v}_{G}-\widetilde{\frak{u}}_{i}\right\} \right) $ is\emph{\ }%
basic\emph{\ }w.r.t. $\widetilde{\Lambda _{G}^{\text{\emph{ext}}}}$\emph{, }
for all\emph{\ }$i$\emph{,} \ $1\leq i\leq \rho +1.\medskip $\newline
\emph{(ii) \ }$\rho =\rho ^{\prime }$ and $\frak{l}_{i}=\frak{v}_{G}-%
\widetilde{\frak{u}}_{i}$\emph{, }for all\emph{\ }$i$\emph{,} \ $0\leq i\leq
\rho +1.\medskip $\newline
\emph{(iii) }$\rho =\rho ^{\prime }$ and $\frak{v}_{G}-\frak{l}_{i}\in \ 
\widetilde{\Lambda _{G}}=\Bbb{Z\,\frak{e}}_{1}+\Bbb{Z\,}\frak{e}_{2}$\emph{, 
} for all\emph{\ }$i$\emph{,} \ $0\leq i\leq \rho +1.$\newline
\end{lemma}

\noindent \textit{Proof. }Since $\Upsilon $ is an affine integral
transformation, it is clear that (i) is equivalent to condition (ii) of
lemma \ref{UNIM2}.\smallskip\ \newline
(i)$\Rightarrow $(ii)$\Rightarrow $(iii): The above mentioned possibilies
are in each case to be excluded because they would for at least one index $i$
contradict the fact that both conv$\left( \left\{ \mathbf{0},\frak{v}_{G}-%
\widetilde{\frak{u}}_{i-1},\frak{v}_{G}-\widetilde{\frak{u}}_{i}\right\}
\right) $ and conv$\left( \left\{ \mathbf{0},\frak{l}_{i-1},\frak{l}%
_{i}\right\} \right) $ are basic w.r.t. $\widetilde{\Lambda _{G}^{\text{%
{\small ext}}}}$. Thus, $\rho =\rho ^{\prime }$ and $\widetilde{\frak{u}}%
_{i}=\frak{v}_{G}-\frak{l}_{i}\in \widetilde{\Lambda _{G}}$\emph{, }for all%
\emph{\ }$i$\emph{,} \ $0\leq i\leq \rho +1$.\smallskip \newline
(iii)$\Rightarrow $(ii)$\Rightarrow $(i): Since conv$\left( \left\{ \mathbf{0%
},\frak{l}_{i-1},\frak{l}_{i}\right\} \right) $ are basic w.r.t. $\widetilde{%
\Lambda _{G}^{\text{{\small ext}}}}$, conv$\left( \left\{ \frak{v}_{G},\frak{%
v}_{G}-\frak{l}_{i-1},\frak{v}_{G}-\frak{l}_{i}\right\} \right) $ will be
basic too, but this time with respect to the ``smaller'' lattice $\widetilde{%
\Lambda _{G}}$. Since $\left\{ \frak{v}_{G}-\frak{l}_{0},\frak{v}_{G}-\frak{l%
}_{1},\ldots ,\frak{v}_{G}-\frak{l}_{\rho +1}\right\} $ determines again a
``lower convex hull'' by means of points belonging to $\widetilde{\Lambda
_{G}}$, using similar arguments, one shows that necessarily $\frak{l}_{i}=%
\frak{v}_{G}-\widetilde{\frak{u}}_{i}$ and conv$\left( \left\{ \mathbf{0},%
\frak{v}_{G}-\widetilde{\frak{u}}_{i-1},\frak{v}_{G}-\widetilde{\frak{u}}%
_{i}\right\} \right) $ are basic w.r.t. $\widetilde{\Lambda _{G}^{\text{%
{\small ext}}}}$. $_{\Box }\bigskip $\newline
$\bullet $\textbf{\ Fourth step}. Maintaining the assumption $\mu _{G}\neq 1$%
, consider a $\Bbb{Z}$-basis $\left\{ \frak{v}_{G}-\widetilde{\frak{u}}_{0},%
\frak{y}\right\} $ of the lattice $\widetilde{\Lambda _{G}^{\text{{\small ext%
}}}}$ (cf. lemma \ref{HNF}), such that 
\[
\frak{v}_{G}-\widetilde{\frak{u}}_{\rho +1}=p\cdot \left( \frak{v}_{G}-%
\widetilde{\frak{u}}_{0}\right) +q\cdot \frak{y}
\]
for two positive integers $p,\ q$ with $0\leq p<q$, gcd$\left( p,q\right) =1$%
, i.e., so that $\mathbf{\tau }_{G}$ becomes a $\left( p,q\right) $-cone
w.r.t. $\left\{ \frak{v}_{G}-\widetilde{\frak{u}}_{0},\frak{y}\right\} $ (in
the sense of \ref{PQC}).

\begin{lemma}
\label{TYPOS1}The multiplicity of the cone $\mathbf{\tau }_{G}$ w.r.t. $%
\widetilde{\Lambda _{G}^{\text{\emph{ext}}}}$ equals 
\begin{equation}
q=\text{\emph{mult}}\left( \mathbf{\tau }_{G};\widetilde{\Lambda _{G}^{\text{%
\emph{ext}}}}\right) =\frac{\left[ \frak{t}_{1}\right] _{\left( r-2\right)
}\cdot \left[ \frak{t}_{2}\right] _{\left( r-2\right) }}{\frak{t}_{1}\cdot 
\frak{t}_{2}}\cdot \frac{l}{\text{\emph{gcd}}\left( \alpha ,\beta ,l\right) }
\label{QMULT}
\end{equation}
where \emph{(}as in \emph{(\ref{MA1}))}, we use the abbreviations\emph{:} 
\[
\frak{t}_{1}:=\text{\emph{\ gcd}}\left( \alpha ,l\right) ,\ \ \frak{t}_{2}:=%
\text{\emph{gcd}}\left( \beta ,l\right) =\text{\emph{\ gcd}}\left( \alpha
+\left( r-2\right) ,l\right). 
\]
\end{lemma}

\noindent \textit{Proof. }Since $\frac{r-2}{\frak{t}_{2}}\cdot $ $\frak{v}%
_{G}$ (resp. $\frac{r-2}{\frak{t}_{1}}\cdot $ $\left( \frak{v}_{G}-\frak{e}%
_{1}\right) $) is the first primitive lattice point of $\left( \Bbb{R}_{\geq
0}\,\frak{v}_{G}\right) \cap \widetilde{\Lambda _{G}}$ (resp. of $\left( 
\Bbb{R}_{\geq 0}\,\left( \frak{v}_{G}-\frak{e}_{1}\right) \right) \cap 
\widetilde{\Lambda _{G}}$), we obtain 
\[
\widetilde{\frak{u}}_{0}=\QOVERD\lfloor \rfloor {\frak{t}_{2}}{r-2}\cdot 
\frac{r-2}{\frak{t}_{2}}\cdot \frak{v}_{G}\Rightarrow \frak{v}_{G}-%
\widetilde{\frak{u}}_{0}=\frac{\left[ \frak{t}_{2}\right] _{\left(
r-2\right) }}{\frak{t}_{2}}\cdot \frak{v}_{G} 
\]
and 
\[
\widetilde{\frak{u}}_{\rho +1}=\QOVERD\lfloor \rfloor {\frak{t}%
_{1}}{r-2}\cdot \frac{r-2}{\frak{t}_{1}}\cdot \left( \frak{v}_{G}-\frak{e}%
_{1}\right) \Rightarrow \frak{v}_{G}-\widetilde{\frak{u}}_{\rho +1}=\frac{%
\left[ \frak{t}_{1}\right] _{\left( r-2\right) }}{\frak{t}_{1}}\cdot \left( 
\frak{v}_{G}-\frak{e}_{1}\right) \ . 
\]
Therefore, 
\[
q=\frac{\left| \text{det}\left( \frak{v}_{G}-\widetilde{\frak{u}}_{0},\frak{v%
}_{G}-\widetilde{\frak{u}}_{\rho +1}\right) \right| }{\text{det}\left( 
\widetilde{\Lambda _{G}^{\text{{\small ext}}}}\right) }=\mu _{G}\cdot \frac{%
\left[ \frak{t}_{1}\right] _{\left( r-2\right) }}{\frak{t}_{1}}\cdot \frac{%
\left[ \frak{t}_{2}\right] _{\left( r-2\right) }}{\frak{t}_{2}}\cdot \left| 
\text{det}\left( \frak{v}_{G},\frak{v}_{G}-\frak{e}_{1}\right) \right| 
\]
and (\ref{QMULT}) follows from the equality $\left| \text{det}\left( \frak{v}%
_{G},\frak{v}_{G}-\frak{e}_{1}\right) \right| =\left| \text{det}\left( \frak{%
e}_{1},\frak{v}_{G}\right) \right| =l\,/\,\left( r-2\right) .\ _{\Box }$

\begin{lemma}
\label{TYPOS2}Suppose \emph{gcd}$\left( \alpha ,\beta ,l\right) =1$ and $%
\left[ \frak{t}_{1}\right] _{\left( r-2\right) }=\left[ \frak{t}_{2}\right]
_{\left( r-2\right) }=1$. \ If we define 
\[
\frak{z}_{1}:=\frac{l}{\frak{t}_{2}},\ \ \ \ \frak{z}_{2}:=\frac{\alpha
+\left( r-2\right) }{\frak{t}_{2}}\smallskip 
\]
\emph{(}as in \emph{(\ref{MA2})) }and if we consider two integers\emph{\ }$%
\frak{c}_{1}\in \Bbb{Z}_{<0}$\emph{, }$\frak{c}_{2}\in \Bbb{N}$, such that 
\begin{equation}
\frak{c}_{1}\cdot \frak{z}_{1}+\frak{c}_{2}\cdot \frak{z}_{2}=1\,
\label{DIOF}
\end{equation}
then we get 
\begin{equation}
q=\text{ }\frac{l}{\frak{t}_{1}\cdot \frak{t}_{2}}\,,\ \ \ \ \ \ \ p=\left[
\,\breve{p}\,\right] _{q}\,\text{ }  \label{PQTYPOI}
\end{equation}
and $\frak{y}$ can be taken to be 
\begin{equation}
\frak{y=}\left( \frac{\,\breve{p}-p\,}{q}\right) \,\left( \frak{v}_{G}-%
\widetilde{\frak{u}}_{0}\right) +\left( -\frak{c}_{1}\cdot \frak{e}_{1}+\,%
\frak{c}_{2}\cdot \frak{e}_{2}\right)   \label{NIKA}
\end{equation}
where 
\[
\breve{p}:=\dfrac{\frak{c}_{1}\cdot l+\frak{c}_{2}\cdot \alpha }{\frak{t}_{1}%
}
\]
\emph{(}as in \emph{(\ref{MA4})).}
\end{lemma}

\noindent \textit{Proof. }Under the above assumption we have obviously 
\[
\text{det}\left( \widetilde{\Lambda _{G}^{\text{{\small ext}}}}\right) =%
\frac{1}{r-2},\ \ \ \ \ q=\text{ }\frac{l}{\frak{t}_{1}\cdot \frak{t}_{2}}, 
\]
and 
\[
\left| \text{det}\left( \frak{v}_{G}-\widetilde{\frak{u}}_{0},\frak{c}%
_{1}\cdot \frak{e}_{1}+\,\frak{c}_{2}\cdot \frak{e}_{2}\right) \right| =%
\frac{1}{\frak{t}_{2}}\ \left| \text{det}\left( \frak{v}_{G}-\widetilde{%
\frak{u}}_{0},\frak{c}_{1}\cdot \frak{e}_{1}+\,\frak{c}_{2}\cdot \frak{e}%
_{2}\right) \right| =\frac{1}{r-2}\ . 
\]
Therefore, $\left\{ \frak{v}_{G}-\widetilde{\frak{u}}_{0},\frak{c}_{1}\cdot 
\frak{e}_{1}+\,\frak{c}_{2}\cdot \frak{e}_{2}\right\} $ is a $\Bbb{Z}$-basis
of the extended lattice $\widetilde{\Lambda _{G}^{\text{{\small ext}}}}$%
\textit{. }Now since 
\[
\frak{v}_{G}-\widetilde{\frak{u}}_{\rho +1}=\frac{1}{\frak{t}_{1}\left(
r-2\right) }\left( \alpha \cdot \frak{e}_{1}+\,l\cdot \frak{e}_{2}\right) 
\]
and 
\[
\frak{v}_{G}-\widetilde{\frak{u}}_{0}=\frac{1}{\frak{t}_{2}\left( r-2\right) 
}\left( \left( \alpha +\left( r-2\right) \right) \cdot \frak{e}_{1}+\,l\cdot 
\frak{e}_{2}\right) =\frac{1}{r-2}\left( \frak{z}_{2}\cdot \frak{e}_{1}+\,%
\frak{z}_{1}\cdot \frak{e}_{2}\right) , 
\]
using (\ref{DIOF}) we deduce: 
\[
\frak{v}_{G}-\widetilde{\frak{u}}_{\rho +1}=\breve{p}\cdot \left( \frak{v}%
_{G}-\widetilde{\frak{u}}_{0}\right) +q\cdot \left( -\frak{c}_{1}\cdot \frak{%
e}_{1}+\,\frak{c}_{2}\cdot \frak{e}_{2}\right) \ . 
\]
Hence, it suffices to take $p$ and $\frak{y}$ to be given by the formulae (%
\ref{PQTYPOI}) and (\ref{NIKA}), respectively. $_{\Box }$

\begin{remark}
\emph{Since }$p$\emph{\ is equal to the non negative remainder }$\left[ \,\breve{%
p}\,\right] _{q}$ \emph{modulo }$q$\emph{, it does not depend on the
particular choice of solutions} $\frak{c}_{1}\in \Bbb{Z}_{<0}$\emph{, }$%
\frak{c}_{2}\in \Bbb{N}$ \emph{of the linear diophantine equation (\ref{DIOF}%
). Nevertheless, in view of what was discussed in remark \ref{dioph}, we
give in the formulation of thm. \ref{MT2} a specific pair }$\left\{ \frak{c}%
_{1},\frak{c}_{2}\right\} $ \emph{for one of the most convenient solutions
of \ (\ref{DIOF}) which can be read off directly from the regular continued
fraction expansion of }$\frak{z}_{1}/\frak{z}_{2}$\emph{.\smallskip }
\end{remark}

\noindent Next lemma shows that our special assumption in lemma \ref{TYPOS2}
is included as a part of a necessary condition for the existence of $%
T_{N_{G}}$-equivariant, crepant, full resolutions.

\begin{lemma}
\label{ONED}Let $\left( X\left( N_{G},\Delta _{G}\right) ,\emph{orb}\left(
\sigma _{0}\right) \right) $ be a Gorenstein cyclic quotient singularity of
type\smallskip\ \emph{(\ref{typab})}$.$ If this singularity admits a $%
T_{N_{G}}$-equivariant, crepant, full resolution, then 
\[
\text{\emph{gcd}}\left( \alpha ,\beta ,l\right) \in \left\{ 1,r-2\right\} \,.
\]
Moreover, in the case in which \emph{gcd}$\left( \alpha ,\beta ,l\right) =1$%
, we have $\left[ \frak{t}_{1}\right] _{\left( r-2\right) }=\left[ \frak{t}%
_{2}\right] _{\left( r-2\right) }=1.$
\end{lemma}

\noindent \textit{Proof. }By definition, $1\leq $ gcd$\left( \alpha ,\beta
,l\right) \leq r-2$. Suppose $r\geq 5$ and gcd$\left( \alpha ,\beta
,l\right) \in \left\{ 2,3,\ldots ,r-3\right\} $. Obviously, 
\[
\frac{1}{l}\left( \stackunder{\left( r-2\right) \text{-times}}{\underbrace{%
\left[ \frac{l}{\text{gcd}\left( \alpha ,\beta ,l\right) }\right]
_{l},\ldots ,\left[ \frac{l}{\text{gcd}\left( \alpha ,\beta ,l\right) }%
\right] _{l}}},\left[ \frac{l}{\text{gcd}\left( \alpha ,\beta ,l\right) }%
\cdot \alpha \right] _{l},\left[ \frac{l}{\text{gcd}\left( \alpha ,\beta
,l\right) }\cdot \beta \right] _{l}\right) ^{\intercal }
\]
equals 
\[
\frac{1}{l}\left( \stackunder{\left( r-2\right) \text{-times}}{\underbrace{%
\frac{l}{\text{gcd}\left( \alpha ,\beta ,l\right) },\ldots ,\frac{l}{\text{%
gcd}\left( \alpha ,\beta ,l\right) }}},0,0\right) ^{\intercal }
\]
and hence it is a lattice point belonging to $\mathbf{Hlb}_{N_{G}}\left(
\sigma _{0}\right) $, because it cannot be written as the sum of two other
elements of $N_{G}\Bbb{r}\left\{ \mathbf{0}\right\} $ (cf. (\ref{Hilbbasis}%
)). On the other hand, 
\[
\frac{1}{l}\stackunder{\left( r-2\right) \text{-times}}{\underbrace{\left( 
\frac{l}{\text{gcd}\left( \alpha ,\beta ,l\right) }+\cdots +\frac{l}{\text{%
gcd}\left( \alpha ,\beta ,l\right) }\right) }}\ =\frac{r-2}{\text{gcd}\left(
\alpha ,\beta ,l\right) }>1,
\]
contradicting thm. \ref{KILLER}. Hence, gcd$\left( \alpha ,\beta ,l\right)
\in \left\{ 1,r-2\right\} .$ Now if gcd$\left( \alpha ,\beta ,l\right) =1$,
using lemma \ref{UNIM3} and the fact that the triangles conv$\left( \left\{ 
\frak{v}_{G},\widetilde{\frak{u}}_{0},\widetilde{\frak{u}}_{1}\right\}
\right) $ and conv$\left( \left\{ \frak{v}_{G},\widetilde{\frak{u}}_{\rho },%
\widetilde{\frak{u}}_{\rho +1}\right\} \right) $ have to be basic w.r.t. $%
\widetilde{\Lambda _{G}^{\text{{\small ext}}}}$ (with $\widetilde{\frak{u}}%
_{0},\widetilde{\frak{u}}_{1},\widetilde{\frak{u}}_{\rho },\widetilde{\frak{u%
}}_{\rho +1}\in \widetilde{\Lambda _{G}}$), we obtain 
\[
\left| \text{det}\left( \frak{v}_{G}-\widetilde{\frak{u}}_{0},\widetilde{%
\frak{u}}_{1}-\widetilde{\frak{u}}_{0}\right) \right| =\left| \text{det}%
\left( \frak{v}_{.G}-\widetilde{\frak{u}}_{\rho },\widetilde{\frak{u}}_{\rho
+1}-\widetilde{\frak{u}}_{\rho }\right) \right| =\text{det}\left( \widetilde{%
\Lambda _{G}^{\text{{\small ext}}}}\right) =\frac{1}{r-2}
\]
if and only if  
\[
\frac{\left[ \frak{t}_{2}\right] _{\left( r-2\right) }}{\frak{t}_{2}}%
\,\left| \text{det}\left( \frak{v}_{G},\widetilde{\frak{u}}_{1}-\widetilde{%
\frak{u}}_{0}\right) \right| =\frac{\left[ \frak{t}_{1}\right] _{\left(
r-2\right) }}{\frak{t}_{1}}\,\left| \text{det}\left( \frak{v}_{G},\widetilde{%
\frak{u}}_{\rho +1}-\widetilde{\frak{u}}_{\rho }\right) \right| =\frac{1}{r-2%
}\ ,
\]
which is equivalent to the existence of integers $\gamma _{1},\gamma
_{2},\gamma _{3},\gamma _{4}$, such that 
\[
\frac{\left[ \frak{t}_{2}\right] _{\left( r-2\right) }}{\frak{t}_{2}}\
\left| \gamma _{1}l+\gamma _{2}\left( \alpha +\left( r-2\right) \right)
\right| =\frac{\left[ \frak{t}_{1}\right] _{\left( r-2\right) }}{\frak{t}_{1}%
}\ \left| \gamma _{3}l+\gamma _{4}\left( \alpha +\left( r-2\right) \right)
\right| =1.
\]
Since necessarily 
\[
\gamma _{1}l+\gamma _{2}\left( \alpha +\left( r-2\right) \right) \equiv 0%
\text{ }\left( \text{mod }\frak{t}_{2}\right) ,\ \ \ \ \gamma _{3}l+\gamma
_{4}\left( \alpha +\left( r-2\right) \right) \equiv 0\text{ }\left( \text{%
mod }\frak{t}_{1}\right) ,
\]
we get $\left[ \frak{t}_{1}\right] _{\left( r-2\right) }=\left[ \frak{t}_{2}%
\right] _{\left( r-2\right) }=1$. $_{\Box }\bigskip $

\noindent $\bullet $\textbf{\ Fifth step}. Now we revert again to the
general setting and make essential use of the statements of \S \ref{BRUCH}
for the $\left( p,q\right) $-cone $\mathbf{\tau }_{G}$. If $\mathbf{\tau }%
_{G}$ is \textit{non-basic} (i.e., if $p\neq 0$), then we expand  $q/p$ as
regular continued fraction
\[
\dfrac{q}{p}=\left[ \lambda _{1},\lambda _{2},\ldots ,\lambda _{\kappa
-1},\lambda _{\kappa }\right] \ ,
\]
($\kappa \geq 2$, $\lambda _{\kappa }\geq 2$). As we explained in theorem 
\ref{KAP}, Kleinian approximations enable us to control completely the
vertices of both cones $\mathbf{\tau }_{G}$ and $\Bbb{R}_{\geq 0}\,\frak{y+}%
\Bbb{R}_{\geq 0}\,\left( \frak{v}_{G}-\widetilde{\frak{u}}_{\rho +1}\right) $
(having the ray $\Bbb{R}_{\geq 0}\,\left( \frak{v}_{G}-\widetilde{\frak{u}}%
_{\rho +1}\right) $ as common face, cf. rem. \ref{AAP}) by introducing the
recurrence relations 
\begin{equation}
\mathbf{v}_{0}:=\frak{v}_{G}-\widetilde{\frak{u}}_{0},\ \mathbf{v}_{1}:=%
\frak{y},\ \ \ \ \ \,\mathbf{v}_{i}:=\lambda _{i-1}\ \mathbf{v}_{i-1}+%
\mathbf{v}_{i-2},\ \forall i,\ \ 2\leq i\leq \kappa +1.  \label{REKUR}
\end{equation}

\begin{lemma}
\label{TELIKO}If $\mu _{G}\neq 1$ and $\mathbf{\tau }_{G}$ non-basic, then
each of the conditions of lemma \emph{\ref{UNIM3}} is equivalent to each of
the following conditions\emph{:\medskip }\newline
\emph{(i) }$\frak{y}\in \widetilde{\Lambda _{G}}$ \ and all points $\left\{ 
\mathbf{v}_{2i+1}\ \left| \ 1\leq i\leq \frac{\kappa -1}{2}\right. \right\} $%
, for $\kappa $ odd $\geq 3$, 
\[
\emph{(}\text{resp. all points \ }\left\{ \mathbf{v}_{2i+1}\ \left| \ 1\leq
i\leq \tfrac{\kappa }{2}-1,\,\kappa \neq 2\right. \right\} \cup \left\{ 
\mathbf{v}_{\kappa +1}-\mathbf{v}_{\kappa }\right\} \text{,\ for }\kappa 
\text{ even }\geq 2\emph{)}
\]
\emph{\ }of the extended lattice $\widetilde{\Lambda _{G}^{\text{\emph{ext}}}%
}$ belong already to $\widetilde{\Lambda _{G}}$.\medskip \newline
\emph{(ii) }$\frak{y}\in \widetilde{\Lambda _{G}}$ \ and all points $\left\{
\lambda _{2i}\left( \frak{v}_{G}-\widetilde{\frak{u}}_{0}\right) \ \left| \
1\leq i\leq \frac{\kappa -1}{2}\right. \right\} $, for $\kappa $ odd $\geq 3$%
, 
\[
\emph{(}\text{resp. all points \ }\left\{ \lambda _{2i}\left( \frak{v}_{G}-%
\widetilde{\frak{u}}_{0}\right) \ \left| \ 1\leq i\leq \tfrac{\kappa }{2}%
-1,\,\kappa \neq 2\right. \right\} \cup \left\{ \left( \lambda _{\kappa
}-1\right) \,\left( \frak{v}_{G}-\widetilde{\frak{u}}_{0}\right) \right\} 
\text{, for }\kappa \text{ \ even }\geq 2\emph{)}
\]
of $\ \widetilde{\Lambda _{G}^{\text{\emph{ext}}}}$ \ belong already to $%
\widetilde{\Lambda _{G}}$.\newline
\end{lemma}

\noindent \textit{Proof. }Since $\frak{l}_{1}=\left( \frak{v}_{G}-\widetilde{%
\frak{u}}_{0}\right) +\frak{y}$, we have 
\[
\frak{v}_{G}-\frak{l}_{1}\in \ \widetilde{\Lambda _{G}}\Longleftrightarrow 
\frak{y}=\mathbf{v}_{1}\in \ \widetilde{\Lambda _{G}}\Longleftrightarrow 
\text{conv}\left( \mathbf{v}_{0},\mathbf{v}_{2}\right) \cap \widetilde{%
\Lambda _{G}^{\text{{\small ext}}}}\subset \widetilde{\Lambda _{G}},
\]
where conv$\left( \mathbf{v}_{0},\mathbf{v}_{2}\right) $ is the ``first''
edge of $\partial \Theta _{\mathbf{\tau }_{G}}^{\mathbf{cp}}$\textit{\ }%
(containing only lattice points which are multiples of $\mathbf{v}_{1}$).
Applying the analogous argument also for the next coming edges of $\partial
\Theta _{\mathbf{\tau }_{G}}^{\mathbf{cp}}$, we obtain 
\[
\left[ 
\begin{array}{c}
\frak{v}_{G}-\frak{l}_{j}\in \ \widetilde{\Lambda _{G}},\text{ } \\ 
\forall \emph{\ }j\emph{,}\ 2\leq j\leq \rho +\newline
1
\end{array}
\right] \Longleftrightarrow \left\{ 
\begin{array}{ll}
\left\{ \mathbf{v}_{2i+1}\ \left| \ 1\leq i\leq \frac{\kappa -1}{2}\right.
\right\} \in \widetilde{\Lambda _{G}}, & \text{for \ }\kappa \notin 2\Bbb{Z}
\\ 
&  \\ 
\left\{ \mathbf{v}_{2i+1}\ \left| \ 1\leq i\leq \tfrac{\kappa }{2}%
-1,\,\kappa \neq 2\right. \right\} \cup \left\{ \mathbf{v}_{\kappa +1}-%
\mathbf{v}_{\kappa }\right\} \in \widetilde{\Lambda _{G}}, & \text{for \ }%
\kappa \in 2\Bbb{Z}
\end{array}
\right. 
\]
Hence, condition (i) is equivalent to condition (iii) of lemma \ref{UNIM3}.
(Note that these $\mathbf{v}_{2i+1}$'s cover the vertex-set of the interior
of the compact part of the support polytope which approximates the ray $\Bbb{%
R}_{\geq 0}\,\left( \frak{v}_{G}-\widetilde{\frak{u}}_{\rho +1}\right) $
``from above'' and is determined by the cone $\Bbb{R}_{\geq 0}\,\frak{y+}%
\Bbb{R}_{\geq 0}\,\left( \frak{v}_{G}-\widetilde{\frak{u}}_{\rho +1}\right) $%
).\medskip \newline
(i)$\Leftrightarrow $(ii) For $\kappa $ odd $\geq 3$ or $\kappa $ even $\geq
4$, and $i\in \left\{ 1,\ldots ,\QOVERD\lfloor \rfloor {\kappa
-1}{2}\right\} $, one easily shows via (\ref{REKUR}) that 
\[
\mathbf{v}_{2i+1}=\left( \sum_{j=0}^{i-1}\xi _{j}\cdot \mathbf{v}%
_{2j+1}\right) +\lambda _{2i}\cdot \mathbf{v}_{0}
\]
for suitable positive integers $\xi _{0},\xi _{1},\ldots ,\xi _{i-1}$.
Analogously for $\kappa $ even $\geq 2$ it is always possible to express $%
\mathbf{v}_{\kappa +1}-\mathbf{v}_{\kappa }$ as a positive integer linear
combination 
\[
\mathbf{v}_{\kappa +1}-\mathbf{v}_{\kappa }=\left( \sum_{j=0}^{\tfrac{\kappa 
}{2}-1}\xi _{j}\cdot \mathbf{v}_{2j+1}\right) +\left( \lambda _{\kappa
}-1\right) \cdot \mathbf{v}_{0}
\]
of the $\mathbf{v}_{2j+1}$'s and $\mathbf{v}_{0}$. This means that it is
enough for examining if $\mathbf{v}_{2i+1}$'s (resp. $\mathbf{v}_{\kappa +1}-%
\mathbf{v}_{\kappa }$) belong to $\widetilde{\Lambda _{G}}$ or not to
restrict ourselves to the consideration of $\lambda _{2i}\cdot \mathbf{v}_{0}
$'s (resp. $\left( \lambda _{\kappa }-1\right) \cdot \mathbf{v}_{0}$). $%
_{\Box }\bigskip \smallskip $

\noindent $\bullet $ \textbf{Proof of thm. \ref{MT2}: }Let $\left( X\left(
N_{G},\Delta _{G}\right) ,\text{orb}\left( \sigma _{0}\right) \right) $
denote a Gorenstein cyclic quotient singularity of type\smallskip\ (\ref
{typab})$.\medskip \smallskip $\textbf{\ }\newline
($\Rightarrow $): If this singularity admits $T_{N_{G}}$-equivariant,
crepant, full resolutions, then at least one of them is necessarily
projective (by cor. \ref{CORPROJ}), and by lemma \ref{ONED} there are only
two possibilities: \textbf{either} $\mu _{G}=1$, i.e., gcd$\left( \alpha
,\beta ,l\right) =r-2$ \textbf{or }gcd$\left( \alpha ,\beta ,l\right) =1$
and $\left[ \frak{t}_{1}\right] _{\left( r-2\right) }=\left[ \frak{t}_{2}%
\right] _{\left( r-2\right) }=1$. In the latter case, the $2$-dimensional
cone $\mathbf{\tau }_{G}\subset \left( \widetilde{\Lambda _{G}^{\text{%
{\small ext}}}}\right) _{\Bbb{R}}$ is a $\left( p,q\right) $-cone with 
\[
q=\text{ }\frac{l}{\frak{t}_{1}\cdot \frak{t}_{2}}\,,\ \ \ p=\left[ \,\breve{%
p}\,\right] _{q}\text{ \ \ and \ \ }\frak{y=}\left( \frac{\,\breve{p}-p\,\,}{%
q}\right) \,\left( \frak{v}_{G}-\widetilde{\frak{u}}_{0}\right) +\left( -%
\frak{c}_{1}\cdot \frak{e}_{1}+\,\frak{c}_{2}\cdot \frak{e}_{2}\right) 
\]
(by lemma \ref{TYPOS2}). If $\mathbf{\tau }_{G}$ is non-basic, then using
the one direction of lemma \ref{TELIKO}, we have: $\frak{y}\in \widetilde{%
\Lambda _{G}}$ \ and all points 
\[
\left\{ \lambda _{2i}\left( \frak{v}_{G}-\widetilde{\frak{u}}_{0}\right) \
\left| \ 1\leq i\leq \frac{\kappa -1}{2}\right. \right\} \text{, \ for \ }%
\kappa \ \text{odd\ }\geq 3,
\]
\[
\text{(resp. all points \ }\left\{ \lambda _{2i}\left( \frak{v}_{G}-%
\widetilde{\frak{u}}_{0}\right) \ \left| \ 1\leq i\leq \frac{\kappa }{2}%
-1,\,\kappa \neq 2\right. \right\} \cup \left\{ \left( \lambda _{\kappa
}-1\right) \,\left( \frak{v}_{G}-\widetilde{\frak{u}}_{0}\right) \right\} 
\text{, for }\kappa \text{ \ even }\geq 2\text{)}
\]
of $\ \widetilde{\Lambda _{G}^{\text{ext}}}$ \ belong already to $\widetilde{%
\Lambda _{G}}$. Since $\left( -\frak{c}_{1}\cdot \frak{e}_{1}+\,\frak{c}%
_{2}\cdot \frak{e}_{2}\right) \in \widetilde{\Lambda _{G}}$, condition $%
\frak{y}\in \widetilde{\Lambda _{G}}$ is equivalent to 
\[
\left( \frac{\,\breve{p}-p\,}{q}\right) \,\left( \frak{v}_{G}-\widetilde{%
\frak{u}}_{0}\right) =\frac{\,\breve{p}-p\,\,}{q\left( r-2\right) }\cdot 
\stackunder{\text{primitive lattice point}}{\underbrace{\left( \frak{z}%
_{2}\cdot \frak{e}_{1}+\,\frak{z}_{1}\cdot \frak{e}_{2}\right) }}\in 
\widetilde{\Lambda _{G}}\Longleftrightarrow \frac{\,\breve{p}-p\,\,}{q}%
\equiv 0\left( \text{mod }\left( r-2\right) \right) .
\]
On the other hand, for $\kappa $ odd $\geq 3$ or $\kappa $ even $\geq 4$,
and $i\in \left\{ 1,\ldots ,\QOVERD\lfloor \rfloor {\kappa -1}{2}\right\} $,
we have 
\[
\lambda _{2i}\left( \frak{v}_{G}-\widetilde{\frak{u}}_{0}\right) =\frac{%
\lambda _{2i}}{r-2}\cdot \left( \frak{z}_{2}\cdot \frak{e}_{1}+\,\frak{z}%
_{1}\cdot \frak{e}_{2}\right) \in \widetilde{\Lambda _{G}}%
\Longleftrightarrow \lambda _{2i}\equiv 0\left( \text{mod }\left( r-2\right)
\right) ,
\]
while for $\kappa $ even $\geq 2$, 
\[
\left( \lambda _{\kappa }-1\right) \,\left( \frak{v}_{G}-\widetilde{\frak{u}}%
_{0}\right) =\frac{\left( \lambda _{\kappa }-1\right) }{r-2}\cdot \left( 
\frak{z}_{2}\cdot \frak{e}_{1}+\,\frak{z}_{1}\cdot \frak{e}_{2}\right) \in 
\widetilde{\Lambda _{G}}\Longleftrightarrow \lambda _{\kappa }\equiv 1\left( 
\text{mod }\left( r-2\right) \right) .
\]
Hence, all conditions (\ref{CON2}) are satisfied, as we have
assserted.\medskip \newline
($\Leftarrow $): Converserly, if \textit{either} gcd$\left( \alpha ,\beta
,l\right) =1$ and $\mathbf{\tau }_{G}$ is basic, \textit{or} gcd$\left(
\alpha ,\beta ,l\right) =r-2$, then there is nothing to be said
(cf. lemma \ref{UNIM1},  remark 
\ref{MI1}). Furthermore, in the case in which gcd$\left( \alpha ,\beta
,l\right) =1$, $\mathbf{\tau }_{G}$ is non-basic, and conditions (\ref{CON2}%
) are fulfilled, we show again by the above arguments that each of the
conditions of lemma \ref{TELIKO} is true too, and then make use of the
``backtracking-method'' for the reverse logical implications of the
conditions of our previous lemmas: 
\[
\begin{array}{ccccccc}
\text{lemma \ref{TELIKO}} &  &  &  &  &  & \text{proposition \ref{CENTRAL}}
\\ 
\Downarrow  &  &  &  &  &  & \Uparrow  \\ 
\text{lemma \ref{UNIM3}} & \Longrightarrow  & \text{lemma \ref{UNIM2}} & 
\Longrightarrow  & \text{lemma \ref{UNIM1}} & \Longrightarrow  & \text{lemma 
\ref{hilmil}}
\end{array}
\]
Thus, also in this case the existence of $T_{N_{G}}$-equivariant, crepant,
full resolutions (and of at least one projective, cf. \ref{CORPROJ}) for the
singularity $\left( X\left( N_{G},\Delta _{G}\right) ,\text{orb}\left(
\sigma _{0}\right) \right) $ is indeed ensured. $_{\Box }$

\begin{remark}
\label{REDUCTION}\emph{Note that considering} $2$\emph{-parameter
singularities of (the most general) type (\ref{gentype}), the whole
``reduction-procedure'' presented in the second part of this section can be
again applied with minor modifications. First of all one has to determine an
appropriate }$\Bbb{Z}$\emph{-basis of the lattice} $\overline{N_{G}}$\emph{;
namely the ``nice'' elements (\ref{EIDIKO}) generating }$\overline{N_{G}}$  
\emph{in the case of type (\ref{typab}) have to be changed. This is always
possible by using Hermitian normal forms. On the other hand, the
corresponding lattice transformations} $\Phi $, $\Psi $ \emph{and} $\Upsilon 
$ \emph{are also definenable, though they become a little bit more
complicated. Finally, the }$\left( p,q\right) $\emph{-cone} $\mathbf{\tau }%
_{G}$ \emph{provides again the extra (somewhat ``wilder'') arithmetical} 
\emph{conditions for the existence of} $T_{N_{G}}$\emph{-equivariant,
crepant, full resolutions. This is why it does not seem to be of conceptual
theoretical interest to deal directly with (\ref{gentype}), as this work can
be done by a simple computer-program. The main and more interesting
conclusion in this context is that also in this most general situation, it
is enough to perform a }\textit{polynomial-time-algorithm}\emph{\ in order
to determine all the above mentioned extra arithmetical conditions, which,
in other words, is an obvious strengthening of theorem \ref{MT1} from the
purely computational point of view.}
\end{remark}

\section{Computing cohomology group dimensions\label{COHOM}}

\noindent To compute the cohomology group dimensions of the overlying spaces
of the crepant resolutions of our $2$-parameter series of Gorenstein cyclic
quotient singularities we need some concepts from enumerative combinatorics
(see e.g. Stanley \cite{Stanley2}, \S 4.6).\medskip \newline
$\bullet $ For a lattice $d$-dimensional polytope $P\subset N_{\Bbb{R}}$,
i.e., for $P$ an ``integral'' polytope w.r.t. a lattice $N\cong \Bbb{Z}%
^{d^{\prime }}$ of rank $d^{\prime }\geq d$, and $\kappa $ a positive
integer, let

\[
\mathbf{Ehr}_{N}\left( P,\nu \right) =\mathbf{a}_{0}\left( P\right) +\mathbf{%
a}_{1}\left( P\right) \ \nu +\cdots +\mathbf{a}_{d-1}\left( P\right) \ \nu
^{d-1}+\mathbf{a}_{d}\left( P\right) \ \nu ^{d}\ \in \ \Bbb{Q}\left[ \nu %
\right] 
\]
denote the \textit{Ehrhart polynomial }of $P$ (w.r.t. $N$), where 
\[
\mathbf{Ehr}_{N}\left( P,\nu \right) :=\#\,\left( \nu \,P\cap \left\{ 
\begin{array}{c}
\text{ the sublattice of\thinspace\ }N\text{ \thinspace of rank \thinspace }d
\\ 
\text{spanned \ by aff}\left( P\right) \cap N
\end{array}
\right\} \right) 
\]
and 
\[
\mathbf{Ehr}_{N}\left( P;\frak{x}\right) :=1+\sum_{\nu =1}^{\infty }\ 
\mathbf{Ehr}_{N}\left( P,\nu \right) \ \frak{x}^{\nu }\in \Bbb{Q}_{\,}\left[
\,\!\left[ \frak{x}\right] \!\,\right] 
\]
\textit{\ }the corresponding \textit{Ehrhart series}. Writing $\mathbf{Ehr}%
_{N}\left( P;\frak{x}\right) $ as 
\[
\mathbf{Ehr}_{N}\left( P;\frak{x}\right) =\frac{\ \mathbf{\delta }_{0}\left(
P\right) +\ \mathbf{\delta }_{1}\left( P\right) \ \frak{x}+\cdots +\ \mathbf{%
\delta }_{d-1}\left( P\right) \frak{x}^{d-1}+\ \mathbf{\delta }_{d}\left(
P\right) \ \frak{x}^{d}}{\left( 1-\frak{x}\right) ^{d+1}} 
\]
we obtain the so-called $\mathbf{\delta }$-\textit{vector} $\mathbf{\delta }%
\left( P\right) =\left( \mathbf{\delta }_{0}\left( P\right) ,\mathbf{\delta }%
_{1}\left( P\right) ,\ldots ,\mathbf{\delta }_{d-1}\left( P\right) ,\mathbf{%
\delta }_{d}\left( P\right) \right) $ of $P.$

\begin{definition}
\emph{For any} $d\in \Bbb{Z}_{\geq 0}$ \emph{we introduce the} \textit{%
transfer} $\mathbf{a}$-$\mathbf{\delta }$-\textit{matrix }$\mathcal{M}%
_{d}\in $ \emph{GL}$\left( d+1,\Bbb{Q}\right) $ \emph{(depending only on }$d$%
\emph{) to be defined as} 
\[
\mathcal{M}_{d}:=\left( \frak{R}_{i,j}\right) _{0\leq i,j\leq d}\ \ \ \text{%
\emph{\ with}}\ \ \ \ \frak{R}_{i,j}:=\dfrac{1}{d!}\ \left\{
\dsum\limits_{\xi =i}^{d}\ \QDATOPD[ ] {d}{\xi }\ \dbinom{\xi }{i}\ \left(
d-j\right) ^{\xi -i}\right\} 
\]
\emph{where} $\QDATOPD[ ] {d}{\xi }$ \emph{denotes the Stirling number (of
the first kind) of} $d$ \emph{over} $\xi $.
\end{definition}

\noindent The following lemma can be proved easily.

\begin{lemma}
For a lattice $d$-polytope $P\subset N_{\Bbb{R}}$ w.r.t. an $N\cong \Bbb{Z}%
^{d^{\prime }}$, $d^{\prime }\geq d$, we have 
\begin{equation}
\left( \mathbf{\delta }_{0}\left( P\right) ,\mathbf{\delta }_{1}\left(
P\right) ,\ldots ,\mathbf{\delta }_{d-1}\left( P\right) ,\mathbf{\delta }%
_{d}\left( P\right) \right) =\left( \mathbf{a}_{0}\left( P\right) ,\mathbf{a}%
_{1}\left( P\right) ,\ldots ,\mathbf{a}_{d-1}\left( P\right) ,\mathbf{a}%
_{d}\left( P\right) \right) \cdot \left( \left( \mathcal{M}_{d}\right)
^{\intercal }\right) ^{-1}  \label{DELTA}
\end{equation}
\end{lemma}

\begin{definition}
\emph{The} $\mathbf{f}$\emph{-vector of} $\mathbf{f}\left( \mathcal{S}%
\right) =\left( \mathbf{f}_{0}\left( \mathcal{S}\right) ,\mathbf{f}%
_{1}\left( \mathcal{S}\right) ,\ldots ,\mathbf{f}_{d-1}\left( \mathcal{S}%
\right) ,\mathbf{f}_{d}\left( \mathcal{S}\right) \right) $ \emph{of a pure }$%
d$\emph{-dimensional simplicial complex }$\mathcal{S}$ \emph{is defined by} 
\[
\,\mathbf{f}_{j}\left( \mathcal{S}\right) :=\#\left\{ j\text{\emph{%
-dimensional simplices of }}\mathcal{S}\right\} \ .
\]
\end{definition}

\begin{lemma}
\label{STAL}If $P\subset N_{\Bbb{R}}$ is a lattice $d$-polytope w.r.t. an $%
N\cong \Bbb{Z}^{d^{\prime }}$, $d^{\prime }\geq d$, admitting a basic
triangulation $\mathcal{T}$,\ \ then 
\[
\mathbf{Ehr}_{N}\left( P,\nu \right) =\sum_{j=0}^{d}\binom{\nu -1}{j}\,\,%
\mathbf{f}_{j}\left( \mathcal{T}\right) 
\]
\end{lemma}

\noindent \textit{Proof. }See Stanley \cite{Stanley1}, cor. 2.5, p. 338. $%
_{\Box }$

\begin{theorem}[Cohomology group dimensions]
\ \ \smallskip \newline
Let $\left( X\left( N_{G},\Delta _{G}\right) ,\emph{orb}\left( \sigma
_{0}\right) \right) $ be a Gorenstein abelian quotient msc-singularity. If $%
X\left( N_{G},\Delta _{G}\right) $ admits crepant, $T_{N_{G}}$-equivariant,
full desingularizations, then only the \textbf{even}-dimensional cohomology
groups of the desingularized spaces $X\left( N_{G},\widehat{\Delta }%
_{G}\right) $ can be non-trivial, and 
\begin{equation}
\fbox{$
\begin{array}{c}
\end{array}
\text{\emph{dim}}_{\Bbb{Q}}H^{2i}\left( X\left( N_{G},\widehat{\Delta }%
_{G}\right) ;\Bbb{Q}\right) =\mathbf{\delta }_{i}\left( \frak{s}_{G}\right) 
\begin{array}{c}
\end{array}
$}  \label{BADA}
\end{equation}
$\forall i$, $0\leq i\leq r-1$, where $\mathbf{\delta }_{i}\left( \frak{s}%
_{G}\right) $'s are the components of \ the $\mathbf{\delta }$-\textit{vector%
} $\mathbf{\delta }\left( \frak{s}_{G}\right) $ of the $\left( r-1\right) $%
-dimensional junior simplex $\frak{s}_{G}$ \emph{(}w.r.t. $N_{G}$\emph{%
).\medskip }\newline
$\bullet $ By \emph{(\ref{BADA}) }and\emph{\ (\ref{DELTA}) }it is now clear
that these cohomology group dimensions \textbf{do not} depend on the basic
triangulations by means of which one constructs the subdivisions $\widehat{%
\Delta }_{G}$ of $\sigma _{0}$, but \textbf{only} on the coefficients of the
Ehrhart polynomial of $\frak{s}_{G}$.\medskip \newline
$\bullet $ In particular, if $\left( X\left( N_{G},\Delta _{G}\right) ,\emph{%
orb}\left( \sigma _{0}\right) \right) $ denotes a Gorenstein cyclic quotient
msc-singularity having type $\tfrac{1}{l}\left( \alpha _{1},\ldots ,\alpha
_{r}\right) $ with $l=\left| G\right| \geq r\geq 4$\emph{, } for which at
least $r-2$ of its defining weights are equal, and if $X\left( N_{G},\Delta
_{G}\right) $ admits crepant, $T_{N_{G}}$-equivariant, full
desingularizations, then the Ehrhart polynomial of $\frak{s}_{G}$ \emph{(}%
whose coefficients lead by \emph{(\ref{DELTA}) }to the computation of the
corresponding non-trivial cohomology group dimensions of $X\left( N_{G},%
\widehat{\Delta }_{G}\right) $'s\emph{)} can be determined as follows\emph{%
:\medskip } \newline
$\bullet $ At first we define the polynomial 
\[
\frak{B}_{G}\left( i;\nu \right) :=
\]
\begin{equation}
\sum_{j=0}^{i+2}\tbinom{\nu -1}{j}\left\{ \tbinom{i}{j+1}+\#\left( \overline{%
\frak{s}_{G}}\cap \overline{N_{G}}\right) \tbinom{i}{j}+\left( 3\,\text{%
\emph{Vol}}\left( \overline{\frak{s}_{G}}\right) +\tfrac{1}{2}\#\left(
\partial \left( \overline{\frak{s}_{G}}\right) \cap \overline{N_{G}}\right)
\right) \tbinom{i}{j-1}+2\,\text{\emph{Vol}}\left( \overline{\frak{s}_{G}}%
\right) \tbinom{i}{j-2}\right\}   \label{BPOL}
\end{equation}
for all $i$, $1\leq i\leq r-3$ \emph{(}where $\frak{Q}_{G}=$ \emph{conv}$%
\left( \frak{s}_{G}\cap \overline{N_{G}}\right) \subseteqq \overline{\frak{s}%
_{G}}=\frak{s}_{G}\cap \frak{L}$ $\subset \left( \overline{N_{G}}\right) _{%
\Bbb{R}}$ as in \emph{\S \ref{MAIN}\textsf{(a)}).} In addition, in the case
in which $\mu _{G}\neq 1$, we define the polynomial 
\begin{equation}
\frak{D}\left( i;\nu \right) :=\sum_{j=0}^{i+1}\tbinom{\nu -1}{j}\left\{ 
\tbinom{i}{j+1}+\left( \rho +2\right) \tbinom{i}{j}+\left( \rho +1\right) 
\tbinom{i}{j-1}\right\}   \label{DPOL}
\end{equation}
for all $i$, $1\leq i\leq r-2$ \emph{(}with $\rho $ as in \emph{\ref{ROS}).}
Using $\,\frak{B}_{G}\left( i;\nu \right) $ and $\,\frak{D}_{G}\left( i;\nu
\right) $ we obtain\emph{:\medskip } \newline
\emph{(i)} If $\mu _{G}=1$, 
\begin{equation}
\mathbf{Ehr}_{N_{G}}\left( \frak{s}_{G},\nu \right) =\sum_{i=1}^{r-3}\
\left( -1\right) ^{r-3-i}\ \binom{r-2}{i}\,\frak{B}_{G}\left( i;\nu \right)
+\left( -1\right) ^{r-3}\ \mathbf{Ehr}_{\overline{N_{G}}}\left( \frak{Q}%
_{G},\nu \right)   \label{EHR1}
\end{equation}
\emph{(ii)} If $\mu _{G}\neq 1$, 
\[
\mathbf{Ehr}_{N_{G}}\left( \frak{s}_{G},\nu \right) =\sum_{i=1}^{r-3}\
\left( -1\right) ^{r-3-i}\ \binom{r-2}{i}\,\frak{B}_{G}\left( i;\nu \right)
+\left( -1\right) ^{r-3}\ \mathbf{Ehr}_{\overline{N_{G}}}\left( \frak{Q}%
_{G},\nu \right) +
\]
\begin{equation}
+\sum_{i=1}^{r-2}\ \left( -1\right) ^{r-2-i}\ \binom{r-2}{i}\,\frak{D}%
_{G}\left( i;\nu \right) +\left( -1\right) ^{r-2}\ \left( \left( \rho
+1\right) \,\nu +1\right)   \label{EHR2}
\end{equation}
Moreover, in both cases we have 
\begin{equation}
\mathbf{Ehr}_{\overline{N_{G}}}\left( \frak{Q}_{G},\nu \right) =\#\left( \nu
\,\frak{Q}_{G}\cap \overline{N_{G}}\right) =\text{\emph{Vol}}\left( \frak{Q}%
_{G}\right) \,\nu ^{2}+\frac{1}{2}\,\#\left( \partial \frak{Q}_{G}\cap 
\overline{N_{G}}\right) \,\nu +1  \label{PICK}
\end{equation}
\end{theorem}

\noindent \textit{Proof. }For the proof of the first statement see
Batyrev-Dais \cite{BD}, thm. 4.4, p. 909. Now let 
\[
\left( X\left( N_{G},\Delta _{G}\right) ,\text{orb}\left( \sigma _{0}\right)
\right) 
\]
be a $2$-parameter Gorenstein cyclic quotient msc-singularity admitting
torus-equivariant crepant resolutions. Since the dimension of $H^{2i}\left(
X\left( N_{G},\widehat{\Delta }_{G}\right) ;\Bbb{Q}\right) $ does not depend
on the choice of the basic triangulations constructing $\widehat{\Delta }%
_{G} $'s, it is enough for its computation to consider just \textit{one}
triangulation of this sort. Hereafter take a $\widehat{\Delta }_{G}=\widehat{%
\Delta }_{G}\left( \mathcal{T}\left[ \frak{T}\right] \right) $ being induced
by a fixed basic triangulation of the special form $\mathcal{T}\left[ \frak{T%
}\right] $ (see \ref{JOI}, \ref{COBAS}). Formula (\ref{PICK}) is nothing but
Pick's theorem applied to the lattice polygon $\frak{Q}_{G}$. Next we shall
treat the two possible cases separately:\medskip \newline
(i) If $\mu _{G}=1$, then $\overline{\frak{s}_{G}}=\frak{Q}_{G}$ and $\frak{s%
}_{G}$ can be written as 
\[
\frak{s}_{G}=\dbigcup\limits_{\left( \xi _{1},\xi _{2},\ldots ,\xi
_{r-3}\right) \in \Xi _{r}}\ \overline{\frak{s}_{G}}\ast \text{conv}\left(
\left\{ e_{\xi _{1}},e_{\xi _{2}},\ldots ,e_{\xi _{r-3}}\right\} \right) 
\]
Using the principle of inclusion-exclusion we deduce 
\[
\mathbf{Ehr}_{N_{G}}\left( \frak{s}_{G},\nu \right) = 
\]
\begin{equation}
\sum_{i=1}^{r-3}\ \left( -1\right) ^{r-3-i}\ \sum_{1\leq \xi _{1}<\xi
_{2}<\cdots <\xi _{i}\leq r-2}\ \mathbf{Ehr}_{N_{G}}\left( \overline{\frak{s}%
_{G}}\ast \text{conv}\left( e_{\xi _{1}},\ldots ,e_{\xi _{i}}\right) ,\nu
\right) +\left( -1\right) ^{r-3}\ \mathbf{Ehr}_{\overline{N_{G}}}\left( 
\overline{\frak{s}_{G}},\nu \right)  \label{ER1}
\end{equation}
By lemma \ref{STAL} we obtain for the join $\overline{\frak{s}_{G}}\ast $conv%
$\left( e_{\xi _{1}},\ldots ,e_{\xi _{i}}\right) $ of dimension $2+\left(
i-1\right) +1=i+2,$ 
\begin{equation}
\mathbf{Ehr}_{N_{G}}\left( \overline{\frak{s}_{G}}\ast \text{conv}\left(
e_{\xi _{1}},\ldots ,e_{\xi _{i}}\right) ,\nu \right) =\sum_{j=0}^{i+2}%
\tbinom{\nu -1}{j}\ \mathbf{f}_{j}\left( \left. \mathcal{T}\left[ \frak{T}%
\right] \,\right| \,_{\overline{\frak{s}_{G}}\ast \text{conv}\left( e_{\xi
_{1}},\ldots ,e_{\xi _{i}}\right) }\right)  \label{ER2}
\end{equation}
with 
\[
\mathbf{f}_{j}\left( \left. \mathcal{T}\left[ \frak{T}\right] \,\right| \,_{%
\overline{\frak{s}_{G}}\ast \text{conv}\left( e_{\xi _{1}},\ldots ,e_{\xi
_{i}}\right) }\right) =\mathbf{f}_{j}\left( \frak{T}\ast \text{conv}\left(
\left\{ e_{\xi _{1}},\ldots ,e_{\xi _{i}}\right\} \right) \right) = 
\]
\[
=\mathbf{f}_{j}\left( \text{conv}\left( \left\{ e_{\xi _{1}},..,e_{\xi
_{i}}\right\} \right) \right) +\mathbf{f}_{0}\left( \frak{T}\right) \cdot 
\mathbf{f}_{j-1}\left( \text{conv}\left( \left\{ e_{\xi _{1}},..,e_{\xi
_{i}}\right\} \right) \right) +\smallskip 
\]
\[
+\mathbf{f}_{1}\left( \frak{T}\right) \cdot \mathbf{f}_{j-2}\left( \text{conv%
}\left( \left\{ e_{\xi _{1}},..,e_{\xi _{i}}\right\} \right) \right) +%
\mathbf{f}_{2}\left( \frak{T}\right) \cdot \mathbf{f}_{j-3}\left( \text{conv}%
\left( \left\{ e_{\xi _{1}},..,e_{\xi _{i}}\right\} \right) \right) , 
\]
where 
\begin{equation}
\mathbf{f}_{j}\left( \text{conv}\left( \left\{ e_{\xi _{1}},..,e_{\xi
_{i}}\right\} \right) \right) =\binom{i}{j+1}  \label{ER3}
\end{equation}
and 
\begin{equation}
\mathbf{f}_{0}\left( \frak{T}\right) =\#\left( \overline{\frak{s}_{G}}\cap 
\overline{N_{G}}\right) ,\ \ \mathbf{f}_{1}\left( \frak{T}\right) =3\,\text{%
Vol}\left( \overline{\frak{s}_{G}}\right) +\tfrac{1}{2}\#\left( \partial
\left( \overline{\frak{s}_{G}}\right) \cap \overline{N_{G}}\right) ,\ \ 
\mathbf{f}_{2}\left( \frak{T}\right) =2\,\text{Vol}\left( \overline{\frak{s}%
_{G}}\right) \ .  \label{ER4}
\end{equation}
Since (\ref{ER3}) is valid for \textit{all} $i$, $1\leq i\leq r-3$, and 
\textit{all} $\binom{r-2}{i}$ $i$-tuples $1\leq \xi _{1}<\cdots <\xi
_{i}\leq r-2$, the formulae (\ref{ER2}), (\ref{ER3}) and (\ref{ER4}) imply 
\begin{equation}
\mathbf{Ehr}_{N_{G}}\left( \overline{\frak{s}_{G}}\ast \text{conv}\left(
e_{\xi _{1}},\ldots ,e_{\xi _{i}}\right) ,\nu \right) =\frak{B}_{G}\left(
i;\nu \right)  \label{ER5}
\end{equation}
and (\ref{EHR1}) follows from (\ref{ER1}) and (\ref{ER5}).\medskip \newline
(ii) If $\mu _{G}\neq 1$, then 
\[
\frak{s}_{G}=\dbigcup\limits_{\left( \xi _{1},..,\xi _{r-3}\right) \in \Xi
_{r}}\ \left( \frak{Q}_{G}\ast \text{conv}\left( \left\{ e_{\xi
_{1}},..,e_{\xi _{r-3}}\right\} \right) \right) \cup \left[ \left(
\dbigcup\limits_{j=0}^{\rho }\text{conv}\left( \left\{ \frak{w}_{j},\frak{w}%
_{j+1}\right\} \right) \right) \ast \text{conv}\left( \left\{
e_{1},e_{2},..,e_{r-2}\right\} \right) \right] 
\]
and applying again inclusion-exclusion-principle we get 
\[
\mathbf{Ehr}_{N_{G}}\left( \frak{s}_{G},\nu \right) = 
\]
\[
=\sum_{i=1}^{r-3}\ \left( -1\right) ^{r-3-i}\ \sum_{1\leq \xi _{1}<\xi
_{2}<\cdots <\xi _{i}\leq r-2}\ \mathbf{Ehr}_{N_{G}}\left( \frak{Q}_{G}\ast 
\text{conv}\left( e_{\xi _{1}},\ldots ,e_{\xi _{i}}\right) ,\nu \right)
+\left( -1\right) ^{r-3}\ \mathbf{Ehr}_{\overline{N_{G}}}\left( \frak{Q}%
_{G},\nu \right) + 
\]
\[
+\,\mathbf{Ehr}_{N_{G}}\left( \left( \dbigcup\limits_{j=0}^{\rho }\text{conv}%
\left( \left\{ \frak{w}_{j},\frak{w}_{j+1}\right\} \right) \right) \ast 
\text{conv}\left( \left\{ e_{1},e_{2},..,e_{r-2}\right\} \right) ,\nu
\right) - 
\]
\begin{equation}
-\,\mathbf{Ehr}_{N_{G}}\left( \dbigcup\limits_{\left( \xi _{1},\xi
_{2},\ldots ,\xi _{r-3}\right) \in \Xi _{r}}\ \left( \left(
\dbigcup\limits_{j=0}^{\rho }\text{conv}\left( \left\{ \frak{w}_{j},\frak{w}%
_{j+1}\right\} \right) \right) \ast \text{conv}\left( \left\{ e_{\xi
_{1}},..,e_{\xi _{r-3}}\right\} \right) \right) ,\nu \right)  \label{ER6}
\end{equation}
By similar arguments to thosed used in (i) we conclude 
\begin{equation}
\mathbf{Ehr}_{N_{G}}\left( \left( \dbigcup\limits_{j=0}^{\rho }\text{conv}%
\left( \left\{ \frak{w}_{j},\frak{w}_{j+1}\right\} \right) \right) \ast 
\text{conv}\left( \left\{ e_{1},e_{2},..,e_{r-2}\right\} \right) ,\nu
\right) =\frak{D}_{G}\left( r-2;\nu \right)  \label{ER7}
\end{equation}
and 
\begin{equation}
\mathbf{Ehr}_{N_{G}}\left( \left( \left( \dbigcup\limits_{j=0}^{\rho }\text{%
conv}\left( \left\{ \frak{w}_{j},\frak{w}_{j+1}\right\} \right) \right) \ast 
\text{conv}\left( \left\{ e_{\xi _{1}},..,e_{\xi _{i}}\right\} \right)
\right) ,\nu \right) =\frak{D}_{G}\left( i;\nu \right)  \label{ER8}
\end{equation}
for \textit{all} $i$, $1\leq i\leq r-3$, and \textit{all} $\binom{r-2}{i}$ $%
i $-tuples $1\leq \xi _{1}<\cdots <\xi _{i}\leq r-2$. Furthermore, 
\begin{equation}
\mathbf{Ehr}_{N_{G}}\left( \left( \dbigcup\limits_{j=0}^{\rho }\text{conv}%
\left( \left\{ \frak{w}_{j},\frak{w}_{j+1}\right\} \right) \right) ,\nu
\right) =\left( \rho +1\right) \,\nu +1  \label{ER9}
\end{equation}
and formula (\ref{EHR2}) follows easily from (\ref{ER6}), (\ref{ER7}), (\ref
{ER8}) and (\ref{ER9}). $_{\Box }$

\begin{corollary}
Let $\left( X\left( N_{G},\Delta _{G}\right) ,\emph{orb}\left( \sigma
_{0}\right) \right) $ be a Gorenstein cyclic quotient singularity of type 
\emph{(\ref{typab}). }If one of the conditions \emph{(i), (ii)} of thm. 
\emph{\ref{MT2} }is fulfilled, then the dimensions of the non-trivial
cohomology groups of all spaces $X\left( N_{G},\widehat{\Delta }_{G}\right) $
desingularizing fully $X\left( N_{G},\Delta _{G}\right) $ by $T_{N_{G}}$%
-equivariant crepant morphisms are given by the formulae \emph{(\ref{BADA}),
(\ref{DELTA}), (\ref{EHR1}),(\ref{EHR2}) }and \emph{(\ref{PICK}) }which
depend \textbf{only }on $\alpha ,\beta ,r$ because\emph{:\medskip }\newline
\emph{(i)} If $\mu _{G}=1$, i.e., if \emph{gcd}$\left( \alpha ,\beta
,l\right) =r-2$, then 
\begin{equation}
\text{\emph{Vol}}\left( \frak{Q}_{G}\right) =\frac{l}{2\left( r-2\right) }=%
\frac{1}{2}+\frac{\alpha +\beta }{r-2}  \label{VOL1}
\end{equation}
and 
\begin{equation}
\#\left( \partial \frak{Q}_{G}\cap \overline{N_{G}}\right) =\text{\emph{gcd}}%
\left( \frac{\alpha }{r-2}+1,\frac{\alpha +\beta }{r-2}+1\right) +\text{%
\emph{gcd}}\left( \frac{\alpha }{r-2},\frac{\alpha +\beta }{r-2}+1\right) +1
\label{RAND1}
\end{equation}
\emph{(ii)} If $\mu _{G}\neq 1$, i.e., if \emph{gcd}$\left( \alpha ,\beta
,l\right) =1$, and $p\neq 0$, then using the continued fraction expansion 
\[
\dfrac{q}{p}=\left[ \lambda _{1},\lambda _{2},\ldots ,\lambda _{\kappa
-1},\lambda _{\kappa }\right] 
\]
\emph{(}as defined in thm. \emph{\ref{MT2}) }we obtain 
\begin{equation}
\text{\emph{Vol}}\left( \frak{Q}_{G}\right) =\left\{ 
\begin{array}{ll}
\frac{1}{2\left( r-2\right) }\left( l-\sum_{i=1}^{\frac{\kappa }{2}}\lambda
_{2i}-1\right) , & \,\text{\emph{for \ }}\kappa \text{ \ \emph{even}} \\ 
&  \\ 
\frac{1}{2\left( r-2\right) }\left( l-\sum_{i=1}^{\frac{\kappa -1}{2}%
}\lambda _{2i}\right) , & \text{\emph{for \ }}\kappa \ \text{\emph{odd}}
\end{array}
\right.   \label{VOL2}
\end{equation}
and 
\begin{equation}
\#\left( \partial \frak{Q}_{G}\cap \overline{N_{G}}\right) =\left\{ 
\begin{array}{cc}
\sum_{i=1}^{\frac{\kappa }{2}}\lambda _{2i}+\QOVERD\lfloor \rfloor {\text{%
\emph{gcd}}\left( \alpha ,l\right) }{r-2}+\QOVERD\lfloor \rfloor {\text{%
\emph{gcd}}\left( \beta ,l\right) }{r-2}+2, & \,\text{\emph{for \ }}\kappa 
\text{ \ \emph{even}} \\ 
&  \\ 
\sum_{i=1}^{\frac{\kappa -1}{2}}\lambda _{2i}+\QOVERD\lfloor \rfloor {\text{%
\emph{gcd}}\left( \alpha ,l\right) }{r-2}+\QOVERD\lfloor \rfloor {\text{%
\emph{gcd}}\left( \beta ,l\right) }{r-2}+3, & \text{\emph{for \ }}\kappa \ 
\text{\emph{odd}}
\end{array}
\right. \smallskip \smallskip   \label{RAND2}
\end{equation}
If $p$ happens to be $=0$, then instead of the formulae \emph{(\ref{VOL2})}
and \emph{(\ref{RAND2}), }use
\begin{equation}
\text{\emph{Vol}}\left( \frak{Q}_{G}\right) =\frac{1}{2\left( r-2\right) }%
\left( l-1\right)   \label{VOL3}
\end{equation}
and
\begin{equation}
\#\left( \partial \frak{Q}_{G}\cap \overline{N_{G}}\right) =\QOVERD\lfloor
\rfloor {\text{\emph{gcd}}\left( \alpha ,l\right) }{r-2}+\QOVERD\lfloor
\rfloor {\text{\emph{gcd}}\left( \beta ,l\right) }{r-2}+2,  \label{RAND3}
\end{equation}
respectively.
\end{corollary}

\noindent \textit{Proof. }The formula (\ref{VOL1}) in case (i) is obvious. (%
\ref{RAND1}) follows from 
\[
\#\left( \text{conv}\left( \frac{\frak{n}_{G}}{\mu_{G}},e_{r}\right) \cap \overline{N_{G}}%
\right) =\text{gcd}\left( \frac{\alpha }{r-2},\frac{l}{r-2}\right) +1,
\]
and
\[
\#\left( \text{conv}\left( \frac{\frak{n}_{G}}{\mu_{G}},e_{r-1}\right) \cap \overline{N_{G}}%
\right) =\text{gcd}\left( \frac{\alpha }{r-2}+1,\frac{l}{r-2}\right) +1.
\]
In case (ii) we obtain (\ref{VOL2}), (\ref{VOL3}) by 
\[
\text{Vol}\left( \frak{Q}_{G}\right) =\text{Vol}\left( \overline{\frak{s}_{G}%
}\right) \,\mathbb{-}\sum_{j=0}^{\rho }\text{Vol}\left( \text{conv}\left(
\left\{ \frac{\frak{n}_{G}}{\mu_{G}},\frak{w}_{j},\frak{w}_{j+1}\right\} \right) \right) =%
\frac{1}{2\left( r-2\right) }\left( l-\rho -1\right) \,.
\]
Finally, the number (\ref{RAND2}) (resp. (\ref{RAND3})) of the lattice
points lying on the boundary of $\frak{Q}_{G}$ equals 
\[
\#\left( \partial \frak{Q}_{G}\cap \overline{N_{G}}\right) =\rho +\#\left( 
\text{conv}\left( \frak{w}^{\prime },e_{r}\right) \cap \overline{N_{G}}%
\right) +\#\left( \text{conv}\left( \frak{w},e_{r-1}\right) \cap \overline{%
N_{G}}\right) ,
\]
where 
\[
\#\left( \text{conv}\left( \frak{w}^{\prime },e_{r}\right) \cap \overline{%
N_{G}}\right) =\QOVERD\lfloor \rfloor {\text{gcd}\left( \beta ,l\right)
}{r-2}+1,\ \ \#\left( \text{conv}\left( \frak{w},e_{r-1}\right) \cap 
\overline{N_{G}}\right) =\QOVERD\lfloor \rfloor {\text{gcd}\left( \alpha
,l\right) }{r-2}+1,
\]
and the proof is completed by expressing $\rho $ by the entries of the above
regular continued fraction (cf. rem. \ref{length2}). $_{\Box }\bigskip
\smallskip $

\noindent\textit{Acknowledgements.} The first author expresses his thanks to
G.~M.~Ziegler for a useful discussion about maximal triangulations. The
third author would like to acknowledge the support by the DFG Leibniz-Preis
of M.~Gr\"{o}tschel. The authors thank the Mathematics Department of
T\"{u}bingen University and Konrad-Zuse-Zentrum in Berlin, respectively, for
ideal working conditions during the writing of this paper.\medskip\

\vspace{2cm}

\hbox{
\vbox{\noindent Dimitrios I.~Dais \\
\noindent Mathematisches Institut \\ 
\noindent Universit{\"a}t T{\"u}bingen \\
\noindent Auf der Morgenstelle 10 \\
\noindent 72076 T{\"u}bingen, Germany \\
\noindent dais@wolga.mathematik.uni-tuebingen.de}
\hspace{-5cm}
\vbox{ \noindent Utz-Uwe Haus, Martin Henk \\
\noindent Konrad-Zuse-Zentrum Berlin \\
\noindent Takustra{\ss}e 7 \\
\noindent 14195 Berlin, Germany \\
\noindent haus\{henk\}@zib.de \\
\noindent \hphantom{space}}
}

\end{document}